\documentclass[reqno,11pt]{amsart}

\usepackage{amsthm, amsmath, amssymb, stmaryrd}
\usepackage[utf8]{inputenc}
\usepackage[T1]{fontenc}

\usepackage[hyphenbreaks]{breakurl}
\usepackage[hyphens]{url}
\usepackage{systeme}
\usepackage[shortlabels]{enumitem}
\usepackage[hidelinks]{hyperref}
\usepackage{microtype}

\usepackage{bm}
\usepackage[margin=1in]{geometry}

\usepackage[textsize=scriptsize,backgroundcolor=orange!5]{todonotes}

\usepackage[noabbrev,capitalize,sort]{cleveref}
\crefname{equation}{}{}
\numberwithin{equation}{section}

\usepackage{mathtools}
\newtheorem{theorem}{Theorem}[section]
\newtheorem{proposition}[theorem]{Proposition}
\newtheorem{lemma}[theorem]{Lemma}
\newtheorem{claim}[theorem]{Claim}
\newtheorem{corollary}[theorem]{Corollary}
\newtheorem{conjecture}[theorem]{Conjecture}

\theoremstyle{definition}
\newtheorem{definition}[theorem]{Definition}

\newtheorem{question}[theorem]{Question}
\newtheorem{example}[theorem]{Example}

\theoremstyle{remark}
\newtheorem*{remark}{Remark}

\newcommand{\abs}[1]{\left\lvert#1\right\rvert}

\newcommand{\avgbeta}{\beta}

\DeclareMathOperator{\mult}{mult}   
\DeclareMathOperator{\vol}{vol}    

\newcommand{\id}{\text{Id}}

\newcommand*{\eqdef}{\stackrel{\mbox{\normalfont\tiny def}}{=}}

\newcommand{\FF}{\mathbb{F}}
\newcommand{\RR}{\mathbb{R}}
\newcommand{\NN}{\mathbb{N}}
\newcommand{\ZZ}{\mathbb{Z}}
\newcommand*{\PP}{\mathbb{P}} 

\newcommand{\cL}{\mathcal L}
\newcommand{\cJ}{\mathcal J}
\newcommand{\cC}{\mathcal C}

\newcommand{\cV}{\mathcal V}
\newcommand{\cE}{\mathcal E}
\newcommand{\cD}{\mathcal D}
\newcommand{\cA}{\mathcal A}
\newcommand{\cM}{\mathcal M}
\newcommand{\cB}{\mathcal B}
\newcommand{\cS}{\mathcal S}

\newcommand{\bF}{\mathbf{F}}
\newcommand{\bJ}{\mathbf{J}}
\newcommand{\fm}{\mathfrak{m}}

\newcommand{\Hasse}{\mathsf H} 

\newcommand\tT{\vcenter{\hbox{\scalebox{0.6}{$T$}}}}

\newlength{\hght}

\newcommand{\halfscript}[2]{\settoheight{\hght}{a}{#1\!\!\:\:}\raisebox{.5\hght}{$\scriptstyle{#2}$}}

\makeatletter
\newcommand\thankssymb[1]{\textsuperscript{\@fnsymbol{#1}}}
\makeatother

\author[Ting-Wei Chao]{Ting-Wei Chao\thankssymb{1}}
\author[Hung-Hsun Hans Yu]{Hung-Hsun Hans Yu\thankssymb{2}}

\thanks{\thankssymb{1}Department of Mathematical Sciences, Carnegie Mellon University, Pittsburgh, PA 15213, USA\@. Supported in part
by U.S. taxpayers through NSF grant DMS-2154063 and NSF CAREER grant DMS-1555149. Email: {\tt tchao2@andrew.cmu.edu}}

\thanks{\thankssymb{2}Department of Mathematics, Princeton University, Princeton, NJ 08544\@.  Email: {\tt hansonyu@princeton.edu}}

\title{Tight Bound and Structural Theorem for Joints}

\begin{document}

\maketitle

\begin{abstract}
	A joint of a set of lines $\mathcal{L}$ in $\mathbb{F}^d$ is a point that is contained in $d$ lines with linearly independent directions. 
    The joints problem asks for the maximum number of joints that are formed by $L$ lines. 
    Guth and Katz showed that the number of joints is at most $O(L^{3/2})$ in $\mathbb{R}^3$ using polynomial method.
    This upper bound is met by the construction given by taking the joints and the lines to be all the $d$-wise intersections and all the $(d-1)$-wise intersections of $M$ hyperplanes in general position.
    Furthermore, this construction is conjectured to be optimal.
 
    In this paper, we verify the conjecture and show that this is the only optimal construction by using a more sophisticated polynomial method argument. 
    This is the first tight bound and structural theorem obtained using this method. 
    We also give a new definition of multiplicity that strengthens the main result of a previous work by Tidor, Zhao and the second author.
    Lastly, we relate the joints problem to some set-theoretic problems and prove conjectures of Bollob\'as and Eccles regarding partial shadows.
\end{abstract}

\section{Introduction}
\subsection{Joints of lines}
The joints problem asks to find the maximum number of joints formed by $L$ lines in $\RR^d$ (or $\FF^d$ where $\FF$ is an arbitrary field), where a \emph{joint} is a point lying on $d$ lines whose directions are linearly independent.
The problem was first raised by Chazelle et al. \cite{CEGPSSS92}, where they provided the following example with many joints.
\begin{example}\label{example:tight}
    Pick $M\geq d$ hyperplanes \emph{in general position}.
    In other words, we pick $M$ hyperplanes whose $d$-wise intersections are all distinct points.
    Take $\cL$ to be the set of $\binom{M}{d-1}$ lines that are the $(d-1)$-wise intersections.
    Then the set $\cJ$ of $d$-wise intersections is the set of joints formed by $\cL$, and $\abs{\cJ} = \binom{M}{d}$.
\end{example}
The problem has drawn lots of attention since then not only because of its own merit, but also because of its connection to the Kakeya problem in harmonic analysis as observed by Wolff \cite{Wol99}.
However, the joints problem remained unsolved for a long time until Guth and Katz \cite{GK10} resolved it up to a constant in $\RR^3$, showing that $L$ lines in $\RR^3$ can form $O(L^{3/2})$ joints.
The technique they used is the polynomial method, introduced by Dvir in his celebrated work on the finite field Kakeya problem \cite{Dvir09}.
The proof was quickly generalized to arbitrary dimensions independently by Kaplan--Sharir--Shustin \cite{KSS10} and Quilodr\'an \cite{Qui09}, which can be extended to arbitrary fields too (see \cite{CI14, Dvir10, Tao14}).

The Guth--Katz joints theorem can also be seen as the precursor of two other important works in the field.
One is the almost resolution of the Erd\H{o}s distinct distanced problem \cite{GK15}, where Guth and Katz used the polynomial method together with many other tools to show that $N$ points in $\RR^2$ determine at least $\Omega(N/\log N)$ different distances.
The other one is Guth's resolution of the so-called endpoint case of the Bennett--Carbert--Tao multilinear Kakeya conjecture \cite{Guth10}, which can be viewed as a ``joints of tubes'' theorem.
This echoes Wolff's observation nicely and shows that techniques developed in the joints setting may help us understand the Kakeya problem.
We refer the readers to \cite[Section 15.8]{Guth-book} for a nice exposition of the endpoint case of the multilinear Kakeya conjecture.

The story of the joints problem does not end here.
Though the maximum number of joints was determined up to a multiplicative constant, its more refined behavior was not known.
Guth \cite[Section 2.5]{Guth-book} conjectured that \cref{example:tight} is the exact optimum, i.e. $\binom{M}{d-1}$ lines form at most $\binom{M}{d}$ joints.
This was recently verified by Zhao and the second author \cite{YZ23} up to an $o(1)$-factor.
The first main result of this paper is a full resolution of Guth's conjecture.
We can actually prove a slightly stronger statement, which can be viewed as the joints analogue of Lovasz's version of the Kruskal--Katona theorem.
In the following, we say that $(\cJ,\cL)$ is a \emph{joints configuration} in $\FF^d$ if $\cL$ is a set of lines in $\FF^d$ and $\cJ$ contains joints formed by $\cL$.

\begin{theorem}\label{theorem:SharpJoints}
    Let $(\cJ,\cL)$ be a joints configuration in $\FF^d$. If $J\eqdef |\cJ|=\binom{x}{d}$ for some real number $x\geq d$, then $L\eqdef |\cL|\geq \binom{x}{d-1}$.
\end{theorem}

Furthermore, we are able to show that \cref{example:tight} gives the only joints configurations that achieve the equality in \cref{theorem:SharpJoints}.

\begin{theorem}\label{theorem:Structural}
Let $(\cJ,\cL)$ be a joints configuration in $\FF^d$ and $M\geq d$ be a real number with $J=\binom{M}{d}$ and $L=\binom{M}{d-1}$. Then $M$ is an integer, and there exist $M$ hyperplanes in general position such that $\cJ$ is the set of $d$-wise intersections of the hyperplanes and $\cL$ is the set of $(d-1)$-wise intersections of the hyperplanes.
\end{theorem}

We stress that \cref{theorem:SharpJoints} and \cref{theorem:Structural} are the first known sharp bound and structural theorem obtained by the polynomial method (in Dvir's sense).
In fact, similar results in incidence geometry are really rare.
The only example that immediate comes to mind is the work on ordinary lines by Green and Tao \cite{GT13}.
In that work, Green and Tao determined the exact minimum number of ordinary lines when the number of points is sufficiently large, and they also proved a strong structural and stability result using the algebraic structure provided by cubic curves.
Though there are some more refined bounds obtained by the polynomial method recently (\cite{YZ23} on the joints problem and \cite{BC21} on the finite field Kakeya problem), the fact that we are able to obtain the sharp bound and the structural theorem here using the polynomial method should be both surprising and exciting.
Of course, the fact that \cref{theorem:Structural} provides absolutely no stability result is slightly disappointing.
We believe that this is the limitation of our own method, and we will discuss this more in detail in \cref{subsection:StructuralDiscussion}.

Lastly, we make a remark here that \cref{theorem:SharpJoints} gives a genuinely new (though way too long) proof of Lovasz's version of the Kruskal--Katona theorem.
It also verifies some conjectures of Bollob\'as and Eccles regarding a generalization of Lovasz's version of the Kruskal--Katona theorem \cite[Conjecture 2 and 3]{BE15}.
We will discuss the conjectures more in detail in \cref{subsection:PartialShadow}.
However, \cref{theorem:Structural}, on the other hand, does not give a new proof of the structural theorem of Lovasz's version of the Kruskal--Katona theorem, as the structural theorem will actually be invoked in the proof of \cref{theorem:Structural}.
It does nonetheless provide an analogous structural result for Bollob\'as and Eccles's conjectures.
We will discuss the connections in detail in \cref{section:Graph}.

\subsection{Joints of curves}
\cref{theorem:SharpJoints} can be generalized to joints of curves.
Here, a point $p$ is a joint in $\FF^d$ if there are $d$ algebraic curves $\ell_1,\ldots,\ell_d$ so that $p$ is a regular (or equivalently in the curve case, smooth) point of $\ell_i$ for any $i=1,\ldots, d$, and the $d$ tangent vectors at $p$ are linearly independent.
To avoid confusion, if $\cJ$ is a set of points and $\cC$ is a set of curves, we call $(\cJ,\cC)$ a \emph{joint-of-curve configuration} if $\cJ$ contains joints formed by curves in $\cC$.

\begin{theorem}\label{theorem:JointsOfCurves}
    Let $(\cJ,\cC)$ be a joint-of-curve configuration in $\FF^d$. If $J\eqdef |\cJ|=\binom{x}{d}$ for some $x\geq d$, then $\deg \cC\eqdef \sum_{\ell\in\cC} \deg\ell\geq \binom{x}{d-1}$.
\end{theorem}

Generalization of the joints problem to joints of curves is not new.
In fact, the original proof technique for the joints problem generalizes to joints of curves (see \cite{KSS10} and \cite{Qui09}).
The sharp bound for joints of curves essentially says that to get as many joints as possible with a given total degree of curves, taking only lines already suffices.
We note here that we do not have a corresponding analogue of \cref{theorem:Structural} in the setting of joints of curves, though we do believe that the same statement holds in this setting.

\subsection{Joints of varieties}
Though it is not hard to generalize the joints theorem to joints of curves, it is not the case for joints of varieties.
Here, as in \cite{TYZ22}, we will assume that each variety is irreducible without loss of generality.
In addition, a point is a joint formed by several varieties if it is a regular point of each variety, and the tangent spaces are linearly independent and span the entire space.
In fact, even the simplest case of joints of $2$-flats in $\RR^6$ was open for a long time.
A recent breakthrough of Tidor, Zhao and the second author \cite{TYZ22} settled all known variants of the joints problems for varieties up to a multiplicative constant.

\begin{theorem}[\cite{TYZ22}]\label{theorem:TYZMainThm}
    Let $k_1,\ldots,k_r,m_1,\ldots, m_r$ be positive integers.
    For each $i=1,\ldots, r$, let $\cV_i$ be a finite multiset of $k_i$-dimensional varieties in $\FF^d$, where $d=k_1m_1+\cdots+k_rm_r$.
    Let $\cJ$ be the set of joints $p$ that are formed by choosing $m_i$ unordered elements from $\cV_i$ for each $i=1,\ldots, r$, and we denote by $M(p)$ the number of such choices for each $p$.
    \begin{enumerate}[(a)]
        \item(without multiplicities) The number of joints can be bounded as
        \[J\eqdef \abs{\cJ}=O_{k_1,\ldots,k_r,m_1,\ldots,m_r}\left(\left(\left(\deg\cV_1\right)^{m_1}\cdots\left(\deg\cV_r\right)^{m_r}\right)^{1/(m_1+\cdots+m_r-1)}\right).\]
        \item (with multiplicities) Summing over all joints,
        \[\sum_{p\in\cJ}M(p)^{1/(m_1+\cdots+m_r-1)}=O_{k_1,\ldots,k_r,m_1,\ldots,m_r}\left(\left(\left(\deg\cV_1\right)^{m_1}\cdots\left(\deg\cV_r\right)^{m_r}\right)^{1/(m_1+\cdots+m_r-1)}\right),\]
        where the implicit constant is slightly larger than the one in (a).
    \end{enumerate}
\end{theorem}
As the main focus of this paper is the sharp bound and the structural theorem for joints of lines, we refer interested readers to \cite{TYZ22} for a more complete review regarding joints of varieties.

After getting \cref{theorem:SharpJoints}, it is natural to see if we can also get a sharp bound for joints of varieties.
However, unlike joints of lines where the optimal constant was already known in \cite{YZ23}, the implicit constants in \cref{theorem:TYZMainThm} are not known to be optimal, and we believe that they are not optimal.
To facilitate the discussion, let $C_{k_1,\ldots, k_r;m_1,\ldots, m_r}$ be the smallest implicit constant that would make (a) in \cref{theorem:TYZMainThm} true, and let $C_{k_1,\ldots,k_r;m_1,\ldots,m_r}^{M}$ be the one for (b).
Even in the case where $r=1$, the exact value of $C_{k_1;m_1}$ is not known when $k_1,m_1>1$.
The main theorem of \cite{TYZ22} shows that
\[C_{k_1;m_1}\leq \left(\frac{(k_1m_1)!}{k_1^{m_1}m_1^{m_1}}\right)^{1/(m_1-1)},\]
whereas priority to this paper, the best known lower bound is
\[C_{k_1;m_1}\geq C_{1;m_1}=\frac{(m_1-1)!^{1/(m_1-1)}}{m_1}.\]
As a consequence, in order to get a sharp bound for joints of varieties, it makes more sense to first determine the constant $C_{k_1,\ldots, k_r;m_1,\ldots, m_r}$.

In this paper, instead of improving the upper bound on $C_{k_1,\ldots, k_r;m_1,\ldots, m_r}$, we demonstrate a result that might explain why it is hard to improve the upper bound using the argument in \cite{TYZ22}.
In particular, we will show that it is possible to extract a new notion of multiplicity $\nu^*$ from the argument so that the argument in \cite{TYZ22} still works and gives the same constant, but this time there is a matching construction that proves the optimality of the constant.
Since the definition of $\nu^*$ would require lots of motivation, we defer the exact definition until \cref{section:multiplicity}, and we only state the main theorem of that section vaguely as following.

\begin{theorem}\label{theorem:IntroMultiplicity}
    In the setting of \cref{theorem:TYZMainThm}, there is a notion of multiplicity $\nu^*$ so that the same argument gives
    \[\sum_{p\in\cJ}\nu^*(p)^{(m_1+\cdots+m_r)/(m_1+\cdots+m_r-1)}\leq C^{\nu^*}_{k_1,\ldots, k_r;m_1,\ldots, m_r}\left(\left(\deg\cV_1\right)^{m_1}\cdots\left(\deg\cV_r\right)^{m_r}\right)^{1/(m_1+\cdots+m_r-1)},\]
    where the constant
    \[C^{\nu^*}_{k_1,\ldots, k_r;m_1,\ldots, m_r}=\left(\frac{d!}{\prod_{i=1}^{r}k_i!^{m_i}m_i^{m_i}}\right)^{1/(m_1+\cdots+m_r-1)}\]
    is optimal.
\end{theorem}

This shows that if one believes the current upper bound on $C_{k_1,\ldots, k_r;m_1,\ldots,m_r}$ needs improvement, then one needs to really distinguish ``simple joints'' (where we only count each joint once) from ``joints with multiplicities'' (where each joint $p$ contributes $\nu^*(p)^{(m_1+\cdots+m_r)/(m_1+\cdots+m_r-1)}$ to the sum).
We still do not know how one might do this in the joints setting, and even the corresponding graph-theoretic problems already seem interesting and challenging.
We will discuss this more in \cref{subsection:GIJointsConfig} and \cref{section:Graph}.

The implication of \cref{theorem:IntroMultiplicity} is not all negative.
In fact, now both (a) and (b) of \cref{theorem:TYZMainThm} come as an immediate consequence with the constants stated in \cite{TYZ22}: the inequality between the two constants $C_{k_1,\ldots, k_r;m_1,\ldots, m_r}\leq C^{\nu^*}_{k_1,\ldots, k_r;m_1,\ldots, m_r}$ is immediate by the fact that $\nu^*(p)\geq 1$, and we will prove that
\[C^{M}_{k_1,\ldots, k_r;m_1,\ldots, m_r}\leq \left(\prod_{i=1}^{r}\frac{m_i^{m_i}}{m_i!}\right)^{1/(m_1+\cdots+m_r-1)}C^{\nu^*}_{k_1,\ldots, k_r;m_1,\ldots, m_r}\]
by comparing $\nu^*(p)^{(m_1+\cdots+m_r)}$ and $M(p)$.
Those recover the main results of \cite{TYZ22}, and thus \cref{theorem:IntroMultiplicity} can be seen as a unification of the two statements.
In addition, the new multiplicity $\nu^*$ allows us to prove some statements that are more refined than the original joints theorem.
We include those at the end of \cref{section:multiplicity} for readers who are interested.

\subsection{Generically induced joints configuration}\label{subsection:GIJointsConfig}
Since the tight configuration \cref{example:tight} is determined by several hyperplanes, it is natural to generalize the construction and see if it gives good lower bounds on the constants $C_{k_1,\ldots, k_r;m_1,\ldots, m_r}$.
This type of constructions has already been considered by Zhao and the second author in \cite{YZ23}, where they call this type of configurations \emph{generically induced}.
Formally, for any hyperplanes $H_1,\dots,H_n$ in $\FF^d$ in general position, the intersection $\cap_{i\in F}H_i$ is a $k$-flat for any $(d-k)$-subset $F\subseteq [n]$. Similarly, each $d$-subset $P$ gives a point in $\FF^d$. A joints configuration is called generically induced if the flats and the joints are constructed in this way. 

Using this language, we may restrict ourselves to generically induced configurations and ask what the best possible bound we can get is. 
This way, the problem is simplified to a purely set-theoretic problem. 
In \cite{YZ23}, Zhao and the second author conjectured that the answer to this question is the same as the answer to the corresponding joints problem in certain cases (multijoints of lines and joints of lines to be more specific).
It is tempting to make the same conjecture for all variants of the joints problem, but it turns out that there is some other small configuration that beats generically induced configurations in some case.
Indeed, we may take a generically induced configuration in some higher-dimensional space and take a general projection back to the original space.
We call this type of configurations \emph{projected generically induced configurations}.
In \cref{section:Graph}, we give an example of a projected generically induced configuration that outperforms all genrically induced configurations with the same number of joints.
This example comes from the corresponding set-theoretic problem that Bollob\'as and Eccles considered \cite{BE15}.
We suspect that this problem has the same answer as the original joints problem.

As the set-theoretic version is conjectured to give the same answer as the joints problem, we investigate the optimal constant in the set-theoretic version.
On the lower bound (construction) side, we show that $C_{2;3}\geq \sqrt{2/7}$ using a generically induced joints configuration, which can be used to disprove \cite[Conjecture 1.11]{TYZ22} saying that $N$ $2$-flats in $\FF^6$ form at most $(\sqrt{2}/3 + o(1))N^{3/2}$ joints.
We also give a lower bound for $C_{k,k,k;1,1,1}$, showing that it grows exponentially.

On the other hand, in a related work \cite{CY23}, we proved some sharp upper bounds for the number of joints in the generically induced configurations via entropy method in some cases, some of which are better than the upper bounds obtained from the polynomial method. 
We refer the readers to that paper for a more in-depth discussion.
\subsection{The partial shadow problem}\label{subsection:PartialShadow}
It turns out that when restricted to projected generically induced joints configuration, the joints problem for lines becomes the partial shadow problem that Bollob\'as and Eccles proposed in \cite{BE15}.
We restate the problem here.
Let $r,m,k$ be three positive integers.
The problem asks for the minimum size of family $\cB\subseteq \binom{\NN}{r-1}$ so that there is a family $\cA\subseteq \binom{\NN}{r}$ of size $m$ satisfying the following: for each $A\in \cA$, there are at most $k$ subsets in $A$ of size $r-1$ missing in $\cB$.
Bollob\'as and Eccles denoted this answer by $f(r,m,k)$.
In their paper, they conjectured that the Lovasz's bound still holds in this setting \cite[Conjecture 3]{BE15}.
Conjecture 2 in their paper is a weaker conjecture, though it suggests what all the tight examples should be.
We can prove both conjectures with \cref{theorem:SharpJoints,theorem:Structural} by applying them to projected generically induced joints configuration.

\begin{theorem}\label{theorem:PartialShadow}
    For any positive integers $r,m$, nonnegative integers $k$ with $k\leq r$ and real number $x\geq r-k$ with $m = \binom{x}{r-k}$, we have $f(r,m,k)\geq \binom{x}{r-k-1}$.
    Moreover, if $\cA, \cB$ are families that certify that $f(r,m,k)=\binom{x}{r-k-1}$, then there exists disjoint $X,Y\subseteq \NN$ such that $\abs{X}=x$, $\abs{Y}=k$, and 
    \begin{align*}
        \cA &=\left\{S\cup Y: S\in \binom{X}{r-k}\right\},\\
        \cB &=\left\{S\cup Y: S\in \binom{X}{r-k-1}\right\}.
    \end{align*}
\end{theorem}

Previously, this was only known for $r-k\leq 3$ by Fitch \cite[Chapter 3]{Fitch18}.
We note that the case $r-k=3$ is implicit in Fitch's argument.
He actually proved a stronger result that the Kruskal--Katona bound still holds when $r-k=3$ and $m$ is sufficiently large.

It is surprising to us that the polynomial method is so far the only known proof of \cref{theorem:PartialShadow}, which only involves set system.
It is quite possible that other proofs exist, especially the ones that use entropy in the way it is used by Fitch \cite{Fitch18} or by us in the other paper \cite{CY23}.

\subsection{Outline}
The paper is structured as follows. 
In \cref{section:Preliminaries}, we recall several definitions and properties we need.
In \cref{section:Sharp,section:Structural}, we prove the tight bound (\cref{theorem:SharpJoints}) and the structural theorem (\cref{theorem:Structural}) for joints of lines respectively. 
We generalize the tight bound to joints of curves in \cref{section:Curve}. 
In \cref{section:multiplicity} we discuss two new notions of multiplicities on joints of varieties.
Lastly, in \cref{section:Graph}, we state the set-theoretic versions of the joints problem and discuss some special cases, including the proof of \cref{theorem:PartialShadow}.

\section{Preliminaries}\label{section:Preliminaries}
In this section, we list out several notations and tools that will be necessary throughout the paper.

\subsection{Notation}
Through out this paper, the notation $\vec{*}$ denotes a vector, and its coordinates are always denoted by $*_i$ for positive integers $i$ if the vector is indexed by positive integers, or $*_s,\,s\in S$ if the vector is indexed by a finite set $S$. 
We will also denote by $\abs{\vec{*}}$ the sum $\sum_i *_i$.

We say that a polynomial vanishes identically on $\FF^d$ (or on a line) if the polynomial on $\FF^d$ (or restricted to the line, resp.) is the zero polynomial.
Note that when $\FF$ is a finite field, this is different from just vanishing everywhere.
Similarly, when $\ell$ is a curve, we say that a polynomial vanishes identically on $\ell$ if the polynomial vanishes on every closed points of $\ell$ when $\ell$ is viewed as a scheme (or equivalently, the polynomial vanishes everywhere on $\ell$ even after extending the base field).

Given a joints configuration $(\cJ,\cL)$, for each joint $p\in \cJ$, we always fix $d$ lines $\ell_{p,1},\dots,\ell_{p,d}$ through $p$ with linearly independent directions. Moreover, for each line $\ell\in\cL$, we fix a vector $\vec{e}_{\ell}$ parallel to $\ell$.

Similarly, given a joint-of-curve configuration $(\cJ,\cC)$, for each joint $p\in\cJ$, we fix labels of the $d$ curves $\ell_{p,1},\ldots, \ell_{p,d}$ that form a joint at $p$.
We will also fix a tangent vector $\vec{e}_{\ell,p}$ of $\ell$ at $p$ for every pair $(p,\ell)\in\cJ\times\cC$ with $p\in\ell$.

In both cases, we write $p\in\ell$ or $\ell\ni p$ if $\ell = \ell_{p,i}$ for some $i\in[d]$.

\subsection{Hasse Derivative}
Since we are working in an arbitrary field $\FF$, we need the notion of Hasse derivative.
Hasse derivative works better than formal derivative in this setting as it allows us to avoid issues in fields with positive characteristic.
For a more detailed introduction of Hasse derivative, we refer the readers to \cite[Section 2]{DKSS13}.

\begin{definition}[Hasse Derivative]
    Let $f\in \FF[x_1,\dots,x_d]$. The Hasse derivatives $(\Hasse^{\vec{\alpha}}f)_{\vec{\alpha}\in\ZZ_{\geq 0}^d}$ are defined via the expansion 
    \[f(x+y)=\sum_{\vec{\alpha}\in\ZZ_{\geq 0}^d}\left(\Hasse^{\vec{\alpha}}f(x)\right)y^{\vec{\alpha}}.\]
\end{definition}
\begin{remark}
    Hasse derivatives are preserved under shifting. Namely, we have $\Hasse^{\vec{\alpha}}f(x)=\Hasse^{\vec{\alpha}}(f\circ A)(x-p)$, where $A$ is the map satisfying $A(x)=x+p$.
\end{remark}

\begin{definition}[Multiplicity]\label{definition:Mult}
    Let $f\in \FF[x_1,\dots,x_d]$, $p\in\FF^d$, and $\ell$ be a line in $\FF^d$. We say that the multiplicity of $f$ at $p$ is at least $v$ if $\Hasse^{\vec{\alpha}}f(p)=0$ holds for all $\vec{\alpha}\in\ZZ_{\geq 0}^d$ with $\abs{\vec{\alpha}}<v$, and denote by $\mult(f,p)$ the largest such $v$.
    Similarly, $\mult(f,\ell)$ is the largest nonnegative integer $v$ such that $\left.\left(\Hasse^{\vec{\alpha}}f\right)\right|_{\ell}=0$ holds for all $\vec{\alpha}\in\ZZ_{\geq 0}^d$ with $\abs{\vec{\alpha}}<v$.
    The multiplicities are set to $\infty$ if the polynomial is zero.
\end{definition}

The most crucial way we will be using the multiplicities is the following proposition that implies vanishing of a polynomial on a line.

\begin{proposition}\label{proposition:BasicVanishing}
    Let $f\in\FF[x_1,\ldots,x_d]$ be a polynomial, and $\ell$ be a line in $\FF^d$.
    If
    \[\sum_{p\in\ell}\mult(f,p)>\deg f,\]
    then $f$ vanishes identically on $\ell$.
\end{proposition}
\begin{proof}
    When restricted to $\ell$, the polynomial $f$ becomes a univariate polynomial.
    The number of roots of $f|_{\ell}$ at $p$ is at least $\mult(f,p)$.
    To see this, we simply expand
    \[f(p+t\cdot \vec{e}_{\ell}) = \sum_{\vec{\alpha}\in\ZZ^d_{\geq 0}}\left(H^{\vec{\alpha}}f(p)\right)t^{\abs{\vec{\alpha}}}\left(\vec{e}_{\ell}\right)^{\vec{\alpha}}.\]
    By the definition of multiplicities, we see that the coefficient of $t^i$ for all $i<\mult(f,p)$ is zero, as desired.

    As the sum of the number of roots (counted with multiplicities) can never exceed the degree of a nonzero univariate polynomial, if $\sum_{p\in\ell}\mult(f,p)>\deg f\geq \deg \left.f\right|_{\ell}$, then $f$ must be identically zero on $\ell$.
\end{proof}

The following proposition says that the concept of multiplicity is preserved under any invertible affine transformation. 

\begin{proposition}
    Assume $f\in \FF[x_1,\dots,x_d]$, $p\in\FF^d$, $\ell$ is a line, and $A$ is an invertible affine transformation. It follows that 
    \begin{align*}
    \mult(f,p)=&\mult(f\circ A, A^{-1}(p)),\\
    \mult(f,\ell)=&\mult(f\circ A, A^{-1}(\ell)).
    \end{align*}
\end{proposition}
\begin{proof}
    The first part is an immediate corollary of \cite[Proposition 6]{DKSS13}. The second part can be proved by expanding both sides. Here, we include a short proof that requires $\ell$ to contain many points. We would thus consider a field extension if needed. 
    
    If $\FF$ is a finite field, we may embed everything into its field extension, and the multiplicities $\mult(f,\ell),\mult(f\circ A,A^{-1}(\ell))$ remain the same. Thus, we may assume without loss of generality that $|\FF|$ is infinite in this proof. From the definition of multiplicity, we know that $\mult(f,p)\geq \mult(f,\ell)$ for every $p\in\ell$. It follows that $\mult(f\circ A,A^{-1}(p))\geq \mult(f,\ell)$ for every $p\in\ell$. From the definition, we know that $\Hasse^{\vec{\alpha}}(f\circ A)(A^{-1}(p))$ is zero for every $p\in\ell$ and every $\vec{\alpha}\in\ZZ_{\geq 0}^d$ with $|\vec{\alpha}|<\mult(f,\ell)$. Since $\ell$ contains infinitely many points, we know that $\Hasse^{\vec{\alpha}}(f\circ A)$ restricting to $A^{-1}(\ell)$ is the zero polynomial by applying \cref{proposition:BasicVanishing}. Thus, $\mult(f\circ A,A^{-1}(\ell))\geq \mult(f,\ell)$. Since $A$ is invertible, we may run the same proof for $f\circ A$ and $A^{-1}$ in place of $f$ and $A$ to get the opposite inequality. Thus, $\mult(f\circ A,A^{-1}(\ell))= \mult(f,\ell)$.
\end{proof}

The following proposition and corollary are also useful, and they are proved in \cite{DKSS13}.
Note that the second part of \cref{corollary:DerivativeDecreasesMulti} is not actually stated there, but the proof for the first part also works for the second part too.

\begin{proposition}\label{proposition:HasseComposition}
    Assume $f\in \FF[x_1,\dots,x_d]$ and $\vec{\alpha},\vec{\alpha}'\in \ZZ_{\geq 0}^d$. It follows that
    \[\binom{\vec{\alpha}+\vec{\alpha}'}{\vec{\alpha}}\Hasse^{\vec{\alpha}+\vec{\alpha}'}f(x)=\Hasse^{\vec{\alpha}}\left(\Hasse^{\vec{\alpha}'}f(x)\right),\]
    where $\binom{\vec{\alpha}+\vec{\alpha}'}{\vec{\alpha}}\eqdef \prod_{i=1}^d\binom{\alpha_i+\alpha'_i}{\alpha_i}$. 
    In particular, Hasse derivatives commute.
    Moreover, we have
    \[\Hasse^{(\alpha_1,\alpha_2,\dots,\alpha_d)}f(x)=\Hasse^{(\alpha_1,0\dots,0)}\left(\Hasse^{(0,\alpha_2,\dots,\alpha_d)}f(x)\right).\]
\end{proposition}

\begin{corollary}\label{corollary:DerivativeDecreasesMulti}
    For any polynomial $f\in\FF[x_1,\ldots, x_d]$ and point $p$ or line $\ell$ in $\FF^d$, we have
    \[\mult(H^{\vec{\alpha}}f,p)\geq \mult(f,p)-\abs{\vec{\alpha}}\]
    and
    \[\mult(H^{\vec{\alpha}}f,\ell)\geq \mult(f,\ell)-\abs{\vec{\alpha}}\]
    for any $\vec{\alpha}\in\ZZ_{\geq 0}^d$.
\end{corollary}

Using Hasse derivative, we may define derivative operators associated to a joints configuration.

\begin{definition}
    Assume $(\cJ,\cL)$ is a joints configuration in $\FF^d$ and $p\in\cJ$. 
    Let $A_p$ be the linear transformation that sends the $i$-th unit vector to $\vec{e}_{\ell_{p,i}}$ for $i\in[d]$. 
    Define $D_p^{\vec{\alpha}}$ to be the operator such that
    \[D_p^{\vec{\alpha}}f(x)=\Hasse^{\vec{\alpha}}(f\circ A_p)(A_p^{-1}(x))\]
    holds for all $f\in\FF[x_1,\dots,x_d]$ and all $x\in\FF^d$.
\end{definition}

Among all the derivative operators we have set at each joint, some of them are actually the same.
To see this, we first prove the following.

\begin{proposition}\label{proposition:ChangeCoordinate}
     Assume $f\in \FF[x_1,\dots,x_d]$, $\alpha_1\in\ZZ_{\geq 0}$, and $A$ is an invertible linear transformation fixing $\vec{e}_1\eqdef(1,0,\dots,0)$.
     It follows that
     \[\Hasse^{(\alpha_1,0,\dots,0)}f(0)=\Hasse^{(\alpha_1,0,\dots,0)}(f\circ A)(0).\]
\end{proposition}
\begin{proof}
    Note that $A$ sends $(x_1,x_2,\dots,x_d)$ to $(y_1,y_2,\dots,y_d)$, where $y_1$ is a linear combination of $x_1,\dots,x_d$ and $y_2,\dots,y_d$ are linear combinations of $x_2,\dots,x_d$. It follows that both sides of the equation are equal to the coefficient of $x_1^{\alpha_1}$ in $f(x_1,0,\dots,0)$. 
\end{proof}

\begin{corollary}\label{corollary:ChangeCoordinate}
    Assume $(\cJ,\cL)$ is a joints configuration and $p,p'$ are two joints on a line $\ell$. Moreover, we may assume without loss of generality that $\ell=\ell_{p,1}=\ell_{p',1}$. It follows that the derivatives at $p$ satisfy
     \[D_p^{(\alpha_1,0,\dots,0)}f(p)=D_{p'}^{(\alpha_1,0,\dots,0)}f(p)\]
    for any $f\in \FF[x_1,\dots,x_d]$ and $\alpha_1\in\ZZ$.
\end{corollary}
\begin{proof}
    Since the Hasse derivatives are preserved under shifting, we may assume with out loss of generality that $p$ is the origin. We may apply \cref{proposition:ChangeCoordinate} with the polynomial $f\circ A_p$ and the linear transformation $A_p^{-1}\circ A_{p'}$. Note that $A_p^{-1}\circ A_{p'}$ indeed fixes $e_1$ since $\ell_{p,1}=\ell_{p',1}$. It follows that 
     \[D_p^{(\alpha_1,0,\dots,0)}f(0)=\Hasse^{(\alpha_1,0,\dots,0)}(f\circ A_p)(0)=\Hasse^{(\alpha_1,0,\dots,0)}(f\circ A_p\circ A_p^{-1}\circ A_{p'})(0)=D_{p'}^{(\alpha_1,0,\dots,0)}f(0)\]
\end{proof}

In the case of joints of curves, we use the tools developed in \cite{TYZ22}.
There, derivative operators along a variety were defined using local coordinates in the ``completion''.
We refer the readers to \cite[Section 4]{TYZ22} for a full discussion on the motivation and the actual definition.

\begin{definition}
    Let $(\cJ,\cC)$ be a joint-of-curve configuration in $\FF^d$, and fix $p\in\cJ$.
    Let $A_p$ be the linear transformation sending the $i$-th unit vector to $\vec{e}_{\ell_{p,i},p}$ and the origin to $p$.
    For every $\alpha\in\ZZ_{\geq 0}^d$, we temporarily define $\Hasse_p^{\vec{\alpha}}$ to be the derivative operator such that
    \[\Hasse^{\vec{\alpha}}_pf(x)=\Hasse^{\vec{\alpha}}(f\circ A_p)(A_p^{-1}(x)).\]
    For any $i\in[d]$ and any positive integer $a$, set
    \[D^{a\cdot \vec{e}_i}_p = \Hasse^{a\cdot \vec{e}_i}_p+\sum_{\substack{\vec{\beta}\in\ZZ_{\geq 0}^d\\\abs{\vec{\beta}}<\abs{\vec{\alpha}}}}c^{a\cdot \vec{e}_i}_{\vec{\beta}}\Hasse^{\vec{\beta}}_p\]
    where $c^{a\cdot \vec{e}_i}_{\vec{\beta}}\in\FF$ is chosen so that if $f$ is identically zero on $\ell_{p,i}$, then $D^{a\cdot \vec{e}_i}_pf$ is also identically zero on $\ell_{p,i}$ (see Section 4 of \cite{TYZ22}).
    For every $\alpha\in\ZZ_{\geq 0}^d$, we now let
    \[D^{\vec{\alpha}}_p = D^{\alpha_1\cdot\vec{e}_1}_p\cdots D^{\alpha_d\cdot\vec{e}_d}_p.\]
\end{definition}

The multiplicity of a polynomial at a curve is similarly defined as in \cref{definition:Mult}.

\begin{definition}[Multiplicity for curves]
Let $\ell$ be a curve in $\FF^d$, and $f$ be a polynomial in $\FF[x_1,\ldots, x_d]$.
The multiplicity $\mult(f,\ell)$ of $f$ at $\ell$ is the largest nonnegative integer $v$ such that $\left.\left(\Hasse^{\vec{\alpha}}f\right)\right|_{\ell}=0$ holds for all $\vec{\alpha}\in\ZZ_{\geq 0}^d$ with $\abs{\vec{\alpha}}<v$.
It is set to $\infty$ if $f$ is zero.
\end{definition}

\cref{proposition:BasicVanishing} can also be generalized to curves.
To do so, we will apply B\'ezout's theorem.

\begin{proposition}\label{proposition:BasicVanishingCurves}
    Let $f\in\FF[x_1,\ldots,x_d]$ be a polynomial, and $\ell$ be an irreducible algebraic curve in $\FF^d$.
    If
    \[\sum_{p\in\ell}\mult(f,p)>\deg\ell\cdot \deg f,\]
    then $f$ is identically zero on $\ell$.
\end{proposition}
\begin{proof}
    Assume for the sake of contradiction that $f$ is not identically zero on $\ell$.
    Then $V(f)\cap \ell$, as a scheme, has degree at most $\deg\ell\cdot \deg f$ by B\'ezout's theorem.
    Now, if we can prove that $V(f)\cap \ell$ has degree at least $\mult(f,p)$ at $p$ for each $p\in\ell$, then we would reach a contradiction.
    Let $m\eqdef \mult(f,p)$, and let $\fm_{p,\FF^d}\subseteq \FF[x_1,\ldots, x_d]$ be the maximal ideal corresponding to $p$ in $\FF^d$.
    Then it is clear that $f\in \fm_{p,\FF^d}^m$, and from here it is not hard to get that $V(f)\cap \ell$ has degree at least $m$ at $p$.
    We include a brief algebro-geometric argument here for completeness.

    Let $R_{\ell}$ be the coordinate ring of $\ell$, and let $\fm_{p,\ell}$ be the image of $\fm_{p,\FF^d}$ under the projection $\FF[x_1,\ldots, x_d]\to R_{\ell}$.
    As the image of $f$ lives in $\fm_{p,\ell}^m$, it suffices to show that $R_{\ell}/\fm_{p,\ell}^m$ is at least $m$-dimensional as a $\FF$-vector space.
    This is true as $\ell$ is $1$-dimensional.
\end{proof}

We end with the following proposition that will be needed in the argument.

\begin{proposition}
    Let $(\cJ,\cC)$ be a joint-of-curve configuration in $\FF^d$.
    Let $p\in\cJ$ be a joint, and $f\in\FF[x_1,\ldots, x_d]$ be a polynomial.
    If $\vec{\alpha}\in\ZZ_{\geq 0}^d$ satisfies that $\abs{\vec{\alpha}}-\alpha_i<\mult(f,\ell_{p,i})$, then $D^{\vec{\alpha}}_pf(p)=0$.
\end{proposition}
\begin{proof}
    Since Hasse derivatives commute, we know that 
    \[D^{\vec{\alpha}}_p = D^{\alpha_i\cdot\vec{e}_i}_pD^{\vec{\alpha}-\alpha_i\cdot\vec{e}_i}_p\]
    as the derivative operators $\left(D^{\alpha_j\cdot \vec{e}_j}_p\right)_{j\in[d]}$ are all linear combinations of Hasse derivatives and thus commute.
    Since $\abs{\vec{\alpha}-\alpha_i\cdot\vec{e}_i}=\abs{\vec{\alpha}}-\alpha_i<\mult(f,\ell_{p,i})$, we know that $D^{\vec{\alpha}-\alpha_i\vec{e}_i}_pf$ is identically zero on $\ell_{p,i}$.
    Therefore $D^{\alpha_i\cdot\vec{e}_i}_pD^{\vec{\alpha}-\alpha_i\cdot\vec{e}_i}_pf$ is identically zero on $\ell_{p,i}$, and in particular it is zero at $p$.
\end{proof}

\subsection{Shannon Entropy}
In \cref{section:multiplicity}, we give two new definitions of multiplicity of joints, and one of the definitions is related to Shannon entropy. 
For completeness, we include the definition of entropy and conditional entropy, and also some properties we need.
\begin{definition}
    For any discrete random variable $X$ with finite support $S$, its entropy is defined as
    \[H(X)\eqdef \sum_{x\in S}-p_X(x)\log_2p_X(x)\]
    where we denote by $p_X(x)$ the probability $\PP(X=x)$.
    
    The entropy $H(X_1,\ldots, X_n)$ of several discrete random variables $X_1,\ldots, X_n$ is defined similarly.
    To be more specific, it is defined as the entropy of $\mathbf{X}=(X_1,\ldots, X_n)$.
\end{definition}

\begin{definition}
    Let $X,Y$ be two random variables with supports $S,T$, respectively, that are both finite.
    Then
    \[H(X\mid Y)\eqdef \sum_{(x,y)\in S\times T}-p_{X,Y}(x,y)\log_2\left(\frac{p_{X,Y}(x,y)}{p_Y(y)}\right)=\sum_{y\in T}p_Y(y)H(X\mid Y=y),\]
    where as before, $p_{X,Y}(x,y)$ denotes the probability $\PP(X=x, Y=y)$ and $p_Y(y)$ denotes the probability $\PP(Y=y)$.
\end{definition}

The only properties we need are the following commonly used properties.

\begin{proposition}[Chain rule]
    For any two random variables $X,Y$ with finite supports, we have $H(X,Y) = H(Y)+H(X\mid Y)$.
    More generally, for any random variables $X_1,\ldots, X_n$ all with finite supports,
    \[H(X_1,\ldots, X_n) = \sum_{i=1}^{n}H(X_i\mid X_1,\ldots, X_{i-1}).\]
\end{proposition}

\begin{proposition}[Subadditivity]
    For any list of random variables $X_1,X_2,\ldots, X_n$ with finite supports,
    \[H(X_1,X_2,\ldots,X_n)\leq \sum_{i=1}^{n}H(X_i).\]
    The equality holds if and only if all the random variables are mutually independent.
\end{proposition}

For a proof of those two statements, we refer the readers to \cite[Section 14.6]{AS00}

\section{Sharp bound}\label{section:Sharp}
In this section, we start with an overview of the method in \cite{YZ23} and see where we might be able to improve.
This would be the first motivating idea, and we will dive into the proof with that idea in the rest of the section.
\subsection{Motivating idea}\label{subsection:KeyIdea}
As one would imagine, the motivating idea is to look into the argument in \cite{YZ23} and see if there is some room for improvement.
Therefore, we will start with a quick review of the proof strategy in that paper.

The first key component is a vanishing lemma tailored to joints configurations.
The vanishing lemma can be thought of as a strengthening of the following basic fact (c.f. \cref{proposition:BasicVanishing}): if $p_1,\ldots, p_m$ are $m$ points on a line $\ell$ in $\FF^d$, $b_1,\ldots, b_m\geq 0$ are real numbers with $\sum_ib_i=1$, and $f$ is a polynomial of degree less than $n$ such that $\mult(f,p_i)\geq b_in$, then $f$ is identically zero on $\ell$.
We will first state the vanishing lemma in our notation, and then we will highlight several components and features of the statement for the discussion.
For now, it is more important to know what the statement looks like than knowing the exact details.

\begin{lemma}[c.f. Lemma 2.1 and 2.3 in \cite{YZ23}]\label{lemma:YZ}
    Let $(\cJ,\cL)$ be a joints configuration in $\FF^d$.
    For each $p\in\cJ$, assign a real number $a_p$.
    Assume that for each $p\in\cJ$ and $i=1,\ldots,d$ we may choose $b_{p,i}\geq 0$ satisfying 
    \begin{enumerate}[(a)]
        \item $b_{p,i}-a_p=b_{p',i'}-a_{p'}$ if $\ell_{p,i}=\ell_{p',i'}$, and
        \item for any fixed $\ell\in\cL$, 
        \[\sum_{\substack{p\in\cJ, i\in[d]\\ \ell_{p,i}=\ell}}b_{p,i}=1.\]
    \end{enumerate}
    Then for every polynomial $f\in \FF[x_1,\ldots, x_d]$ with degree less than $n$, if 
    \[D^{\vec{\alpha}}_pf(p)=0\]
    holds for every $\vec{\alpha}\in\ZZ_{\geq 0}^d$ with $\alpha_i< b_{p,i}n$ for any $i=1,\ldots,d$, then $f$ must be zero.
\end{lemma}

The $b_{p,i}$'s play the roles of the $b_i$'s in the basic fact stated earlier, and the connection would be clearer once we go through the proof.
One other important feature is that there are parameters (the $a_p$'s) that we can vary.
This will be useful in the later steps.

The second component is usually called ``parameter counting''.
Parameter counting is the core idea of the polynomial method and can be found in most, if not all, proofs regarding the joints problem.
In this case, the parameter counting step is done as follows.
The space of polynomials in $\FF[x_1,\ldots, x_d]$ with degree less than $n$ is an $\FF$-vector space of dimension $\binom{n+d-1}{d}$.
For each $p\in\cJ$ and $\vec{\alpha}\in\ZZ_{\geq 0}^d$, we can verify that $D^{\vec{\alpha}}_pf(p)$ is a linear functional on the space of $f$.
The vanishing lemma imposes a certain number of linear conditions on $f$, and those force $f$ to be zero.
It is then a simple linear algebra fact that the number of linear conditions is at least the dimension of the space.
By counting the number of constraints, the vanishing lemma gives
\[\sum_{p\in\cJ}\prod_{i=1}^{d}\lceil b_{p,i}n\rceil \geq \binom{n+d-1}{d},\]
and by taking $n$ to infinity and comparing the leading coefficients, we get
\[\sum_{p\in\cJ}\prod_{i=1}^{d}b_{p,i}\geq \frac{1}{d!}.\]
Therefore, by parameter counting, we get an inequality from the vanishing lemma.

The very last component of the argument in \cite{YZ23} is to choose the parameters $a_p$'s nicely so that the inequality allows us to bound the number of joints in terms of the number of lines.
This is why we emphasized the role the $a_p$'s play in the vanishing lemma earlier.
In \cite{YZ23}, under a mild assumption, they chose the $a_p$'s so that the product $\prod_{i=1}^{d}b_{p,i}$ is the same for every joint $p\in\cJ$.
This way, each joint contributes the same to the left hand side of the inequality, and the rest is a computation and some clever applications of the AM-GM inequality.

Now, to see where we can improve in the argument above, consider the tight configuration \cref{example:tight}.
Then by a simple counting argument, we may conclude that there are $M-d+1$ joints on each line.
Therefore, if we set $a_p=0$ for any $p\in\cJ$, we get $b_{p,i}=(M-d+1)^{-1}$ in \cref{lemma:YZ}.
The inequality parameter counting gives implies that
\[J\cdot (M-d+1)^{-d}\geq \frac{1}{d!},\]
and so $J\geq (M-d+1)^d/d!$.
This is slightly off from the truth $J= \binom{M}{d}$.
This suggests that some of the linear conditions we put might be redundant.
To see what the redundant ones are, we will need to look at the proof of \cref{lemma:YZ}.
We will rewrite the proof using the notations set up in this paper.

\begin{proof}[Proof of \cref{lemma:YZ}]
    If $f$ is nonzero, we may choose $p\in \cJ$ and $\vec{\alpha}\in\ZZ_{\geq 0}^d$ with the smallest $\abs{\vec{\alpha}}-a_pn$ so that $D^{\vec{\alpha}}_pf(p)\neq 0$.
    By the minimality of $(p,\vec{\alpha})$, we know that $\mult(f,p')\geq \abs{\vec{\alpha}}+(a_{p'}-a_p)n$ for every $p'\in\cJ$.
    By the assumption, there must exist some $i\in[d]$ such that $\alpha_i\geq b_{p,i}n$.
    Let $\vec{\alpha}' = \vec{\alpha}-\alpha_i\cdot\vec{e}_i$.
    Then for every $p'\in\ell_{p,i}$, if $\ell_{p,i}=\ell_{p',i'}$, we have
    \[\mult(D^{\vec{\alpha}'}_pf,p')\geq \mult(f,p')-\abs{\vec{\alpha}'}\geq \abs{\vec{\alpha}}+(a_{p'}-a_p)n-\abs{\vec{\alpha}}+\alpha_i=\alpha_i+(a_{p'}-a_p)n.\]
    Note that by the conditions, we have $\alpha_i+(a_{p'}-a_p)n\geq (b_{p,i}-a_p+a_{p'})n=b_{p',i'}n$.
    This shows that $\sum_{p'\in\ell_{p,i}}\mult(D^{\vec{\alpha}'}_pf,p')\geq n> \deg f\geq \deg D^{\vec{\alpha}'}_pf$.
    By \cref{proposition:BasicVanishing}, this shows that $D^{\vec{\alpha}'}_pf$ vanishes identically on $\ell_{p,i}$.
    In particular, we have $D^{\vec{\alpha}}_pf(p) = D^{\alpha_i\cdot \vec{e}_i}_pD^{\vec{\alpha}'}_pf(p)=0$, which is a contradiction.
    Therefore, $f$ must be zero.
\end{proof}

The inequality $\deg f\geq \deg D^{\vec{\alpha}'}_pf$ looks wasteful: if $\vec{\alpha}'$ is large, maybe we can weaken the condition $\alpha_i\geq b_{p,i}n$ and still reach the conclusion that $D^{\vec{\alpha}'}_pf$ vanishes identically on $\ell_{p,i}$.
Indeed, when $(\cJ,\cL)$ is the tight configuration and $a_p=0$ for all $p\in\cJ$, to have $\sum_{p'\in\ell_{p,i}}\mult(D^{\vec{\alpha}}_pf,p')>\deg D^{\vec{\alpha}'}_pf$ hold in the proof, we just need
\[(M-d+1)\alpha_i\geq n-\abs{\vec{\alpha}'},\]
or equivalently,
\[\alpha_i+\frac{1}{M-d+1}(\abs{\vec{\alpha}}-\alpha_i)\geq \frac{n}{M-d+1}.\]
What we can conclude from this discussion is that we can still prove that $f$ is zero if we only assume that $D^{\vec{\alpha}}_pf(p)=0$ whenever $\alpha_i+(\abs{\vec{\alpha}}-\alpha_i)/(M-d+1) < n/(M-d+1)$ holds for all $i=1,\ldots, d$.
To simplify the notation, we make the following definition.

\begin{definition}\label{definition:SpVanishing}
    Let $(\cJ,\cL)$ be a joints configuration in $\FF^d$ and $f$ be a polynomial on $\FF^d$.
    For a subset $S\subseteq \ZZ_{\geq 0}^d$ and a joint $p\in\cJ$, $f$ is said to be \emph{$S$-vanishing at $p$} if $D^{\vec{\alpha}}_pf(p)=0$ holds for any $\vec{\alpha}\in S$.
    The polynomial $f$ is said to be \emph{$(S_p)_{p\in\cJ}$-vanishing} if it is $S_p$-vanishing at every $p\in\cJ$.
\end{definition}

With this definition, we can rewrite what \cref{lemma:YZ} and the discussion above imply when $(\cJ,\cL)$ is the joints configuration in \cref{example:tight} and $a_p=0$ for all $p\in\cJ$.
Let $B$ be the box $\{\vec{\alpha}\in\ZZ_{\geq 0}^d: \alpha_i< n/(M-d+1)\,\forall i\in[d]\}$, and let $B'$ be the ``shaved'' box 
\[\left\{\vec{\alpha}\in\ZZ_{\geq 0}^d: \alpha_i+\frac{1}{M-d+1}(\abs{\vec{\alpha}}-\alpha_i)< \frac{n}{M-d+1}\,\forall i\in[d]\right\}.\]
\cref{lemma:YZ} says that if $f$ is a polynomial of degree less than $n$ that is $(B)_{p\in\cJ}$-vanishing, then it is zero.
The computation above shows that we actually just need $f$ to be $(B')_{p\in\cJ}$-vanishing in order for it to be true.
It is clear that $B'$ is strictly smaller than $B$, so the inequality we get from parameter counting would be better.
In fact, one can show that $\abs{B'} = \left(\left(M(M-1)\cdots (M-d+1)\right)^{-1}+o(1)\right)n^d$, and so the inequality obtained by parameter counting now gives $J\geq \binom{M}{d}$, which is tight now.
Therefore, the key idea is that instead of the ``box-vanishing'' conditions put in \cite{YZ23}, we will ``shave'' the box by a bit and show that the vanishing lemma still holds.
This allows us to prove a slightly stronger inequality with parameter counting.

Though the motivating idea is really simple, carrying out the detail is more technical compared to \cite{YZ23}.
The main reason is that the size of a box is easy to compute and control, while the size of a shaved box is much harder to compute in general.
To circumvent the technicalities, we will actually prove a much more general version of \cref{lemma:YZ} that incorporates one of the key ideas of \cite{TYZ22}.
In particular, we will vastly generalize the concept of ``priority order'' that was used in \cite{TYZ22}.
There will be several new ideas needed to tackle various technicalities that come up in the proof, and we will mention them when they appear.

\subsection{Main proof}
We begin with proving a general and refined vanishing lemma, and we will gradually add in more constraints to actually determine the vanishing conditions we would like to impose.
To motivate many of the notations we will use, we start with the following basic vanishing lemma.

\begin{lemma}\label{lemma:BasicVanishing}
    Let $n$ be a nonnegative integer.
    Let $(\cJ,\cL)$ be a joints configuration in $\FF^d$, and let $\vec{r}\in \ZZ_{\geq 0}^{\cJ}$.
    If $f\in \FF[x_1,\ldots, x_d]$ is a polynomial of degree at most $n$ such that
    \[\mult(f,p)\geq r_p\quad\forall p\in \cJ,\]
    then
    \[\mult(f,\ell)\geq \frac{\sum_{p\in \ell}r_p-n}{\#\{p\in\ell\}-1}\quad\forall \ell\in\cL.\]
\end{lemma}

\begin{proof}
    For any $\vec{\alpha}\in \ZZ_{\geq 0}^d$, we have $\deg \Hasse^{\vec{\alpha}}f\leq n-\abs{\vec{\alpha}}$ and $\mult(\Hasse^{\vec{\alpha}}f,p)\geq r_p-\abs{\vec{\alpha}}$.
    Note that if $\sum_{p\in \ell} \mult(\Hasse^{\vec{\alpha}}f,p)>\deg\Hasse^{\vec{\alpha}}f$, then we must have $\left.\Hasse^{\vec{\alpha}}f\right|_{\ell}\equiv 0$.
    Rearranging, we have that $\left.\Hasse^{\vec{\alpha}}f\right|_{\ell}\equiv 0$ whenever
    \[\abs{\vec{\alpha}} < \frac{\sum_{p\in \ell}r_p-n}{\#\{p\in\ell\}-1}.\]
    The desired statement then follows.
    
\end{proof}

\begin{remark}
    In \cite{YZ23} (or in the proof of \cref{lemma:YZ}), though implicitly, only a weaker version is needed where $\#\{p\in\ell\}-1$ is replaced with $\#\{p\in \ell\}$.
    We emphasize here that this extra $-1$ is crucial in the computation later.
\end{remark}

Recall that we would like to impose some vanishing conditions on a polynomial $f$ of degree at most $n$ to force it to be zero.
Inspired by \cite{TYZ22}, we will do so by visiting each joint for $n+1$ times in some order, and increasing by one the lower bound on multiplicity of $f$ at that joint.
Eventually, we would be able to force $\mult(f,p)>n$ for some joint $p$, which immediately shows that $f$ is zero.
Note that we have the freedom of choosing the order we visit each joint.
The following notation allows us to record a choice of the order.

\begin{definition}
    A function $T:\cJ\times \{0,1,\ldots,n\}\to [(n+1)\abs{\cJ}]$ is a \emph{priority timestamp} if it is bijective and also increasing in the second argument.
\end{definition}

\begin{remark}
    The timestamp induces a total order on the set $\cJ\times \{0,1,\ldots, n\}$, which is reminiscent of the priority order defined in \cite{TYZ22}.
    We denote by the timestamp instead of the total order simply for our own convenience.
\end{remark}

Given a priority timestamp $T$, we define $\halfscript{\vec{v}}{\tT}(t)\in \{0,1\ldots,n+1\}^{\cJ}$ for each $t\in \{0,1,\ldots, (n+1)\abs{\cJ}\}$ so that $v^{\tT}_p(t)$ records the next $r$ such that $T(p,r)>t$, i.e.
\[v^{\tT}_p(t)\eqdef \min\{r\in \ZZ_{\geq 0}: T(p,r)>t\}\quad\forall p\in\cJ.\]
Here we set $T(p,n+1)=\infty$ so that $v^{\tT}_p(t)=n+1$ if $T(p,n)\leq t$.
This will be the vanishing order of the polynomial at stage $t$ at the joint $p$.
In light of \cref{lemma:BasicVanishing}, we also define
\[v^{\tT}_{\ell}(t)\eqdef \frac{\sum_{p\in\ell}v^{\tT}_p(t)-n}{\#\{p\in\ell\}-1}\quad\forall \ell\in \cL .\]
This records the lower bound of the vanishing order at the line $\ell$ guaranteed by \cref{lemma:BasicVanishing} at stage $t$.
In both notations, if the timestamp $T$ is clear from context, we will drop $T$ from the notation.

For a nonnegtaive integer $n$, a priority timestamp $T$, and a terminating time $t_{\textup{f}}$, let
\[S_p(n,T,t_{\textup{f}})\eqdef\left\{\vec{\alpha}\in \ZZ_{\geq 0}^{d}:t\eqdef T(p,\abs{\vec\alpha})\leq t_{\textup{f}},\,\abs{\vec\alpha}-\alpha_i\geq v_{\ell_{p,i}}(t-1)\,\,\forall i\in [d]\right\}\quad\forall p\in \cJ.\]
Intuitively, this collects all $\vec{\alpha}$ such that right before stage $T(p,\abs{\vec{\alpha}})$, we still cannot show $D^{\vec{\alpha}}_pf(p)=0$ simply using \cref{lemma:BasicVanishing}.
We will drop $t_{\textup{f}}$ if $t_{\textup{f}}=(n+1)\abs{\cJ}$.
We will also drop $n$ and $T$ from the notation if they are clear from context.

\begin{lemma}[General vanishing lemma]\label{lemma:GeneralVanishing}
    Let $n$ be a nonnegative integer.
    Let $(\cJ,\cL)$ be a joints configuration in $\FF^d$ with a priority timestamp $T$.
    Let $t_{\textup{f}}\in \{0,\ldots, (n+1)\abs{\cJ}\}$ be a terminating time.
    Then for any polynomial $f\in \FF[x_1,\ldots, x_d]$ with degree at most $n$, if it is $(S_p(n,T,t_{\textup{f}}))_{p\in\cJ}$-vanishing, then $\mult(f,p)\geq v_p(t_{\textup{f}})$ for any $p\in\cJ$ and $\mult(f,\ell)\geq v_{\ell}(t_{\textup{f}})$ for any $\ell\in\cL$.
\end{lemma}
\begin{proof}
    The proof is a simple induction and an application of \cref{lemma:BasicVanishing}.
    We first note that it suffices to show that $\mult(f,p)\geq v_p(t_{\textup{f}})$ for any $p\in \cJ$ by simply applying \cref{lemma:BasicVanishing}.
    When $t_{\textup{f}}=0$, we have $v_p(t_{\textup{f}})=0$, and so the statement holds vacuously.

    Now suppose that the statement already holds for $t_{\textup{f}}-1$, let $(p,r)\in \cJ\times \{0,\ldots, n\}$ be such that $T(p,r) = t_{\textup{f}}$.
    To show the statement for $t_{\textup{f}}$, note that $S_{p'}(n,T,t_{\textup{f}}-1)\subseteq S_{p'}(n,T,t_{\textup{f}})$ for every $p'\in\cJ$.
    Therefore by the inductive hypothesis, it then suffices to show that $D^{\vec{\alpha}}_pf(p)=0$ for any $\abs{\vec{\alpha}}=r$ as $\vec{v}(t_{\textup{f}}) = \vec{v}(t_{\textup{f}}-1)+\vec{e}_p$.
    
    If $\vec{\alpha}\in S_p(n,T,t_{\textup{f}})$, then $D^{\vec{\alpha}}_pf(p)=0$ follows immediately from the assumption that $f$ is $S_p(n,T,t_{\textup{f}})$-vanishing at $p$.
    Therefore we may assume that $\vec{\alpha}\not\in S_p(n,T,t_{\textup{f}})$.
    As $T(p,\abs{\vec{\alpha}})= t_{\textup{f}}$, this means that there exists $i\in [d]$ with $\abs{\vec{\alpha}}-\alpha_i<v_{\ell_{p,i}}(t_{\textup{f}}-1)$.
    Let $\vec{\alpha}'=\vec{\alpha}-\alpha_i\cdot \vec{e_i}$.
    Then the inductive hypothesis implies that $D^{\vec{\alpha}'}_pf$ vanishes identically on $\ell_{p,i}$, and in particular, $D^{\vec{\alpha}}_pf(p) = D^{\alpha_i\cdot\vec{e_i}}_pD^{\vec{\alpha}'}_pf(p)=0$.
    This completes the induction.
\end{proof}

By parameter counting, we immediately get the following inequality.

\begin{corollary}\label{corollary:ParamCounting}
    Let $n$ be a nonnegative integer.
    Let $(\cJ,\cL)$ be a joints configuration in $\FF^d$ with a priority timestamp $T$.
    Then
    \[\sum_{p\in\cJ}\abs{S_p(n,T)}\geq \binom{n+d}{d}\]
\end{corollary}
\begin{proof}
    By \cref{lemma:GeneralVanishing} (where $t_{\textup{f}}=(n+1)\abs{\cJ}$), we see that if $f$ is $(S_p(n,T))_{p\in\cJ}$-vanishing, then $\mult(f,p)\geq n+1$ for any $p\in\cJ$.
    In particular, $f$ must be zero.
    The condition of being $(S_p(n,T))_{p\in\cJ}$-vanishing is a system of $\sum_{p\in \cJ}\abs{S_p(n,T)}$ linear equations on the coefficients of $f$.
    As those equations suffice to force $f$ to be zero, the number of linear equations must be at least the dimension of the vector space $\{f\in \FF[x_1,\ldots, x_d]:\deg f\leq n\}$, which is $\binom{n+d}{d}$.
    The statement then follows.
\end{proof}

Having deduced the general vanishing lemma and the corresponding inequality obtained by parameter counting, we will now specialize to priority timestamps of a certain type.
Fix an arbitrary total order on $\cJ$.
Previously in \cite{TYZ22} and implicitly in \cite{YZ23}, the priority timestamp corresponding to the priority order is the ranking function of $(p,r)\in\cJ\times\ZZ_{\geq 0}$ based on the lexicographical order of $(r-\alpha_p,p)$ for some given ``handicap'' $\vec{\alpha}\in\ZZ_{\geq 0}^{\cJ}$.
One of the new ingredients is a different cleverly chosen priority timestamp that makes the computation easier.

Let $\vec{z}\in \RR_{>0}^{\cJ}$, and let $T_{\vec{z}}$ be the ranking function of $(p,r)\in \cJ\times \ZZ_{\geq 0}$ based on the lexicographical order of $((r-n)/z_p, p)$.
That is, we order the pairs by the value of $(r-n)/z_p$ and break the tie with a prescribed order on the joints if there is a tie.
We call such priority timestamp \emph{the associated priority timestamp of $\vec{z}$}.
In this case, we can compute the size of $\abs{S_p(n,T_{\vec{z}}, t_{\textup{f}})}$ asymptotically by a normalization described in the following paragraph.

In the following, the $O(1)$ and $o(1)$ factors may depend on the joints configuration and the parameter $\vec{z}$, but they do not depend on $n$.
Consider the normalized vector $\vec{a} = \vec{\alpha}/n$, and suppose that $t_{\textup{f}}$ is chosen for each $n$ so that $v_p(t_{\textup{f}}) = (r_p+o(1))n$ for any $p\in\cJ$ for some $\vec{r}\in [0,1]^{\cJ}$.
In this case, if we set $t=T(p,\abs{\vec{\alpha}})$ and compute $v_{p'}(t)$ for any $p'\in \cJ$, then 
\[v_{p'}(t) = (\abs{\vec\alpha}-n)\frac{z_{p'}}{z_p}+n+O(1)=\left(\left(\abs{\vec{a}}-1\right)\frac{z_{p'}}{z_p}+1+o(1)\right)n.\]
Plugging this in, the condition that $\vec{\alpha}\in S_p(n,T_{\vec{z}},t_{\textup{f}})$ can be reformulated as $\abs{\vec{a}}\leq r_p+o(1)$ and 
\[\abs{\vec{a}}-a_i\geq \frac{\sum_{p'\in \ell_{p,i}}z_{p'}}{\left(\#\{p'\in\ell_{p,i}\}-1\right)z_p}\cdot\abs{\vec{a}}+\left(1-\frac{\sum_{p'\in \ell_{p,i}}z_{p'}}{\left(\#\{p'\in\ell_{p,i}\}-1\right)z_p}\right)+o(1)\]
by plugging in the definition of $v_{\ell_{p,i}}(t)$.
To simplify, given $\vec{z}\in \ZZ_{\geq 0}^{\cJ}$, for any $p\in\cJ$ and $i\in [d]$, let
\[\beta_{p,i}\eqdef 1-\left(\#\{p'\in \ell_{p,i}\}-1\right)\frac{z_p}{\sum_{p'\in \ell_{p,i}}z_{p'}}.\]
The second condition can now be written as
\[\abs{\vec{a}}-a_i\geq \frac{1}{1-\beta_{p,i}}\abs{\vec{a}}-\frac{\beta_{p,i}}{1-\beta_{p,i}}+o(1),\]
or equivalently,
\[a_i+\beta_{p,i}(\abs{\vec{a}}-a_i)\leq \beta_{p,i}+o(1).\]
Motivated by this, for any $\beta_1,\ldots, \beta_d,r\leq 1$, define
\[\cS(\beta_1,\ldots, \beta_d,r)\eqdef \left\{\vec{a}\in \RR_{\geq 0}^{d}: \abs{\vec{a}}\leq r,\, a_i+\beta_i(\abs{\vec{a}}-a_i)\leq \beta_i \,\,\forall i\in [d]\right\}.\]
We will drop $r$ from the notation if $r=1$.
Then in fact, the volume of $\cS(\beta_{p,1},\ldots, \beta_{p,d},r_p)$ determines how fast $\abs{S_p(n,T,t_{\textup{f}})}$ grows as $n$ tends to infinity.

\begin{lemma}[Continuous approximation of $S_p$ with $\cS$]\label{lemma:ContiApprox}
    Given a joints configuration $(\cJ,\cL)$ in $\FF^d$, let $\vec{z}\in \RR^{\cJ}_{>0}$ be a vector and $T_{\vec{z}}$ be its associated priority timestamp.
    Define $\beta_{p,i}$ as above for any $p\in\cJ$ and $i\in [d]$.
    Suppose that $\vec{r}\in [0,1]^{\cJ}$ and $0\leq t_{\textup{f}}\leq (n+1)\abs{\cJ}$ are such that $v_p(t_{\textup{f}}) = (r_p+o(1))n$ for any joint $p$, then
    \[\abs{S_p(n,T_{\vec{z}},t_{\textup{f}})} = \left(\vol\cS(\beta_{p,1},\ldots, \beta_{p,d}, r_p)+o(1)\right)\cdot n^d\]
    for any $p\in \cJ$.
\end{lemma}
\begin{proof}
    We fix $p\in \cJ$ throughout the proof.
    It suffices to show that for any $\varepsilon>0$, there exists a sufficiently large $N$ such that for any $n\geq N$, we have $\abs{\abs{S_p(n,T_{\vec{z}}, t_{\textup{f}})}-\vol\cS(\beta_{p,1},\ldots, \beta_{p,d},r_p)}\leq \varepsilon n^d$.
    Let $\vec{\alpha}\in\ZZ_{\geq 0}^d$ and $\vec{a}=\vec{\alpha}/n$.
    From the computation above, there exists $\delta(N)>0$ that tends to $0$ when $N$ goes to infinity that satisfies the following for all $n\geq N$:
    \begin{enumerate}
        \item If $\abs{\vec{a}}\leq r_p-\delta(N)$ and $a_i+\beta_i(\abs{\vec{a}}-a_i)\leq \beta_{p,i}-\delta(N)$ for all $i\in [d]$, then $\vec{\alpha}\in S_p(n,T_{\vec{z}},t_{\textup{f}})$.
        \item If $\vec{\alpha}\in S_p(n,T_{\vec{z}},T_{\textup{f}})$, then $\abs{\vec{a}}\leq r_p+\delta(N)$ and $a_i+\beta_i(\abs{\vec{a}}-a_i)\leq \beta_{p,i}+\delta(N)$ for all $i\in [d]$.
    \end{enumerate}
    Note that $\delta$ may depend on the joints configuration and also $\vec{z}$.
    We make a temporary notation 
    \[\cS^{\pm}(N)\eqdef \left\{\vec{a}\in \RR_{\geq 0}^d: \abs{\vec{a}}\leq r\pm \delta(N),\, a_i+\beta_i(\abs{\vec{a}}-a_i)\leq \beta_i\pm \delta(N)\,\forall i\in [d]\right\}.\]
    Then from the observation above, we have
    \[n\cS^-(N)\cap \ZZ_{\geq 0}^{d} \subseteq S_p(n,T_{\vec{z}},t_{\textup{f}})\subseteq n\cS^+(N)\cap \ZZ_{\geq 0}^{d}\]
    for any $n\geq N$.
    Now we pick $N_0$ sufficiently large so that $\abs{\vol \cS^{\pm}-\vol \cS(\beta_1,\ldots,\beta_d,r)}\leq \varepsilon/2$, and we pick $N\geq N_0$ sufficiently large so that $\abs{\abs{n\cS^{\pm}(N)\cap \ZZ_{\geq 0}^d}-n^d\vol \cS^{\pm}}\leq \varepsilon n^d/2$ for any $n\geq N$.
    The desired statement follows immediately from triangle inequality.
\end{proof}

Combining \cref{corollary:ParamCounting} and \cref{lemma:ContiApprox} and noting that $\binom{n+d}{d} = (\frac{1}{d!}+o(1))n^d$, we get the following corollary.

\begin{corollary}\label{corollary:VanishingIneq}
    Let $(\cJ,\cL)$ be a joints configuration in $\FF^d$ with a vector $\vec{z}\in \RR_{>0}^{\cJ}$.
    Compute $\beta_{p,i}$ accordingly for any $p\in \cJ$ and $i\in [d]$ as above.
    Then
    \[\sum_{p\in\cJ}\vol\cS(\beta_{p,1},\ldots, \beta_{p,d})\geq \frac{1}{d!}.\]
\end{corollary}

We will show that the volume $\vol\cS(\beta_{p,1},\ldots, \beta_{p,d})$ satisfies a good concavity property, and this will help us to estimate the left hand side of the inequality in \cref{corollary:VanishingIneq}. To be more specific, we have the following lemma.

\begin{lemma}\label{lemma:Vol}
Given real numbers $\beta_1,\dots,\beta_d\leq 1$. Let $\avgbeta=(\beta_{1}+\dots+\beta_{d})/d$. We have
    \begin{align*}
        \vol\cS(\beta_{1},\ldots, \beta_{d})\leq\vol\cS(\avgbeta,\ldots, \avgbeta)= \begin{cases}
        \frac{1}{d!\binom{\avgbeta^{-1}+d-1}{d}}&\text{ if $\avgbeta>0$},\\
        0 &\text{ if $\avgbeta\leq 0$}.\end{cases}
    \end{align*}
    When $\avgbeta>0$, the equality of the first part holds if and only if all the $\beta_i$ are the same. 
\end{lemma}

We will postpone the proof of this lemma and also the proof of the next lemma to the following subsections since these two lemmas are rather technical. 

Since we want to count the number of joints, we want to make the volume $\vol\cS(\beta_{p,1},\dots,\beta_{p,d})$ the same for all $p\in\cJ$. In \cite{YZ23}, they also considered the system of equations that the numbers of linear conditions at all the joints are the same. However, in both their proof and our case, it is easier to solve the equation for $\avgbeta_p\eqdef(\beta_{p,1}+\dots+\beta_{p,d})/d$ being the same for all $p\in\cJ$. In their case, the system becomes a linear system rather than a complicated system. In our case, there is no simple formula for $\cS(\beta_{p,1},\dots,\beta_{p,d})$, so we solve the easier system instead.

We know that 
\begin{align*}
 \sum_{p\in \cJ,i\in [d]} \beta_{p,i}=&\sum_{p\in\cJ}\sum_{\ell\ni p}\left(1-\left(\#\{p'\in \ell\}-1\right)\frac{z_p}{\sum_{p'\in \ell} z_{p'}}\right)\\
 =&\sum_{\ell\in \cL}\left(\#\{p\in \ell\}-\left(\#\{p'\in \ell\}-1\right)\frac{\sum_{p\in\ell}z_p}{\sum_{p'\in \ell} z_{p'}}\right)\\
 =& \sum_{\ell\in \cL} 1=L.
\end{align*}
Thus, by this equation and the definition of $\beta_{p,i}$, we shall solve the system of equations 
\[\sum_{i=1}^d \beta_{p,i}=\sum_{i=1}^d \left(1-\left(\#\{p'\in \ell_{p,i}\}-1\right)\frac{z_p}{\sum_{p'\in \ell_{p,i}}z_{p'}}\right)=\frac{L}{J}~\forall p\in\cJ.\]
Let $b_{p,i}\eqdef\frac{z_p}{\sum_{p'\in \ell_{p,i}}z_{p'}}$. We can write the system as
\[\sigma_p\eqdef\sum_{i=1}^d \left(\#\{p'\in \ell_{p,i}\}-1\right)b_{p,i}=\frac{dJ-L}{J}~\forall p\in\cJ.\]
We wish we could solve this equation on all the joints configurations. Unfortunately, we can only solve it on some configurations.
\begin{definition}
A joints configuration $(\cJ,\cL)$ is \emph{critical} if for any joints configuration $(\cJ',\cL')$ with $J'<J$, we have $J'/L'<J/L$.
\end{definition}
\begin{lemma}\label{lemma:SolveEq}
    If $(\cJ,\cL)$ is a critical joints configuration, then the system  
    \begin{align*}
    \sigma_p=\frac{dJ-L}{J}~\forall p\in\cJ
    \end{align*}
    has a solution $\vec{z}$ with $z_p>0$ for all $p\in\cJ$.
\end{lemma}

We will first show how to prove \cref{theorem:SharpJoints} using \cref{lemma:Vol} and \cref{lemma:SolveEq}.
After that, we will give the proofs of the two technical lemmas.

\begin{proof}[Proof of \cref{theorem:SharpJoints} assuming \cref{lemma:Vol} and \cref{lemma:SolveEq}]
Assume first that $(\cJ,\cL)$ is critical. \\
From \cref{lemma:SolveEq}, we know that there is a solution $\vec{z}$ with $z_p>0$ for all $p\in\cJ$ such that $\avgbeta_p=L/dJ$. Combining \cref{corollary:VanishingIneq} and \cref{lemma:Vol}, we have 
\[\frac{J}{d!\binom{dJ/L+d-1}{d}}\geq\frac{1}{d!}.\]
That is, $J\geq \binom{dJ/L+d-1}{d}$. If $J=\binom{x}{d}$, then $x\geq dJ/L+d-1$, and hence $L\geq dJ/(x-d+1)=\binom{x}{d-1}$.

If $(\cJ,\cL)$ is not critical, then by the definition we can always find a critical configuration $(\cJ',\cL')$ with fewer joints and $J/L\leq J'/L'$. Say $J'=\binom{y}{d}$. By the previous part, we know that $L'\geq \binom{y}{d-1}$. Therefore, $J/L\leq J'/L'=(y-d+1)/d$. Since $J'<J$, we know that $y<x$. It follows that $J/L< (x-d+1)/d$, and hence $L>\binom{x}{d-1}$.
\end{proof}
\subsection{Proof of \cref{lemma:Vol}}
\begin{proof}
    We first prove the inequality part. Notice that, if $\beta_i\leq 0$ for some $i$, the left hand side is zero and hence the statement holds. 
    Thus, we may assume that $\beta_1,\dots,\beta_d >0$. We will bound the volume by bounding the length of the intersection of the polytope with lines in certain direction.
    \begin{claim}\label{claim:VolumeSegment}
        Fix any $\beta_3,\dots,\beta_d >0$. For any  $\beta_1,\beta_2,\gamma_3,\dots,\gamma_d\in [0,1]$, let $\eta (\beta_1,\beta_2,\gamma_3,\dots,\gamma_d)$ be the (two-dimensional) area of the polygon $\cS(\beta_{1},\ldots, \beta_{d})\cap \{(a_1,\dots,a_d)\in\RR_{\geq 0}^d\mid a_3=\gamma_3,\dots,a_d=\gamma_d,\}$. It follows that, conditioning on $\beta_1+\beta_2$ being fixed, the area $\eta(\beta_1,\beta_2,\gamma_3,\dots,\gamma_d)$ is maximized when $\beta_1=\beta_2$. Moreover, if $\gamma_3=\dots=\gamma_d=0$ and $\beta_1+\beta_2>0$, the maximum value is only attended at $\beta_1=\beta_2$.
    \end{claim}
    \begin{proof}
        For convenience, we will focus on the plane $a_3=\gamma_3,\dots,a_d=\gamma_d$ and drop the $i$-th coordinate for $i\geq 3$. Also, we will write $\eta(\beta_1,\beta_2)$ instead of $\eta(\beta_1,\beta_2,\gamma_3,\dots,\gamma_d)$ through out the proof of this claim. The polygon is the set bounded by the lines 
        \begin{align}
        a_1,a_2\geq &0,\nonumber\\
        a_1+\beta_1(a_2+\gamma)\leq &\beta_1,\label{Ineq:polygon2}\\
        a_2+\beta_2(a_1+\gamma)\leq &\beta_2,\label{Ineq:polygon3}\\
        \gamma_i+\beta_i(a_1+a_2+\gamma-\gamma_i)\leq & \beta_i,~i=3,\dots,d,\label{Ineq:polygon4}
        \end{align}
        where $\gamma=\gamma_3+\dots+\gamma_d$.
        For each $i=3,\dots,d$, we can rewrite inequality \cref{Ineq:polygon4} as $a_1+a_2\leq r_i$ for some $r_i\in\RR$. Assume the intersection of the equations given by \cref{Ineq:polygon2,Ineq:polygon3} is $(a^*_1,a^*_2)$ and $r^*=a^*_1+a^*_2$. Then \cref{Ineq:polygon4} can be replaced by
        \begin{align}
            a_1+a_2\leq r,\label{Ineq:polygon5}
        \end{align}
        where $r=\min(r^*,r_3,\dots,r_d)$. 
        We know that $r\leq r^*\leq 1-\gamma$ since $a^*_1+a^*_2+\gamma\leq a^*_1/\beta_1+a^*_2+\gamma\leq 1$. If $r\leq 0$ or $\gamma\geq 1$, then clearly $\eta(\beta_1,\beta_2)=0$ and the claim holds. Thus, we may assume $r> 0$ and $0\leq \gamma <1$. 
        
        Let $\ell_1,\ell_2,\ell_3$ denote the lines of the equality cases of inequalities \cref{Ineq:polygon2,Ineq:polygon3,Ineq:polygon5}. We know that the intersection of $\ell_1$ and $\ell_3$ is $(\frac{\beta_1(1-r-\gamma)}{1-\beta_1},\frac{r-\beta_1(1-\gamma)}{1-\beta_1})$, the intersection between $\ell_1$ and $\{a_2=0\}$ is $(\beta_1(1-\gamma),0)$, and the intersection between $\ell_3$ and $\{a_2=0\}$ is $(r,0)$. Thus, the volume of the area satisfying \cref{Ineq:polygon5}, $a_1,a_2\geq 0$, but not $\cref{Ineq:polygon2}$ is 
        \[\frac{\bigl(r-\beta_1(1-\gamma)\bigr)^{+2}}{2(1-\beta_1)}.\]
        where $x^{+k}\eqdef\max(0,x)^k$. Note that since $r\leq r^*$, this area will not overlap with its counterpart obtained by switching the role $a_1,a_2$ and $\beta_1,\beta_2$. 
        It follows that
        \begin{align*}
            \eta(\beta_1,\beta_2)=\frac{r^2}{2}-\frac{\bigl(r-\beta_1(1-\gamma)\bigr)^{+2}}{2(1-\beta_1)}-\frac{\bigl(r-\beta_2(1-\gamma)\bigr)^{+2}}{2(1-\beta_2)}.
        \end{align*}
        Note that similarly we have
        \begin{align*}
            \eta\left(\frac{\beta_1+\beta_2}{2},\frac{\beta_1+\beta_2}{2}\right)\geq \frac{r^2}{2}-\frac{\bigl(r-\frac{\beta_1+\beta_2}{2}(1-\gamma)\bigr)^{+2}}{(1-\frac{\beta_1+\beta_2}{2})},
        \end{align*}
        since either $a_1+a_2\leq r^*$ cuts some part of the polygon defining $\eta(\frac{\beta_1+\beta_2}{2},\frac{\beta_1+\beta_2}{2})$ or some parts are overlapped and double counted in the negative term.
        
        Thus, it is sufficient to show that $\eta'(\beta)\eqdef\frac{\bigl(r-\beta(1-\gamma)\bigr)^{+2}}{2(1-\beta)}$ is a convex function for $\beta\in [0,1]$ and strictly convex when $\gamma=0$ and $\beta\in [0,1)$. Note that 
        when $\beta<\frac{r}{1-\gamma}$, the first derivative of $\eta'(\beta)$ is
        \[\frac{d}{d\beta}\frac{\bigl(r-\beta(1-\gamma)\bigr)^{2}}{2(1-\beta)}=\frac{(r-\beta(1-\gamma))(r-(2-\beta)(1-\gamma)))}{2(1-\beta)^2}.\]
        This is negative since $\beta\leq 1, r>\beta (1-\gamma)$, and $r\leq 1-\gamma\leq (2-\beta)(1-\gamma)$.
        The second derivative of $\eta'(\beta)$ when $\beta<\frac{r}{1-\gamma}$ is
        \[\frac{d^2}{d\beta^2}\frac{\bigl(r-\beta(1-\gamma)\bigr)^{2}}{2(1-\beta)}=\frac{(r+\gamma-1)^2}{4(1-\beta)^3}> 0.\]
        Thus, $\eta'(\beta)$ is non-negative and continuous, strictly convex and monotonely decreasing on $[0,\frac{r}{1-\gamma}]$, and zero when $\beta\geq \frac{r}{1-\gamma}$, and it is therefore convex on $[0,1]$. 
        
        For the second part of the claim, we have $\gamma=0$. For the sake of contradiction, assume $\beta_1\neq \beta_2$ maximize $\eta(\beta_1,\beta_2)$.
        Then it suffices to prove that $r>\frac{\beta_1+\beta_2}{2}$.
        This is because that $\eta$ is strictly convex on the interval $[0,r)$.
        Note that \cref{Ineq:polygon4} becomes
        \[\beta_i(a_1+a_2)\leq \beta_i\]
        for each $i=3,\ldots, d$, or equivalently,
        \[a_1+a_2\leq 1\]
        for each $i=3,\ldots,d$.
        Hence $r_3=\cdots =r_d= 1$, and it remains to show that $r^*>\frac{\beta_1+\beta_2}{2}$.
        Note that $r^*=a_1^*+a_2^*\geq a_1^*+\beta_1a_2^*=\beta_1$ and similarly $r^*\geq \beta_2$. Thus, $r^*\geq\max(\beta_1,\beta_2)>\frac{\beta_1+\beta_2}{2}$ if $\beta_1\neq \beta_2$.
    \end{proof}
    We know that $\vol\cS(\beta_{1},\ldots, \beta_{d})$ is the integration of $\eta (\beta_1,\beta_2,\gamma_3,\dots,\gamma_d)$ over $\gamma_3,\dots,\gamma_d\in [0,1]$. Thus, by using \cref{claim:VolumeSegment}, we know that $\vol\cS(\beta_{1},\ldots, \beta_{d})\leq \vol\cS(\frac{\beta_{1}+\beta_2}{2},\frac{\beta_{1}+\beta_2}{2},\ldots, \beta_{d})$. Moreover, by looking at the values around $\gamma_3=\dots=\gamma_d=0$, the inequality is strict if $\beta_1\neq \beta_2$ and $\beta_1,\dots,\beta_d>0$. Since $\vol\cS(\beta_{1},\ldots, \beta_{d})$ is symmetric and continuous on the compact region $[0,1]^d$, the first part of the lemma follows.

    It remains to compute the volume $\vol\cS(\avgbeta,\ldots, \avgbeta)$. Note that, if $\avgbeta\leq 0$, the volume is clearly zero. 

    Assume that $\avgbeta>0$. We can partition $\cS(\avgbeta ,\ldots, \avgbeta )$ into $d!$ parts in the following way. For each $\sigma$ in the symmetric group of $[d]$, define 
    \[\cS_\sigma=\cS(\avgbeta ,\ldots, \avgbeta )\cap\{(a_1,\dots,a_d)\in\RR_{\geq 0}^d\mid a_{\sigma(1)}\geq a_{\sigma(2)}\geq \dots\geq a_{\sigma(d)}\}.\]
    It follows that all $\cS_\sigma$ form a partition of $\cS(\avgbeta ,\ldots, \avgbeta )$ (up to boundaries). By symmetry, the volumes $\vol\cS_{\sigma}$ are all the same, so it is enough compute the $\vol\cS_{\id}$, where $\id$ is the identity permutation. Note that $\cS_{\id}$ is the simplex defined by 
    \[a_1\geq \dots\geq a_d\geq 0, a_1+\avgbeta (a_2+\dots+a_d)\leq \avgbeta ,\]
    and hence the vertices of $\cS_{\id}$ are
    \[p_0\eqdef (0,\dots,0)\text{ and } p_i\eqdef\Big(\underbrace{\frac{1}{\avgbeta^{-1} +i-1},\dots,\frac{1}{\avgbeta^{-1} +i-1}}_{i \text{ terms}},0,\dots,0\Big)\text{ for }i\in [d].\]
    Since the height from the subspace spanned by $p_0\dots,p_{i-1}$ to $p_i$ is $\frac{1}{1/\avgbeta +i-1}$, the volume of $\cS_{\id}$ is
    \[\frac{1}{d!}\prod_{i=1}^d\frac{1}{\avgbeta^{-1} +i-1}.\]
    Hence, we have
    \[\vol\cS(\avgbeta ,\ldots, \avgbeta )=d!\vol\cS_{\id}=\prod_{i=1}^d\frac{1}{\avgbeta^{-1} +i-1}=\frac{1}{d!\binom{\avgbeta^{-1} +d-1}{d}}.\]
   
\end{proof}
\subsection{Proof of \cref{lemma:SolveEq}}
Let $S$ be a finite index set.
We say that a function $g:\RR_{>0}^S\to \RR^S$ is \emph{homogeneous} if it is scaling-invariant: for any $r>0$ and any $\vec{z}\in \RR_{>0}^S$, we always have $g(r\vec{z}) = g(\vec{z})$.
The following statement and proof are inspired by Nachmias \cite[Chapter 3.2]{Nach-book}, where he proved the circle packing theorem using an iterative approximation algorithm.

\begin{lemma}\label{lemma:EquationTechnical}
    Let $S$ be a finite set, and let $g:\RR_{>0}^S\to \RR^S$ be continuous with $g(\vec{z})=\nobreak(g_s(\vec{z}))_{s\in S}$ such that $\sum_{s\in S}g_s(\vec{z})$ does not depend on $\vec{z}\in \RR_{>0}$.
    Suppose furthermore that for any nonempty proper $S'\subsetneq S$ and continuous $\vec{z}:(0,1]\to \RR_{>0}^S$ with
    \[z_s(r)=\begin{cases}
        rz_s(1)& \textup{if }s\in S',\\
        z_s(1)& \textup{if } s\not\in S',
    \end{cases}\]
    we have the following:
    \begin{enumerate}
        \item The function $g_s\circ \vec{z}$ is increasing if $s\in S'$, and decreasing if $s\not\in S'$.
        \item We have
        \[\max_{s\not\in S'}\lim_{r\to 0}g_s(\vec{z}(r)) > \min_{s\in S'}\lim_{r\to 0}g_s(\vec{z}(r)).\]
    \end{enumerate}
    Then there exists an infinite sequence $\vec{z}^{(1)},\vec{z}^{(2)},\ldots\in \RR_{>0}^S$ such that $\lim_{n\to\infty}g_s(\vec{z}^{(n)})$ exists for any $s\in S$, and the limits are all equal.
\end{lemma}
\begin{proof}
    We will construct the sequence inductively.
    We start with an arbitrary $\vec{z}^{(1)}$.
    Given $\vec{z}^{(n)}$, we list out all the values of $g_s(\vec{z}^{(n)})\, (s\in S)$ on the real line.
    If they are all the same, then we are already done by simply setting $\vec{z}^{(n+1)} = \vec{z}$.
    Otherwise, we look at the biggest gap between the $\abs{S}$ points on the real line (if there is tie, choose any among the biggest ones).
    Collects all the $s\in S$ with $g_s(\vec{z}^{(n)})$ being above the gap into a set $A^{(n)}$.
    Now, with a slight abuse of notation, let $\vec{z}^{(n)}$ be a continuous function where $\vec{z}^{(n)}(1)$ is the original $\vec{z}^{(n)}$, and
    \[z^{(n)}_s(r)=\begin{cases}
        rz^{(n)}_s(1)& \textup{if }s\in A^{(n)},\\
        z^{(n)}_s(1)& \textup{if } s\not\in A^{(n)}.
    \end{cases}\]
    By the continuity of $g\circ \vec{z}^{(n)}$, the conditions and the intermediate value theorem, there exists $r^{(n)}\in (0,1)$ such that $\max_{s\not\in A^{(n)}}g_s(\vec{z}^{(n)}(r^{(n)})) = \min_{s\in A^{(n)}}g_s(\vec{z}^{(n)}(r^{(n)}))\eqdef x^{(n)}$.
    We then set $\vec{z}^{(n+1)} = \vec{z}^{(n)}(r^{(n)})$.
    We also suppose that $s_B^{(n)}\not\in A^{(n)}$ maximizes $g_s(\vec{z}^{(n)}(r))$ for $s\not\in A^{(n)}$ and $s_A^{(n)}\in A^{(n)}$ minimizes $g_s(\vec{z}^{(n)}(r^{(n)}))$ for $s\in A^{(n)}$.

    To analyze how close the values of $g_s(\vec{z}^{(n)})\, (s\in S)$ are to one another, let $\gamma=\sum_{s\in S}g(\vec{z})/\abs{S}$ for any $\vec{z}\in \RR^S_{>0}$ and $\cE^{(n)} = \sum_{s\in S}\left(g_s(\vec{z}^{(n)})-\gamma\right)^2$.
    Recall from the assumption that $\gamma$ does not depend on the choice of $\vec{z}$.
    We will show that $\lim_{n\to\infty}\cE^{(n)}=0$ by comparing $\cE^{(n+1)}$ to $\cE^{(n)}$.

    Let $\delta^{(n)}$ be the gap between $\max_{s\not\in A^{(n)}}g_s(\vec{z}^{(n)}(1))$ and $\min_{s\in A^{(n)}}g_s(\vec{z}^{(n)}(1))$.
    We first rewrite
    \[\cE^{(n)}=\sum_{s\in S}\left(g_s(\vec{z}^{(n)}(1))-x^{(n)}\right)^2-\abs{S}\left(x^{(n)}-\gamma\right)^2\]
    as $\sum_{s\in S}g_s(\vec{z}^{(n)}(1))=\gamma\abs{S}$.
    Similarly
    \[\cE^{(n+1)} = \sum_{s\in S}\left(g_s(\vec{z}^{(n)}(r^{(n)}))-x^{(n)}\right)^2-\abs{S}\left(x^{(n)}-\gamma\right)^2.\]
    Note that by definition and the monotonicity of $g_s\circ \vec{z}^{(n)}$ for each $s\in S$, we have $s\in A^{(n)}$ if and only if the value $g_s(\vec{z}^{(n)}(1))$ is above $x^{(n)}$.
    Again, by the monotonicity, we know that $\left(g_s(\vec{z}^{(n)}(1))-x^{(n)}\right)^2\geq \left(g_s(\vec{z}^{(n)}(r^{(n)}))-x^{(n)}\right)^2$.
    Therefore
    \[\cE^{(n)}-\cE^{(n+1)}\geq \left(g_{s_A^{(n)}}(\vec{z}^{(n)}(1))-x^{(n)}\right)^2+\left(g_{s_B^{(n)}}(\vec{z}^{(n)}(1))-x^{(n)}\right)^2\geq \frac{1}{4}\left(\delta^{(n)}\right)^2\]
    as $g_{s_A^{(n)}}(\vec{z}^{(n)}(1))-g_{s_B^{(n)}}(\vec{z}^{(n)}(1))\geq \delta^{(n)}$.
    However, as $\delta^{(n)}$ is the largest gap between $g_s(\vec{z}^{(n)}(1))\,(s\in S)$, we can simply bound $\abs{g_s(\vec{z}^{(n)}(1))-\gamma}$ by $\abs{S}\delta^{(n)}$.
    Thus,
    \[\cE^{(n)}\leq \abs{S}^3\left(\delta^{(n)}\right)^2,\]
    and so
    \[\cE^{(n+1)}\leq \left(1-\frac{1}{4\abs{S}^3}\right)\cE^{(n)}.\]
    This shows that $\cE^{(n)}$ tends to zero exponentially as $n$ tends to infinity.
    The statement then follows immediately.
    
\end{proof}

\begin{lemma}\label{lemma:EquationTechnical2}
    Let $(\cJ,\cL)$ be a critical joints configuration.
    Then the function $g$ who takes $\vec{z}\in \RR_{>0}^{\cJ}$ as input and outputs $\vec{\sigma}$ as in \cref{lemma:SolveEq} satisfies all the conditions in \cref{lemma:EquationTechnical}.
\end{lemma}
\begin{proof}
    First of all, as $b_{p,\ell}$ is homogeneous for any $p\in \ell$, it is clear that $\vec{\sigma}$ is homogeneous.
    It remains to check the two conditions (1) and (2).
    
    Fix any nonempty proper subset $\cJ'\subsetneq \cJ$.
    For each $p\in\cJ'$ and $p\in \ell\in \cL$, we have 
    \[b_{p,\ell}(\vec{z}(r)) = \frac{rz_p(1)}{r\sum_{\substack{p'\in \ell\\ p'\in \cJ'}}z_{p'}(1)+\sum_{\substack{p'\in \ell \\ p'\not\in \cJ'}}z_{p'}(1)}=\frac{z_p(1)}{\sum_{\substack{p'\in \ell\\ p'\in \cJ'}}z_{p'}(1)+r^{-1}\sum_{\substack{p'\in \ell \\ p'\not\in \cJ'}}z_{p'}(1)},\]
    which is increasing in $r$.
    Therefore $\sigma_p$ is also increasing in $r$.
    On the other hand, if $p\not\in \cJ'$, then
    \[b_{p,\ell}(\vec{z}(r)) = \frac{z_p(1)}{r\sum_{\substack{p'\in \ell\\ p'\in \cJ'}}z_{p'}(1)+\sum_{\substack{p'\in \ell \\ p'\not\in \cJ'}}z_{p'}(1)},\]
    which is decreasing in $r$.
    Thus, condition (1) is satisfied.
    Now if (2) fails, we extend the definition of $\vec{z}$ continuously to $0$.
    We also extend $b_{p,\ell}\circ \vec{z}$ continuously to $0$ as follows.
    If $p\not\in \cJ'$ then it is clear that 
    \[\lim_{r\to 0}b_{p,\ell}(\vec{z}(r))=\frac{z_p(1)}{\sum_{\substack{p'\in \ell\\ p'\not\in \cJ'}}z_{p'}(1)}.\]
    If $p\in \cJ'$ and $\ell$ contains some other joint that is not in $\cJ'$, then $\lim_{r\to 0}b_{p,\ell}(\vec{z}(r))=0$.
    Otherwise, we have that $b_{p,\ell}(\vec{z}(r))$ is constant, and thus $b_{p,\ell}(\vec{z}(0))=b_{p,\ell}(\vec{z}(1))$.
    
    Let $\cL'$ be the set of lines that only contains joints in $\cJ'$.
    Remove all lines not in $\cL'$, and for every $(p,\ell)\in \cJ'\times (\cL\backslash \cL')$ with $p\in \ell$, add a dummy line that only passes through $p$ in general position.
    Let $\Tilde{\cL}$ be the set of dummy lines.
    Then $(\cJ',\cL'\cup \Tilde{\cL})$ is a joints configuration.
    Therefore
    \begin{align*}
        dJ'-(L'+\Tilde{L})=&\sum_{p\in \cJ'}\#\{\ell\in \cL':p\in\ell\}-L'\\
        =&\sum_{\ell\in\cL'}\left(\#\{p'\in\cJ':p'\in\ell\}-1\right)\\
        =&\sum_{\ell\in \cL'}\sum_{p\in \cJ'\cap \ell}\left(\#\{p'\in \cJ': p'\in \ell\}-1\right)\frac{z_p(1)}{\sum_{p'\in \cJ'\cap \ell}z_{p'}(1)}\\
        =&\sum_{p\in \cJ'}\sum_{\substack{\ell\in \cL'\\
        \ell\ni p}}\left(\#\{p'\in \cJ': p'\in \ell\}-1\right)b_{p,\ell}(\vec{z}(1))\\
        =&\sum_{p\in \cJ'}\sigma_p(\vec{z}(0))\\
        \geq& \frac{J'}{J}\sum_{p\in\cJ}\sigma_p(\vec{z}(0))\\
        =&\frac{J'}{J}\cdot \left(dJ-L\right).
    \end{align*}
    After rearranging, we get $J'/(L'+\Tilde{L})\geq J/L$.
    Since $J'<J$, this is a contradiction with the assumption that $(\cJ,\cL)$ is critical.
\end{proof}

\begin{proof}[Proof of \cref{lemma:SolveEq}]
    Since $(\cJ,\cL)$ is critical, by \cref{lemma:EquationTechnical2} and \cref{lemma:EquationTechnical}, there exists an infinite sequence $\vec{z}^{(1)},\vec{z}^{(2)},\ldots \in\RR^{\cJ}_{>0}$ such that 
    \[\lim_{n\to \infty}\sigma_p(\vec{z}^{(n)})=\frac{dJ-L}{J}\]
    for all $p\in\cJ$.
    We would like to take the limit of $\vec{z}^{(n)}$, but as $\RR^{\cJ}_{>0}$ is not compact, there is not necessarily a subsequence that is convergent.
    However, as $\vec{\sigma}$ is homogeneous, we may focus on the ratios between the components of $\vec{z}^{(n)}$ instead.

    To record the ratios, for each $\ell\in \cL$ let
    \[s_{\ell}^{(n)}\eqdef \frac{\sum_{p\in \ell}z_p^{(n)}}{\sum_{p\in \cJ}z_p^{(n)}},\]
    and for every $(p,\ell)\in \cJ\times \cL$ with $p\in \ell$, let
    \[b_{p,\ell}^{(n)}\eqdef \frac{z_p^{(n)}}{\sum_{p\in \ell}z_p^{(n)}}.\]
    It is clear that $\sum_{\ell\in\cL}s_{\ell}^{(n)}=d$ and $\sum_{p\in\ell}b_{p,\ell}^{(n)}=1$ for all fixed $\ell$.
    Therefore, we may pass to a subsequence so that after relabeling, the limits $\lim_{n\to\infty} s_{\ell}^{(n)}$ and $\lim_{n\to\infty} b_{p,\ell}^{(n)}$ exist for all $p\in\ell$.
    We will denote the limits by $s_{\ell}^{(\infty)}$ and $b_{p,\ell}^{(\infty)}$.

    Let $\ell$ and $\ell'$ be two lines that pass through $p$.
    Then it is clear that $b_{p,\ell}^{(n)}s_{\ell}^{(n)} = b_{p,\ell'}^{(n)}s_{\ell'}^{(n)}$, and thus $b_{p,\ell}^{(\infty)}s_{\ell}^{(\infty)}=b_{p,\ell'}^{(\infty)}s_{\ell'}^{(\infty)}$.
    We will first show that $s_{\ell}^{(\infty)}>0$ for all $\ell\in\cL$.
    Otherwise, let $\cJ'$ be the set of joints with $s_{\ell}^{(\infty)}=0$ for some $\ell\ni p$.
    Then $\cJ'$ is nonempty.
    It is also proper: take $\ell^+$ to be some $\ell$ with $s_{\ell}^{(\infty)}>0$, and on $\ell^+$ take $p^+$ to be some $p$ with $b_{p,\ell^+}^{(\infty)}>0$.
    Then for any $\ell\ni p^+$, we have $b_{p^+,\ell}^{(\infty)}s_{\ell}^{(\infty)} = b_{p^+,\ell^+}^{(\infty)}s_{\ell^+}^{(\infty)}>0$, and so $s_{\ell}^{(\infty)}\neq 0$.
    This shows that $p^+\not\in \cJ'$.

    Let $\cL'$ be the set of lines with $s_{\ell}^{(\infty)}=0$.
    Then for any $\ell\not\in \cL'$ that intersects $\ell'\in\cL'$ at $p$, we have that $b_{p,\ell}^{(\infty)} = b_{p,\ell'}^{(\infty)}s_{\ell'}^{(\infty)}/s_{\ell}^{(\infty)}=0$.
    As in the proof of \cref{lemma:EquationTechnical2}, we remove all lines not in $\cL'$, and for every $p\in \cJ'$ lying on $\ell\in\cL\backslash \cL'$, we add a line through $p$ in a general direction so that $p$ remains a joint.
    Let $\Tilde{\cL}$ be the set of added lines, making $(\cJ',\cL'\cup\Tilde{\cL})$ a joints configuration.
    Then, similar to the argument used in \cref{lemma:EquationTechnical2},
    \begin{align*}
        dJ'-(L'+\Tilde{L})&=\sum_{\ell\in \cL'\cup \Tilde{\cL}}(\#\{p'\in \cJ': p'\in \ell\}-1)\\
        &=\sum_{\ell\in \cL'}\sum_{p\in\cJ'\cap \ell}\left(\#\{p'\in \cJ': p'\in \ell\}-1\right)b_{p,\ell}^{(\infty)}\\
        &=\sum_{p\in \cJ'}\sum_{\substack{\ell\ni p\\ \ell\in\cL'}}\left(\#\{p'\in \cJ': p'\in \ell\}-1\right)b_{p,\ell}^{(\infty)}\\
        &=\sum_{p\in\cJ'}\sigma_p^{(\infty)}\\
        &=\frac{J'}{J}\cdot (dJ-L).
    \end{align*}
    Therefore $J'/(L'+\Tilde{L})=J/L$.
    Note that $J'<J$, and so we reach a contradiction with the assumption that $(\cJ,\cL)$ is critical.
    As a consequence, $s_{\ell}^{(\infty)}>0$ for all $\ell\in\cL$.

    To wrap up, we can simply take, for all $p\in\cJ$, $z_p^{(\infty)} = b_{p,\ell}^{(\infty)}s_{\ell}^{(\infty)}$ for an arbitrary $\ell\ni p$.
    If $z_p^{(\infty)}>0$ for all $p\in\cJ$, then $b_{p,\ell}(\vec{z}^{(\infty)}) = b_{p,\ell}^{(\infty)}$ and so we are done.
    It remains to deal with the case where $z_p^{(\infty)}=0$ for some $p\in\cJ$.
    However, if this occurs, then $b_{p,\ell}^{(\infty)}=0$ for all $\ell\ni p$ as $s_{\ell}^{(\infty)}\neq 0$.
    This shows that $(dJ-L)/J = \sigma_p(\vec{z}^{(\infty)})=0$, and thus $L=dJ$.
    By the criticality of $(\cJ,\cL)$, it must be the case that $J=1, L=d$ and the joints configuration is $d$ linearly independent lines through a single joint.
    In this case we can simply choose $z_p$ to be an arbitrary positive real.
\end{proof}

\section{Structural theorem}\label{section:Structural}
In this section, we will first discuss some idea that will be useful.
We will then prove \cref{theorem:Structural}.
At the end, we briefly discuss what we expect for some stronger structural and stability statements for future research.

\subsection{Motivating idea}
Let $(\cJ,\cL)$ be the joints configuration in \cref{example:tight} once again for the discussion.
Then the solution to the equations in \cref{lemma:SolveEq} is simply taking the $z_p$'s to be the same.
With this $\vec{z}$, we get that $\beta_{p,i} = (M-d+1)^{-1}$.
Therefore $S_p(n,T_{\vec{z}})$ is approximately 
\[B'\eqdef \left\{\vec{\alpha}\in\ZZ_{\geq 0}^d:\alpha_i+\frac{1}{M-d+1}(\abs{\vec{\alpha}}-\alpha_i)<\frac{n}{M-d+1}\,\forall i\in[d]\right\},\]
the shaved box that appeared in \cref{subsection:KeyIdea}.
\cref{lemma:GeneralVanishing} states that if $f$ is a polynomial of degree at most $n$ that is $(S_p(n,T_{\vec{z}}))_{p\in\cJ}$-vanishing, then $f$ is zero.
One can ask: is there any $f$ such that $f$ is not zero, but $f$ is ``almost'' $(S_p(n,T_{\vec{z}}))_{p\in\cJ}$-vanishing?
Indeed, it is possible to construct such $f$ given the structure of $(\cJ,\cL)$.
Let $s_1(x_1,\ldots, x_d),\ldots,s_M(x_1,\ldots, x_d)$ be the $M$ linear functionals on $\FF^d$ that define the $M$ hyperplanes in the construction of $(\cJ,\cL)$.
Suppose that $M\mid n$ for our convenience.
Then we can take $f=(s_1\cdots s_M)^{n/M}$.
Though it is definitely not $(S_p(n,T_{\vec{z}}))_{p\in\cJ}$-vanishing by \cref{lemma:GeneralVanishing}, one can verify that it is indeed $(B')_{p\in\cJ}$-vanishing.
In some sense, we can view $f$ as being almost $(S_p(n,T_{\vec{z}}))_{p\in\cJ}$-vanishing.

Now note that the polynomial $f$ actually reveals the structure of $(\cJ,\cL)$: the linear factors that appear in $f$ exactly determine the $M$ hyperplanes.
Our strategy in general is thus to find polynomials $f$ that are ``almost'' $(S_p(n,T_{\vec{z}}))_{p\in\cJ}$-vanishing and proceed from there.
We will not be able to prove such a strong global statement as $f$ has exactly $M$ distinct linear factors, but in some sense, we will be able to factorize $f$ locally around the joint and utilize this factorization.

\subsection{Main proof}

We call a joints configuration $(\cJ,\cL)$ in $\FF^d$ \emph{tight} if there exists a real number $M\geq d$ such that $\abs{\cJ} = \binom{M}{d}$ and $\abs{\cL} = \binom{M}{d-1}$.

We will show the structural theorem by showing that, for any two joints on the same line, they share a common structure.

\begin{definition}
    For any two joints $p,p'$ on the same line $\ell$, we say that they \emph{share} a hyperplane $H$ if $H$ contains $\ell$, $(d-2)$ other lines through $p$, and $(d-2)$ other lines through $p'$. Moreover, we say that $p$ and $p'$ \emph{share all the hyperplanes} if they share $(d-1)$ distinct hyperplanes.
\end{definition}

The key lemma is the following structural result.

\begin{lemma}\label{lemma:Structural}
    Let $(\cJ,\cL)$ be a tight joints configuration in $\FF^d$.
    For any two joints $p,p'$ on the same line $\ell$, the two joints $p,p'$ share all the hyperplanes.
\end{lemma}

Once this is done, the main structural theorem will follow immediately by the Kruskal--Katona theorem.

\begin{proof}[Proof of \cref{theorem:Structural} assuming \cref{lemma:Structural}]
    Let $X$ be the set of all hyperplanes in $\FF^d$. For each joint $p\in \cJ$, let $P_p\in\binom{X}{d}$ be the set of hyperplanes spanned by $(d-1)$ of the lines $\ell_{p,1},\dots,\ell_{p,d}$. Let $\bJ$ be the family $\{P_p:p\in\cJ\}$. Set $\partial\bJ$ to be the family consists of all $(d-1)$-sets $F\subset P_p$ for some $p$. From \cref{lemma:Structural}, we know that there is a bijection between $\partial\bJ$ and $\cL$. Thus, we may apply the structural theorem of Lovasz's version of the Kruskal--Katona theorem (see \cite[Problem 13.31(b)]{Lovasz93}). The structural theorem is not explicitly written in the problem statement. Thus, we give a brief remark here on the required modification to obtain the structural result.
    
    In problem 13.31(a), the equality holds if and only if $u=v+1,v=w$. Therefore, in the solution of problem 13.31(b), $H_1$ is a complete $r$-uniform graph on $v$ vertices and $H_2$ is a complete $(r-1)$-uniform graph on $v$ vertices by also inducting on the structural statement. Note that the degree of each vertex must be the same since $x$ is chosen arbitrarily. Thus, the only way to combine $H_1$ and $H_2$ is to make $H$ also complete.
    
    Using the structural theorem of Lovasz's version of the Kruskal--Katona theorem , it follows that $M$ is an integer, and there must be a set $X'$ of $M$ hyperplanes such that $\bJ=\binom{X'}{d}$ and $\partial\bJ=\binom{X'}{d-1}$. That is, $\cJ$ is the set of $d$-wise intersections of the hyperplanes in $X'$ and $\cL$ is the set of $(d-1)$-wise intersections of the hyperplanes in $X'$.
\end{proof}

In the rest of this subsection, we will prove \cref{lemma:Structural}.

\begin{lemma}\label{lemma:PrelimStructure}
    Suppose that $(\cJ,\cL)$ is a tight joints configuration in $\FF^d$ where $\abs{\cJ}=\binom{M}{d}$ and $\abs{\cL} = \binom{M}{d-1}$.
    Then $M$ is a positive integer, the joints configuration is connected, and there are exactly $M-d+1$ joints on each line in $\cL$. Here, we say that a joints configuration $(\cJ,\cL)$ is connected if the hypergraph with vertex set $\cJ$ and edge set induced by $\cL$ is connected.
\end{lemma}
\begin{proof}
    From the proof of the main theorem, we know the configuration $(\cJ,\cL)$ must be critical. This implies the configuration is connected. Otherwise, one of its connected component $(\cJ',\cL')$ is a configuration with fewer joints and $J'/L'\geq J/L$, which yields a contradiction.

    Since $(\cJ,\cL)$ is a tight joints configuration, we know the estimate in \cref{lemma:Vol} must also be tight. Therefore, $\beta_{p,i}$ must be the same for all $p\in \cJ$ and $i\in [d]$. For each $\ell\in L$, we have
\begin{align*}
\sum_{(p,i):\ell_{p,i}=\ell}\beta_{p,i}=&\sum_{p\in\ell}\left(1-\left(\#\{p'\in \cJ:p'\in \ell\}-1\right)\frac{z_p}{\sum_{p'\in \ell} z_{p'}}\right)\\
=&\left(\#\{p\in \cJ:p\in \ell\}-\left(\#\{p'\in \cJ:p'\in \ell\}-1\right)\frac{\sum_{p\in\ell}z_p}{\sum_{p'\in \ell} z_{p'}}\right)\\
=&1.
\end{align*}
Therefore, $\beta_{p,i}$ is the inverse of the number of points on $\ell$ for each $\ell\in\cL$. In particular, all the lines contain the same number of joints. 

On the other hand, we know that a line contains $dJ/L=M-d+1$ joints on average. Thus, in this case, all the lines contain exactly $M-d+1$ joints and hence $M$ is an integer.
\end{proof}

As a consequence, the solution to the equations in \cref{lemma:SolveEq} is $\vec{z}$ where $z_p=z_{p'}$ for any $p,p'\in \cJ$.
Its associated priority timestamp $T_{\vec{z}}$ thus visits every joint in order for $n+1$ times.
In particular, with this priority timestamp, we always have $\abs{v_p(t)-v_{p'}(t)}\leq 1$ for any $p,p'\in\cJ$ and $0\leq t\leq (n+1)\abs{\cJ}$.
We will fix this choice of $\vec{z}$ and priority timestamp in the rest of this section.
With this choice, we know that $\beta_{p,i}=1-(M-d)/(M-d+1) = 1/(M-d+1)$ for every $p\in\cJ$ and $i\in[d]$.

To prove \cref{lemma:Structural}, we will use crucially a nonzero polynomial $f$ that vanishes to really high order at each joint.
We will first show that such polynomial always exists, and then we will investigate to which order the polynomial does not vanish.
Those data would eventually help us in showing the structural result.

\begin{lemma}\label{lemma:StructuralPolynomial}
    Fix $\varepsilon \in (0,1)$.
    Let $(\cJ,\cL)$ be a tight joints configuration, and $p\in\cJ$ be a joint.
    For any positive integer $n$, set $t_{\textup{f}}$ so that $v_{p'}(t_{\textup{f}})=\lceil(1-\varepsilon)dn/M\rceil$ for all $p'\in\cJ$.
    Then for all sufficiently large $n$, there exists a nonzero polynomial $f$ of degree at most $n$ that is $S_p(n,T_{\vec{z}},t_{\textup{f}})$-vanishing at $p$ and $S_{p'}(n,T_{\vec{z}})$-vanishing at $p'$ for any $p'\in\cJ\backslash \{p\}$.
\end{lemma}
\begin{proof}
    It suffices to show that 
    \[\abs{S_p(n,T_{\vec{z}},t_{\textup{f}})}+\sum_{p'\in\cJ\backslash \{p\}}\abs{S_{p'}(n,T_{\vec{z}})}<\binom{n+d}{d}\]
    as the left hand side is the number linear conditions required to force $f$ to be $S_p(n,T_{\vec{z}},t_{\textup{f}})$-vanishing at $p$ and $S_{p'}(n,T_{\vec{z}})$-vanishing at $p'$ for any $p'\in\cJ\backslash \{p\}$.
    By \cref{lemma:ContiApprox}, it suffices to show that
    \[\vol\cS\left(\frac{1}{M-d+1},\ldots, \frac{1}{M-d+1},(1-\varepsilon)\frac{d}{M}\right)+(J-1)\vol P<\frac{1}{d!}\]
    where $P\eqdef \cS(1/(M-d+1),\ldots, 1/(M-d+1))$.
    Note that $J = \binom{M}{d}$ and $\vol P= \left(d!\binom{M}{d}\right)^{-1}$ by \cref{lemma:Vol}.
    Therefore it suffices to show that
    \[\vol P -\vol \cS\left(\frac{1}{M-d+1},\ldots, \frac{1}{M-d+1},(1-\varepsilon)\frac{d}{M}\right)>0,\]
    or equivalently, $P\cap \{\vec{a}\in\RR^d_{\geq 0}:\abs{\vec{a}}>(1-\varepsilon)d/M\}$ has positive volume.
    This is true as the apex of the polytope $P$ is $(1/M, \ldots, 1/M)$.
\end{proof}

In the remaining of the section, we will fix this choice of $p,\varepsilon, t_{\textup{f}}$ and a nonzero $f$ that satisfies the vanishing conditions.
We first derive some vanishing and non-vanishing result from \cref{lemma:GeneralVanishing}.

\begin{lemma}\label{lemma:StructuralVanishing}
    Let $(\cJ,\cL)$ be a tight joints configuration, and consider $\varepsilon, t_{\textup{f}}, f$ chosen above.
    Then for every $p'\in\cJ$ and $\ell\in\cL$, we have
    \[\mult(f,p')\geq (1-\varepsilon)\frac{dn}{M}\]
    and
    \[\mult(f,\ell)\geq \left(\frac{d-1}{M}-O_{d,M}(\varepsilon)\right)n.\]
    Moreover, there exists $\vec{\alpha}\in S_p(n,T_{\vec{z}})\backslash S_p(n,T_{\vec{z}},t_{\textup{f}})$ such that $D^{\vec{\alpha}}_pf(p)\neq 0$.
\end{lemma}
\begin{proof}
    Since $f$ is $(S_{p'}(n,T_{\vec{z}},t_{\textup{f}}))_{p'\in\cJ}$-vanishing, by \cref{lemma:GeneralVanishing}, we have $\mult(f,p')\geq v_p(t_{\textup{f}})$ and $\mult(f,\ell)\geq v_{\ell}(t_{\textup{f}})$ for every $p'\in\cJ$ and $\ell\in\cJ$.
    By the choice of $\varepsilon$, we have $v_p(t_{\textup{f}})\geq (1-\varepsilon)dn/M$.
    Plugging this in to the definition of $v_{\ell}$, we get
    \[v_{\ell}(t_{\textup{f}}) \geq \frac{(M-d+1)(1-\varepsilon)\frac{dn}{M}-n}{M-d} = \left(\frac{d-1}{M}-\frac{d(M-d+1)}{M(M-d)}\varepsilon\right)n.\]
    Lastly, if for every $\vec{\alpha}\in S_p(n,T_{\vec{z}})\backslash S_p(n,T_{\vec{z}},t_{\textup{f}})$ we have $D^{\vec{\alpha}}_pf(p)=0$, then $f$ is $(S_{p'}(n,T_{\vec{z}}))_{p'\in\cJ}$-vanishing.
    Therefore by \cref{lemma:GeneralVanishing}, we must have $f\equiv 0$, which is a contradiction.
    Thus there must exist $\vec{\alpha}\in S_p(n,T_{\vec{z}})\backslash S_p(n,T_{\vec{z}},t_{\textup{f}})$ such that $D^{\vec{\alpha}}_pf(p)\neq 0$.
\end{proof}

\cref{lemma:StructuralVanishing} shows that there is some $\vec{\alpha}$ in $S_p(n,T_{\vec{z}})$ that is ``close to the apex'' of $S_p(n,T_{\vec{z}})$ so that $D^{\vec{\alpha}}_pf(p)\neq 0$.
To make the discussion more concise, we call a vector $\vec{\alpha}$ \emph{$\varepsilon$-close to the apex} if $\abs{\alpha_i-n/M}\leq \varepsilon n/M$ for all $i\in [d]$.
The statement below shows a way to deduce $O(\varepsilon)$-closeness from a weaker condition that $\abs{\vec{\alpha}}$ is close to $dn/M$ assuming $D^{\vec{\alpha}}_pf(p)\neq 0$.

\begin{lemma}\label{lemma:ClostToApex}
    If $\vec{\alpha}\in\ZZ_{\geq 0}^d$ satisfies that $\abs{\abs{\vec{\alpha}}-dn/M}= O_{d,M}(\varepsilon n)$ and $D^{\vec{\alpha}}_pf(p)\neq 0$ for some $p\in\cJ$, then $\vec{\alpha}$ is $O_{d,M}(\varepsilon)$-close to the apex.
\end{lemma}
\begin{proof}
    Since $D^{\vec{\alpha}}_pf(p)\neq 0$, this shows that for any $i\in[d]$ we have
    \[\abs{\vec{\alpha}}-\alpha_i \geq \mult(f,\ell_{p,i}) \geq \left(\frac{d-1}{M}-O_{d,M}(\varepsilon)\right)n.\]
    Rearranging, we get
    \[\alpha_i\leq \left(\frac{1}{M}+O_{d,M}(\varepsilon)\right)n\]
    for all $i\in[d]$.
    Therefore
    \begin{align*}
        \alpha_i =& \abs{\vec{\alpha}}-(\abs{\vec{\alpha}}-\alpha_i)\\
        \geq& \left(\frac{d}{M}-O_{d,M}(\varepsilon)\right)n-(d-1)\left(\frac{1}{M}+O_{d,M}(\varepsilon)\right)n\\
        \geq& \left(\frac{1}{M}-O_{d,M}(\varepsilon)\right)n
    \end{align*}
    for every $i\in[d]$, showing that $\vec{\alpha}$ is $O_{d,M}(\varepsilon)$-close to the apex.
\end{proof}

\begin{proof}[Proof of \cref{lemma:Structural}]
    After relabeling, we may assume that $\ell = \ell_{p,1}=\ell_{p',1}$.
    We first take $\vec{\alpha}$ as in \cref{lemma:StructuralVanishing}.
    Take $\vec{\alpha}^*$ such that $\alpha_1^*$ is maximized given $\abs{\vec{\alpha}^*}=\abs{\vec{\alpha}}$ and $D^{\vec{\alpha}^*}_pf(p)\neq 0$.
    Note that $\abs{\abs{\vec{\alpha}^*}-dn/M}=\abs{\abs{\vec{\alpha}}-dn/M}\leq \varepsilon dn/M$.
    From now on, we will fix this choice of $\vec{\alpha}^*$ and recycle the notation $\vec{\alpha}$.
    By \cref{lemma:ClostToApex}, if $\vec{\alpha}\in\ZZ^d_{\geq 0}$ satisfies $\abs{\vec{\alpha}}=\abs{\vec{\alpha}^*}$ and $D^{\vec{\alpha}}_pf(p)\neq 0$, then $\vec{\alpha}$ is $O_{d,M}(\varepsilon)$-close to the apex. 

    Let $\vec{e}_{\ell_{p,i}}$ be the unit vector pointing in the direction of $\ell_{p,i}$.
    Then we may expand
    \[f(p+x_1\vec{e}_{\ell_{p,1}}+\cdots +x_d\vec{e}_{\ell_{p,d}})=\sum_{\vec{\alpha}\in \ZZ^d_{\geq 0}}\left(D^{\vec{\alpha}}_pf(p)\right)x_1^{\alpha_1}\cdots x_d^{\alpha_d}.\]
    Let $g$ be a polynomial of degree $\abs{\vec{\alpha}^*}$ such that
    \[g(p+x_1\vec{e}_{\ell_{p,1}}+\cdots +x_d\vec{e}_{\ell_{p,d}})=\sum_{\substack{\vec{\alpha}\in\ZZ^d_{\geq 0}\\ \abs{\vec{\alpha}} = \abs{\vec{\alpha}^*}}}\left(D^{\vec{\alpha}}_pf(p)\right)x_1^{\alpha_1}\cdots x_d^{\alpha_d}.\]
    In other words, $g$ is the degree-$\abs{\vec{\alpha}^*}$ part of $f$ when expanded at $p$.
    By the choice of $\vec{\alpha}^*$, the summand on the right hand side is $0$ if $\min_{i\in[d]}\alpha_i < (1-O_{d,M}(\varepsilon))n/M$.
    By the maximality of $\alpha_1^*$, we know that $D^{(\alpha_1^*,0,\ldots,0)}g(p+x_1\vec{e}_{\ell_{p,1}}+\cdots +x_d\vec{e}_{\ell_{p,d}})$ is not identically zero, but it does not depend on $x_1$.
    Therefore there exists a nonzero polynomial $h\in \FF[x_2,\ldots, x_d]$ such that
    \[D^{(\alpha_1^*,0,\ldots,0)}_pg(p+x_1\vec{e}_{\ell_{p,1}}+\cdots +x_d\vec{e}_{\ell_{p,d}})=\left(x_2\cdots x_d\right)^{\left\lceil(1-O_{d,M}(\varepsilon))n/M\right\rceil}h(x_2,\ldots, x_d).\]

    We will now rewrite everything in the coordinate system given by $p'$ at $p$ by setting $y_1,\ldots, y_d$ be some linear combinations of $x_1,\ldots, x_d$ so that $p+x_1\vec{e}_{\ell_{p,1}}+\cdots +x_d\vec{e}_{\ell_{p,d}} = p+y_1\vec{e}_{\ell_{p',1}}+\cdots +y_d\vec{e}_{\ell_{p',d}}$.
    The key claim to show is that $y_2,\ldots, y_d$ is a permutation of $x_2,\ldots,x_d$.
    To this end, we will expand $D^{(\alpha_1^*,0,\ldots,0)}_pg$ in this new coordinate system.
    Since it does not depend on $x_1$ (and thus $y_1$), we get
    \[D^{(\alpha_1^*,0,\ldots,0)}_pg(p+y_1\vec{e}_{\ell_{p',1}}+\cdots +y_d\vec{e}_{\ell_{p',d}})=\sum_{\substack{\vec{\alpha}\in \ZZ_{\geq 0}^d\\ \abs{\vec{\alpha}}=\abs{\vec{\alpha}^*}-\alpha_1^*\\\alpha_1=0}}\left(D^{\vec{\alpha}}_{p'}D^{(\alpha_1^*,0,\ldots,0)}_pg(p)\right)y_2^{\alpha_2}\cdots y_d^{\alpha_d}.\]
    Note that $D^{(\alpha_1^*,0,\ldots,0)}_p=D^{(\alpha_1^*,0,\ldots,0)}_{p'}$ as $\ell_{p,1} = \ell_{p',1}$ by \cref{corollary:ChangeCoordinate}, and $D^{\vec{\alpha}}_{p'}D^{(\alpha_1^*,0,\ldots,0)}_{p'}=D^{(\alpha_1^*,\alpha_2,\ldots, \alpha_d)}_{p'}$ if $\alpha_1=0$ by \cref{proposition:HasseComposition}.
    Therefore, we may rewrite the identity as 
    \[D^{(\alpha_1^*,0,\ldots,0)}_pg(p+y_1\vec{e}_{\ell_{p',1}}+\cdots +y_d\vec{e}_{\ell_{p',d}})=\sum_{\substack{\vec{\alpha}\in \ZZ_{\geq 0}^d\\ \abs{\vec{\alpha}}=\abs{\vec{\alpha}^*}\\\alpha_1=\alpha_1^*}}\left(D^{\vec{\alpha}}_{p'}g(p)\right)y_2^{\alpha_2}\cdots y_d^{\alpha_d}=\sum_{\substack{\vec{\alpha}\in \ZZ_{\geq 0}^d\\ \abs{\vec{\alpha}}=\abs{\vec{\alpha}^*}\\\alpha_1=\alpha_1^*}}\left(D^{\vec{\alpha}}_{p'}f(p)\right)y_2^{\alpha_2}\cdots y_d^{\alpha_d}.\]
    Now if $\abs{\vec{\alpha}} = \abs{\vec{\alpha}^*}$ and $\alpha_1= \alpha_1^*$ but $D^{\vec{\alpha}}_{p'}f(p)\neq 0$, let $\vec{\alpha}' = (0,\alpha_2,\ldots, \alpha_d)$.
    Then $D^{\vec{\alpha}'}f$ does not vanish identically on $\ell$.
    Moreover, for any $p''\in\ell$, we know that
    \[\mult\left(\left.D^{\vec{\alpha}'}_{p'}f\right|_{\ell},p''\right)\geq \mult(f,p'')-\abs{\vec{\alpha}'}\geq \left(\frac{1}{M}-O_{d,M}(\varepsilon)\right)n\]
    when $n$ is sufficiently large.
    Here, we use the fact that $\mult(f,p'')\geq (1-\varepsilon)dn/M$ by \cref{lemma:StructuralVanishing} and $\abs{\vec{\alpha}'}=\alpha_2^*+\cdots+\alpha_d^*\leq (1+O_{d,M}(\varepsilon))(d-1)n/M$ by the $O_{d,M}(\varepsilon)$-closeness of $\vec{\alpha}^*$.
    When $p''=p'$, we have another lower bound
    \[\mult\left(\left.D^{\vec{\alpha}'}_{p'}f\right|_{\ell},p'\right)\geq \mult(f,\ell_{p',i})-\left(\abs{\vec{\alpha}'}-\alpha_i'\right)\geq \alpha_i'-O_{d,M}(\varepsilon) n\]
    for every $i=2,\ldots, d$, where the second inequality comes from \cref{lemma:StructuralVanishing} and the upper bound on $\abs{\alpha}'$ we just used.
    The inequality
    \[\mult\left(\left.D^{\vec{\alpha}'}_{p'}f\right|_{\ell},p'\right)+\sum_{p''\in\ell\backslash \{p'\}}\mult\left(\left.D^{\vec{\alpha}'}_{p'}f\right|_{\ell},p''\right)\leq \deg\left.D^{\vec{\alpha}'}_{p'}f\right|_{\ell} = n - \abs{\vec{\alpha}'}\]
    yields
    \[\alpha_i'+(M-d)\left(\frac{1}{M}-O_{d,M}(\varepsilon)\right)n\leq n-\left(\alpha_2^*+\cdots+\alpha_d^*\right)\leq n-(d-1)\left(\frac{1}{M}-O_{d,M}(\varepsilon)\right)n\]
    for every $i=2,\ldots, d$.
    Rearranging, we get $\alpha_i'\leq (1/M +O_{d,M}(\varepsilon))n$ for $i=2,\ldots, d$.
    Therefore for any $i=2,\ldots, d$, we have
    \[\alpha_i' = \abs{\vec{\alpha}'}-\left(\abs{\vec{\alpha}'}-\alpha_i'\right)\geq (d-1)\left(1-O_{d,M}(\varepsilon)\right)\frac{n}{M}-(d-2) \left(\frac{1}{M}+O_{d,M}(\varepsilon)\right)n\geq \left(\frac{1}{M}-O_{d,M}(\varepsilon)\right)n.\]
    As a consequence, we have that
    \[D^{(\alpha_1^*,0,\ldots,0)}_pg(p+x_1\vec{e}_{\ell_{p,1}}+\cdots +x_d\vec{e}_{\ell_{p,d}})=\left(x_2\cdots x_d\right)^{\left\lceil(1-O_{d,M}(\varepsilon))n/M\right\rceil}h(x_2,\ldots, x_d)\]
    is divisible by $(y_2\cdots y_d)^{\left\lceil(1-O_{d,M}(\varepsilon))n/M\right\rceil}.$
    Note that $\deg h = \abs{\vec{\alpha}^*}-\alpha_1^*-(d-1)\left\lceil(1-O_{d,M}(\varepsilon))n/M\right\rceil$, which can bounded by $O_{d,M}(\varepsilon n)$.
    Therefore, as long as $\varepsilon$ is sufficiently small with respect to $M$ and $d$, we have that $y_i$ divides $(x_2\cdots x_d)$ for every $i=2,\ldots, d$.
    Therefore $y_i = x_{\sigma(i)}$ for some $\sigma:\{2,\ldots, d\}\to\{2,\ldots, d\}$.
    Since $y_2,\ldots, y_d$ must be linearly independent, $\sigma$ must be a bijection.
    Note that $x_i=0$ defines the hyperplane defined by the $(d-1)$ lines through $p$ that are not $\ell_{p,i}$, and $y_i=0$ defines the hyperplane defined by the $(d-1)$ lines through $p'$ that are not $\ell_{p,i}$.
    The desired statement then follows.
\end{proof}

\subsection{Further discussion}\label{subsection:StructuralDiscussion}
From our method, we can only get the structural theorem when the configuration is tight. There are some obstacles when the configuration is not tight but also not far away from tight.

First, if the number of joints is $o(1)$ away from tight, then we can not expect that \cref{lemma:Structural} holds on the whole configuration. Here is an example. 
Consider a tight joints configuration and add a point $p$ on one of the lines $\ell$, then pick $d-1$ lines through $p$ with arbitrary directions that make $p$ a joint. 
In this case, we can pick the directions so that $p$ does not share any hyperplane with any other $p'\in\ell$.  
This configuration indicates that we should consider a sub-configuration instead.

There is another configuration that is nearly tight without such a strong structure.
In fact, we will consider the following projected generically induced configuration.
The example comes from an construction by Bollob\'as and Eccles \cite{BE15}.

\begin{example}\label{example:BE}
    Take $6$ hyperplanes $H_1,\ldots, H_6$ in $\FF^5$ in general position.
    Consider the $5$-wise intersections $P_1,\ldots, P_6$ of the six hyperplanes, where $P_i$ is the intersection of $(H_j)_{j\neq i}$.
    Collect all $4$-wise intersections of $H_1,\ldots, H_6$ except $H_1\cap H_2\cap H_3\cap H_4$, $H_1\cap H_2\cap H_5\cap H_6$ and $H_3\cap H_4\cap H_5\cap H_6$.
    Those are $\binom{6}{4}-3 = 12$ lines, and we call the set of those lines $\cL_0$.
    Note that it is easy to verify that each $P_i$ lies on precisely four lines in $\cL_0$, whose directions are linearly independent.

    Take a general projection $\pi:\FF^5\to \FF^4$ so that for it maps the $P_i$'s to different points, the lines in $\cL_0$ to different lines, and for each $i=1,\ldots, 6$, the images of the lines through $P_i$ are still linearly independent.
    Let $\cJ$ be the set of images of $P_1,\ldots, P_6$, and $\cL$ be the collection of images of lines in $\cL_0$.
    Then $(\cJ,\cL)$ is a joints configuration with six joints and twelve lines in $\FF^4$.
    Note that \cref{theorem:SharpJoints} shows that if $\abs{\cJ}=6$, then $x\approx 5.1458$ and hence $\abs{\cL}>11.1848$.
    Therefore this configuration is optimal.

    To see if the joints share all the hyperplanes, we first determine the hyperplnaes containing three out of four lines passing through a given joint.
    We take $\pi(P_1)$ for example.
    The four lines through $\pi(P_1)$ are $\pi\left(H_2\cap H_4\cap H_5\cap H_6\right)$, $\pi\left(H_2\cap H_3\cap H_5\cap H_6\right)$, $\pi\left(H_2\cap H_3\cap H_4\cap H_6\right)$ and $\pi\left(H_2\cap H_3\cap H_4\cap H_5\right)$.
    Therefore, the four hyperplanes determined by any three of them are $\pi\left(H_2\cap H_3\right)$, $\pi\left(H_2\cap H_4\right)$, $\pi\left(H_2\cap H_5\right)$ and $\pi\left(H_2\cap H_6\right)$.
    A similar argument works for all other joints.

    Note that $\pi(P_1)$ and $\pi(P_2)$ do not lie on the same line, and all other joints are connected to $\pi(P_1)$ via a line.
    We may take $\pi(P_3)$ for example, as the other cases are the same up to symmetry.
    Then $\pi(P_1)$ and $\pi(P_3)$ only share the hyperplane $\pi\left(H_2\cap H_4\right)$ as long as the projection $\pi$ is general enough.
    We may prove similarly that any two distinct joints on the same line share exactly one hyperplane.
\end{example}

The example above shows that if the number of joints is too small, we may lose the structure we are looking for even on optimal configurations.
We suspect that similar phenomenon does not occur when the number of joints is large enough, and we will elaborate on this in \cref{section:Graph}.
Motivated by the previous two examples, we make the following conjecture.
\begin{conjecture}
    Assume $(\cJ,\cL)$ is a joints configuration in $\FF^d$ with $J=(1-o(1))\binom{x}{d}$ and $L=\binom{x}{d-1}$ for some $x\geq d$. Then there exists a subset of joints $\cJ'\subseteq \cJ$ with $|\cJ'|=(1-o(1))J$ satisfying \cref{lemma:Structural} with $\cL$.
\end{conjecture}
Note that the $o(1)$-factor prevents \cref{example:BE} from being a counterexample.

The next obstacle we have in mind is that the structural theorem definitely fails when the configuration is a multiplicative constant away from being tight.
\begin{example}
    We first pick some distinct real numbers $c_1,\dots,c_n,d_1,\dots,d_n>0$ such that the unique solution $y=y_{ij}$ to the equation $c_iy+d_i=c_jy+d_j$ are all distinct for $1\leq i<j\leq n$, and $kc_i$ are all distinct for all $i,k\in [n]$.
    In $\RR^3$, define reguli $R_i\eqdef\{(x,y,z)\in\RR^3:z=x(c_iy+d_i)\}$ for $i=1,\dots,n$. Also, let $H_k$ denote the hyperplane $x=k$ for $k=1,\dots,n$. It follows that, for any $1\leq i<j\leq n$, the intersection $R_i\cap R_j$ is the union of two lines
    \[\{(t,y_{ij},y_{ij}t):t\in\RR\}\cup\{(0,t,0):t\in\RR\}.\]
    Set $\ell_{ij}$ to be the first line above.
    Also, note that, for any $i,k\in [n]$, the intersection $R_i\cap H_k$ is a line
    \[\ell'_{ik}\eqdef\{(k,t,k(c_it+d_i)):t\in\RR\}.\]
    It is easy to check that, for any $1\leq i<j\leq n,k\in [n]$, the intersection $R_i\cap R_j\cap H_k$ is the point $p_{ijk}\eqdef(k,y_{ij},k(c_iy_{ij}+d_i))$, which is a joint formed by $\ell_{ij},\ell'_{ik},\ell'_{jk}$.
    This is because the directions of $\ell_{ij},\ell'_{ik},\ell'_{jk}$ are 
    \[v_{ij}\eqdef (1,0,y_{ij}),v'_{ik}\eqdef (0,1,kc_i),v'_{jk}\eqdef (0,1,kc_j)\]
    respectively, which are linearly independent. 
    
    Let $\cL\eqdef\{\ell_{ij}:1\leq i<j\leq n\}\cup\{\ell'_{ik}:i,k\in [n]\}$ and $\cJ\eqdef\{p_{ijk}:1\leq i<j\leq n,k\in [n]\}$. We claim that $(\cJ,\cL)$ is a joints configuration that any two joints on the same line do not share all the hyperplanes. For any two joints $p_{ijk},p_{ijk'}$ on the line $\ell_{ij}$, we know that $v_{ij}$ is not linearly dependent with any two of the vectors $v'_{ik},v'_{jk},v'_{ik'},v'_{jk'}$. Thus, $p_{ijk}$ and $p_{ijk'}$ do not share any hyperplane. For any two joints on $\ell'_{ik}$, we may without loss of generality assume they are $p_{ijk},p_{ij'k}$. We know that they share the hyperplane $H_k$ since $v'_{ik},v'_{jk},v'_{j'k}$ are linearly dependent. However, since $v'_{ik},v_{ij},v_{ij'}$ are linearly independent, they do not share the other hyperplane.
    
    In this example, we have $L=\binom{n}{2}+n^2$ and $J=\binom{n}{2}n=(1-o(1))\sqrt{\frac{2}{27}}L^{3/2}$, which is only a multiplicative constant away from being tight. 
\end{example}

It is interesting to ask how close to the optimal $J$ should be in order to guarantee the existence of two points on the same line that share all the hyperplanes.
The contrapositive of the question is the following.
\begin{question}
    Let $(\cJ,\cL)$ be a joint configuration so that for any $p\neq p'$ that lie on the same line, they do not share all the hyperplanes.
    What is the maximum of $J$ given $L$?
\end{question}

Another main issue in the structural theorem is that we really rely on the polynomial method being absolutely tight. 
This allows us to control the number of joints on each line and to solve for the equations in \cref{lemma:SolveEq}.
However, we do not even have a tight bound when $J$ is not of the form $\binom{M}{d}$ for some integer $M\geq d$. 
We believe that there is a chance to get the structural results for all $J$'s if we could prove the optimal bound using polynomial method for each $J$.
We will discuss what we might guess about the optimal bound in \cref{section:Graph}.

\section{Joints of Curves}\label{section:Curve}
In this section, we briefly discuss how to generalize the argument in \cref{section:Sharp} to curves.
We will essentially replace each part in \cref{section:Sharp} with an appropriate analogue for curves.
We start with \cref{lemma:BasicVanishing}.

\begin{lemma}\label{lemma:CurveBasicVanishing}
    Let $n$ be a nonnegative integer.
    Let $(\cJ,\cC)$ be a joint-of-curve configuration in $\FF^d$, and let $\vec{r}\in \ZZ_{\geq 0}^{\cJ}$.
    If $f\in \FF[x_1,\ldots, x_d]$ is a polynomial of degree at most $n$ such that
    \[\mult(f,p)\geq r_p\quad\forall p\in \cJ,\]
    then
    \[\mult(f,\ell)\geq \frac{\sum_{p\in \ell}r_p-n\deg \ell}{\#\{p\in\ell\}-\deg \ell}\quad\forall \ell\in\cL.\]
\end{lemma}

\begin{proof}
    For any $\vec{\alpha}\in \ZZ_{\geq 0}^d$, we have $\deg \Hasse^{\vec{\alpha}}f\leq n-\abs{\vec{\alpha}}$ and $\mult(\Hasse^{\vec{\alpha}}f,p)\geq r_p-\abs{\vec{\alpha}}$.
    Note that if $\sum_{p\in \ell} \mult(\Hasse^{\vec{\alpha}}f,p)>\deg\Hasse^{\vec{\alpha}}f\cdot\deg \ell$, then by \cref{proposition:BasicVanishingCurves} we must have $\left.\Hasse^{\vec{\alpha}}f\right|_{\ell}\equiv 0$.
    Rearranging, we have that $\left.\Hasse^{\vec{\alpha}}f\right|_{\ell}\equiv 0$ whenever
    \[\abs{\vec{\alpha}} < \frac{\sum_{p\in \ell}r_p-n\deg \ell}{\#\{p\in\ell\}-\deg \ell}.\]
    The desired statement then follows.
\end{proof}

Again, we may define the \textit{priority timestamp} and the functions $v_p^{\tT}(t)$ in the same way, and define
\[v^{\tT}_{\ell}(t)\eqdef \frac{\sum_{p\in\ell}v^{\tT}_p(t)-n\deg\ell}{\#\{p\in\ell\}-\deg\ell}\quad\forall \ell\in \cL .\]
The definition of $S_p(n,T,t_{\textup{f}})$ remains the same. 
The definition of being $(S_p(n,T,t_{\textup{f}}))_{p\in\cJ}$-vanishing is the same as \cref{definition:SpVanishing}.
We can also prove a curve version of \cref{lemma:GeneralVanishing}.

\begin{lemma}[General vanishing lemma]\label{lemma:CurveGeneralVanishing}
    Let $n$ be a nonnegative integer.
    Let $(\cJ,\cC)$ be a joint-of-curve configuration in $\FF^d$ with a priority timestamp $T$.
    Let $t_{\textup{f}}\in \{0,\ldots, (n+1)\abs{\cJ}\}$ be the terminating time.
    Then for any polynomial $f\in \FF[x_1,\ldots, x_d]$ with degree at most $n$, if it is $(S_p(n,T,t_{\textup{f}}))_{p\in\cJ}$-vanishing, then $\mult(f,p)\geq v_p(t_{\textup{f}})$ for any $p\in\cJ$ and $\mult(f,\ell)\geq v_{\ell}(t_{\textup{f}})$ for any $\ell\in\cC$.
\end{lemma}
\begin{proof}
    The proof is almost identical to the proof of \cref{lemma:GeneralVanishing}.
    The two modifications we need make are the following.
    First, we now apply \cref{lemma:CurveBasicVanishing} instead of \cref{lemma:BasicVanishing}.
    In addition, previously $D^{\vec{\alpha}}_p$ was a linear combination of Hasse derivatives of degree $\abs{\vec{\alpha}}$, but now $D^{\vec{\alpha}}_p$ is a linear combination of Hasse derivatives of degrees at most $\abs{\vec{\alpha}}$.
    This does not hurt us, as in the inductive step, we still have that $D^{\vec{\alpha}-\alpha_i\vec{e}_i}_pf$ vanishes identically on $\ell_{p,i}$ if $\abs{\vec{\alpha}}-\alpha_i<v_{\ell_{p,i}}(t_{\textup{f}}-1)$.
    The rest of the argument works verbatim.
\end{proof}

As in \cref{section:Sharp}, this immediately gives the following inequality.

\begin{corollary}\label{corollary:CurveParamCounting}
    Let $n$ be a nonnegative integer.
    Let $(\cJ,\cC)$ be a joint-of-curve configuration in $\FF^d$ with a priority timestamp $T$.
    Then
    \[\sum_{p\in\cJ}\abs{S_p(n,T)}\geq \binom{n+d}{d}\]
\end{corollary}

As in \cref{section:Sharp}, let $\vec{z}\in \RR_{>0}^{\cJ}$, and let $T_{\vec{z}}$ be the associated priority timestamp of $\vec{z}$.
We still have
\[v_{p'}(t) = (\abs{\vec\alpha}-n)\frac{z_{p'}}{z_p}+n+O(1)=\left(\left(\abs{\vec{a}}-1\right)\frac{z_{p'}}{z_p}+1+o(1)\right)n.\]
Plugging this in, the condition that $\vec{\alpha}\in S_p(n,T_{\vec{z}})$ can be reformulated as $\abs{\vec{a}}\leq 1$ and 
\[\abs{\vec{a}}-a_i\geq \frac{\sum_{p'\in \ell_{p,i}}z_{p'}}{\left(\#\{p'\in\ell_{p,i}\}-\deg\ell_{p,i}\right)z_p}\cdot\abs{\vec{a}}+\left(1-\frac{\sum_{p'\in \ell_{p,i}}z_{p'}}{\left(\#\{p'\in\ell_{p,i}\}-\deg\ell_{p,i}\right)z_p}\right)+o(1)\]
by plugging in the definition of $v_{\ell_{p,i}}(t)$.
Once again, to simplify, given $\vec{z}\in \ZZ_{\geq 0}^{\cJ}$, for any $p\in\cJ$ and $i\in [d]$, let
\[\beta_{p,i}\eqdef 1-\left(\#\{p'\in \ell_{p,i}\}-\deg\ell_{p,i}\right)\frac{z_p}{\sum_{p'\in \ell_{p,i}}z_{p'}}.\]
The second condition can now be written as
\[\abs{\vec{a}}-a_i\geq \frac{1}{1-\beta_{p,i}}\abs{\vec{a}}-\frac{\beta_{p,i}}{1-\beta_{p,i}}+o(1),\]
or equivalently,
\[a_i+\beta_{p,i}(\abs{\vec{a}}-a_i)\leq \beta_{p,i}+o(1).\]
With all the notations set up, the following two statements follow with the same proof as in \cref{section:Sharp}.

\begin{lemma}[Continuous approximation of $S_p$ with $\cS$]\label{lemma:CurveContiApprox}
    Given a joint-of-curve configuration $(\cJ,\cC)$ in $\FF^d$, let $\vec{z}\in \RR^{\cJ}_{>0}$ be a vector and $T_{\vec{z}}$ be its associated priority timestamp.
    Define $\beta_{p,i}$ as above for any $p\in\cJ$ and $i\in [d]$.
    Then
    \[\abs{S_p(n,T_{\vec{z}})} = \left(\vol\cS(\beta_{p,1},\ldots, \beta_{p,d})+o(1)\right)\cdot n^d\]
    for any $p\in \cJ$.
\end{lemma}

\begin{corollary}\label{corollary:CurveVanishingIneq}
    Let $(\cJ,\cC)$ be a joint-of-curve configuration in $\FF^d$ with a vector $\vec{z}\in \RR_{>0}^{\cJ}$.
    Compute $\beta_{p,i}$ accordingly for any $p\in \cJ$ and $i\in [d]$ as above.
    Then
    \[\sum_{p\in\cJ}\vol\cS(\beta_{p,1},\ldots, \beta_{p,d})\geq \frac{1}{d!}.\]
\end{corollary}

Similar as before, we know that $\sum_{p\in\cJ,i\in [d]}\beta_{p,i}=\deg\cC$. Again, if we set $b_{p,i}=\frac{z_p}{\sum_{p'\in \ell_{p,i}}z_{p'}}$, we want to solve the system 
\[\sigma_p\eqdef\sum_{i=1}^d \left(\#\{p'\in \ell_{p,i}\}-\deg\ell_{p,i}\right)b_{p,i}=\frac{dJ-\deg\cC}{J}~\forall p\in\cJ.\]

As in \cref{section:Sharp}, we will do so on critical configurations, but this time with an additional technical assumption that each curve $\ell$ contains at least $\deg \ell$ points.

\begin{definition}
A joint-of-curve configuration $(\cJ,\cC)$ is \emph{critical} if for any joint-of-curve configuration $(\cJ',\cL')$ with $J'<J$, we have $J'/\deg\cL'<J/\deg\cL$.
\end{definition}
\begin{lemma}\label{lemma:CurveSolveEq}
    If $(\cJ,\cC)$ is a critical joint-of-curve configuration where each curve $\ell\in\cC$ contains at least $\deg\ell$ joints, then the system  
    \begin{align*}
    \sigma_p=\frac{dJ-\deg\cL}{J}~\forall p\in\cJ
    \end{align*}
    has a solution $\vec{z}$ with $z_p>0$ for all $p\in\cJ$.
\end{lemma}
\begin{proof}
    We first show that the function $g$ taking $\vec{z}$ as the input and $\vec{\sigma}$ as its output satisfies all conditions in \cref{lemma:EquationTechnical}.
    The proof is identical to the proof of \cref{lemma:EquationTechnical2}, except that $\{p'\in J': p'\in \ell\}-1$ should be replaced by $\{p'\in J': p'\in \ell\}-\deg \ell$ and $L,L',\Tilde{L}$ should be replaced by $\deg\cC, \deg \cC'$ and $\deg\Tilde{\cC}$.
    Note that $\{p'\in J': p'\in \ell\}-\deg \ell$ is nonnegative by the assumption, so the proof really works verbatim after the modification.

    After showing this, the rest of the proof is identical to the proof of \cref{lemma:SolveEq} with the same modification.
\end{proof}
Now, we are ready to prove \cref{theorem:JointsOfCurves}.
\begin{proof}[Proof of \cref{theorem:JointsOfCurves}]
First of all, we may assume that $(\cJ,\cC)$ is critical by the argument we used for \cref{theorem:SharpJoints} in \cref{section:Sharp}.
We may also assume that each curve $\ell\in\cC$ contains at least $\deg\ell$ joints, as otherwise, we may remove $\ell$ and simply add a line to each $\cJ\cap \ell$ so that each of them is still a joint, while the total degree of the curves decreases.
From \cref{lemma:CurveSolveEq}, we know that there is a solution $\vec{z}$ with $z_p>0$ for all $p\in\cJ$ such that $\avgbeta_p\eqdef(\beta_{p,1}+\dots+\beta_{p,d})/d=\deg\cC/dJ$. Combining \cref{corollary:CurveVanishingIneq} and \cref{lemma:Vol}, we have 
\[\frac{J}{d!\binom{dJ/\deg\cL+d-1}{d}}\geq\frac{1}{d!}.\]
That is, $J\geq \binom{dJ/\deg\cC+d-1}{d}$. If $J=\binom{x}{d}$, then $x\geq dJ/\deg\cC+d-1$, and hence we get the desired bound $\deg\cC\geq dJ/(x-d+1)=\binom{x}{d-1}$.
\end{proof}
\section{Optimal Constant with New Multiplicities}\label{section:multiplicity}
In \cite{YZ23}, Zhao and the second author showed that $N$ lines in $\FF^d$ form at most $(d-1)!^{1/(d-1)}/d \cdot N^{d/(d-1)}$ joints, and this leading constant is optimal as shown by \cref{example:tight} where we take $M$ hyperplanes in general position and their $\binom{M}{d-1}$ $(d-1)$-wise intersections to be the lines.
The situation for joints of varieties is slightly different.
Take joints of $2$-flats in $\FF^6$ for example.
The bound that was proved in \cite{TYZ22} states that $N$ $2$-flats in $\FF^6$ form at most $\sqrt{10/3}N^{3/2}$ joints, though there is still no known matching construction.
One na\"ive attempt is to take $M$ hyperplanes in general positions and their $\binom{M}{4}$ $4$-wise intersections, giving $\binom{M}{4}$ $2$-flats and $\binom{M}{6}$ joints.
Note, however, that if $N=\binom{M}{4}$, then $\binom{M}{6} = (4!N)^{3/2}/6! = \sqrt{2/75}N^{3/2}$, which is $5^{3/2}$ away from the upper bound $\sqrt{10/3}N^{3/2}$ we are aiming for.

Though the construction does not quite match the upper bound, it is clear that each joint is formed in many ways.
Indeed, if a joint is the intersection of six hyperplanes $H_1,\ldots, H_6$, then it is a joint formed by $H_1\cap H_2\cap H_3\cap H_4, H_1\cap H_2\cap H_5\cap H_6$ and $H_3\cap H_4\cap H_5\cap H_6$, and this statement is also true with $H_1,\ldots, H_6$ permuted arbitrarily.
If there is some way of defining multiplicities so that each joint in this construction contributes $5^{3/2}$ to the count, and we can still prove the same upper bound, then the upper bound would be tight.
Note that in this case, for each joint, there are exactly $15$ $2$-flats through the joint, and they can be split into $5$ triples that form a joint.
This inspires us to define a new notion of multiplicities based on the number of disjoint triplets that form a joint.

In this section, we will expand the idea above.
In the first subsection, we define a new multiplicity and reprove the same upper bound as in \cite{TYZ22}.
We will also show that this time, the constant in the upper bound is optimal.
In the second subsection, we examine the proof in the first subsection and extract yet another notion of multiplicity directly from the proof.
We will then show that this second notion of multiplicity is in some sense the strongest---the new upper bound implies directly the upper bound in the first subsection and also the general theorem on joints with multiplicities \cite[Theorem 1.10 (b)]{TYZ22}.

\subsection{First definition of multiplicity}
In this subsection, we will keep the following notation.
Let $k_1,\ldots,k_r,m_1,\ldots, m_r$ be some positive integers and let $d=m_1k_1+\cdots+m_rk_r$.
For each $i=1,\ldots, r$, let $\cV_i$ be a finite multiset of $k_i$-dimensional varieties in $\FF^{d}$.
Let $\cJ$ be the set of joints formed by $m_i$ elements from $\cV_i$ for each $i\in[r]$.
We will now define the multiplicity.

\begin{definition}
    For each $p\in\cJ$, let $G_p$ be an $(m_1+\cdots+m_r)$-uniform hypergraph on the vertex (multi)set $\cV_1\cup\cdots\cup\cV_r$, where we add an edge $S = S_1\cup\cdots\cup S_r$ for any $(S_1,\ldots, S_r)\in \binom{\cV_1}{m_1}\times\cdots\times\binom{\cV_r}{m_r}$ such that $S_1\cup\cdots\cup S_r$ form a joint at $p$.
    Let $\nu(p)$ be the packing number of $G_p$, i.e. the maximum number of disjoint edges in $G_p$.
\end{definition}

\begin{example}
    Going back to the joints-of-$2$-flats example, if $r=1, k_1=2, m_1=3$ and $\cV_1$ is the set of $4$-wise intersections of $k$ hyperplanes in general positions, then it is easy to verify that $\nu(p)=5$ for any $p\in\cJ$.
    It is thus natural to guess that $\nu(p)^{3/2}$ is the right notion of multiplicities.
\end{example}

We can now state and prove the following.

\begin{theorem}\label{theorem:FirstMult}
    Let $k_1,\ldots,k_r,m_1,\ldots, m_r$ be some positive integers and $d=m_1k_1+\cdots+m_rk_r$.
    For each $i=1,\ldots, r$, let $\cV_i$ be a finite multiset of $k_i$-dimensional varieties in $\FF^{d}$.
    Let $\cJ$ be the set of joints formed by $m_1+\cdots +m_r$ varieties, where there are $m_i$ of them from $\cV_i$ for each $i=1,\ldots, r$.
    Then
    \[\sum_{p\in\cJ}\nu(p)^{(m_1+\cdots+m_r)/(m_1+\cdots+m_r-1)}\leq C^{\nu}_{k_1,\ldots, k_r;m_1,\ldots, m_r}\left(\left(\deg\cV_1\right)^{m_1}\cdots \left(\deg\cV_r\right)^{m_r}\right)^{1/(m_1+\cdots+m_r-1)},\]
    where
    \[C^{\nu}_{k_1,\ldots, k_r;m_1,\ldots, m_r}=\left(\frac{d!}{\prod_{i=1}^{r}k_i!^{m_i}m_i^{m_i}}\right)^{1/(m_1+\cdots+m_r-1)}.\]    
    Moreover, this constant cannot be improved.
\end{theorem}

The proof relies on lemmas and notations from \cite{TYZ22}.
As going through all necessary definitions and motivations behind them would take too long, we make no effort here to be self-contained and simply refer the readers to \cite{TYZ22}.
We will only discuss the necessary modification to prove this stronger result.
\begin{proof}
    We first show that this constant cannot be improved.
    Let $d=m_1k_1+\cdots +m_rk_r$ be the dimension of the underlying space, and $s=m_1+\cdots +m_r$ be the number of varieties needed to form a joint.
    Take $M$ hyperplanes in general position.
    For each $i$, let $\cV_i$ be $m_i\cdot k_i!(d-k_i)!$ copies of $(d-k_i)$-wise intersections of the $M$ hyperplanes.
    Then the joints are the $d$-wise intersections of the $M$ hyperplanes.
    Suppose that $H_1,\ldots, H_d$ are $d$ hyperplanes that intersect at $p$.
    Let $S = \{H_1,\ldots, H_d\}$, and for any $S'\subseteq S$ define $F(S')$ to be the flat defined by the intersection $\bigcap_{H\in S\backslash S'}H$.
    For each permutation $\sigma\in S_d$, we consider the sequence $H_{\sigma(1)},\ldots, H_{\sigma(d)}$ and divide it into $m_1+\cdots +m_r$ chunks, where the first $m_1$ chunks contain $k_1$ terms, the next $m_2$ chunks contain $k_2$ terms, and so on, until the last $m_r$ chunks each containing $k_r$ terms.
    This way, we get subsets $A^{(\sigma)}_{1,1},\ldots, A^{(\sigma)}_{1,m_1},\ldots, A^{(\sigma)}_{r,1},\ldots, A^{(\sigma)}_{r,m_r}$ of $S$.
    Note that for any $i\in[r]$ and $j\in[m_i]$, we have that $F(A^{(\sigma)}_{i,j})\in \cV_i$.
    Moreover, we see that via $F$, those subsets correspond to flats that form a point at $p$, and thus they correspond to a hyperedge in $G_p$.
    For each flat in $\cV_i$ containing $p$ and for every fixed $j\in[m_i]$, the flat is equal to $F(A^{(\sigma)}_{i,j})$ for $k_i!(d-k_i)!$ different permutations $\sigma$.
    Therefore each flat in $\cV_i$ containing $p$ is counted for exactly $m_i\cdot k_i!(d-k_i)!$ times, and so the $d!$ hyperedges we have found in $G_p$ form a packing, showing that $\nu(p)\geq d!$.
    Therefore
    \begin{align*}\sum_{p\in\cJ}\nu(p)^{s/(s-1)}\geq& \binom{M}{d}(d!)^{s/(s-1)}\\
    =& \left(d!^{1/(s-1)}-o(1)\right)M^d.
    \end{align*}
    On the other hand, we have
    \begin{align*}\left(\left(\deg\cV_1\right)^{m_1}\cdots \left(\deg\cV_r\right)^{m_r}\right)^{1/(s-1)}=&\left(\prod_{i=1}^{r}\left(m_ik_i!(d-k_i)!\binom{M}{d-k_i}\right)^{m_i}\right)^{1/(s-1)}\\
    =&\left(\prod_{i=1}^{r}\left(m_i^{m_i}k_i!^{m_i}-o(1)\right)M^{(d-k_i)m_i}\right)^{1/(s-1)}\\
    =&\left(\left(\prod_{i=1}^{r}m_i^{m_i}k_i!^{m_i}\right)^{1/(s-1)}-o(1)\right)M^d,
    \end{align*}
    where the last equality holds as $\sum_{i=1}^{r}(d-k_i)m_i = sd-\sum_{i=1}^{r}m_ik_i = \left(s-1\right)d$.
    Therefore the constant $C^{\nu}_{k_1,\ldots, k_r;m_1,\ldots, m_r}$ cannot be improved.

    To show the theorem, for any $p\in\cJ$, write each edge in $E(G_p)$ as $S_1\cup\cdots \cup S_r$ where $S_i\subseteq \cV_i$ has size exactly $m_i$.
    For each $p\in\cJ$, choose $\cM(p)\subseteq E(G_p)$ to be a packing of size $\nu(p)$, and let $\omega(p)$ be $\nu(p)^{s/(s-1)}$.
    We first suppose that the joints configuration is connected.
    Then by \cite[Lemma 5.10]{TYZ22} with this choice of $\cM(p)$, we have that for any $n\in\NN$, there exists a handicap $\vec{\alpha}$ such that 
    \[\abs{\frac{1}{\omega(p)}\left[\prod_{S_1\cup\cdots\cup S_r\in \cM(p)}\prod_{V\in S_1\cup\cdots\cup S_r}\frac{\abs{\cD_{p,V}(\vec{\alpha},n)}}{\binom{n}{\dim V}}\right]^{1/\nu(p)}-W}=o(1),\]
    where the $o(1)$ term tends to zero as $n$ goes to infinity for any fixed joints configuration $(\cJ,\cV_1,\ldots, \cV_r)$.
    From now on, we will fix the choice of $\vec{\alpha}$ for each $n$ and drop them from the notation of $\cD_{p,V}$ accordingly.

    For each $p\in\cJ$, choose $S_1\cup\cdots\cup S_r\in \cM(p)$ with the smallest product $\prod_{V\in S_1\cup\cdots \cup S_r}\abs{\cD_{p,V}}$, and set $V_i(p) = S_i$ for each $i\in[r]$.
    Then by \cite[Lemma 5.9]{TYZ22}, we have 
    \[\sum_{p\in\cJ}\omega(p)W\geq \frac{\binom{n+d}{d}}{\prod_{i=1}^{r}\binom{n}{k_i}^{m_i}}-o(1)\]
    and so
    \[W\geq \frac{\prod_{i=1}^{r}(k_i!)^{m_i}}{J_{\omega}\cdot d!}-o(1)\]
    where $J_{\omega}=\sum_{p\in\cJ}\omega(p)$.
    We also have that for all $p\in\cJ$,
    \begin{align*}
        \omega(p)W^{1/s}&=\nu(p)\left[\prod_{S_1\cup\cdots\cup S_r\in \cM(p)}\prod_{i=1}^{r}m_i^{-m_i}(\deg\cV_i)^{m_i}\prod_{V\in S_i}\frac{m_i\abs{\cD_{p,V}}}{\binom{n}{\dim V}\deg \cV_i}\right]^{1/s\nu(p)}+o(1)\\
        &\leq \frac{\prod_{i=1}^{r}(\deg \cV_i)^{m_i/s}}{s\prod_{i=1}^{r}m_i^{m_i/s}}\sum_{S_1\cup\cdots\cup S_r\in \cM(p)}\sum_{i=1}^{r}\sum_{V\in S_i}\frac{m_i\abs{\cD_{p,V}}}{\binom{n}{\dim V}\deg \cV_i}+o(1)
    \end{align*}
    by AM-GM.
    By summing over all $p\in\cJ$ and rearranging, we get
    \begin{align*}
        sJ_{\omega}W^{1/s}\prod_{i=1}^{r}m_i^{m_i/s}(\deg \cV_i)^{-m_i/s}\leq&\sum_{p\in\cJ}\sum_{S_1\cup\cdots\cup S_r\in \cM(p)}\sum_{i=1}^{r}\sum_{V\in S_i}\frac{m_i\abs{\cD_{p,V}}}{\binom{n}{\dim V}\deg \cV_i}+o(1)\\
        \leq&\sum_{i=1}^{r}\sum_{V\in \cV_i}\sum_{p\in \cJ\cap V}\frac{m_i\abs{\cD_{p,V}}}{\binom{n}{\dim V}\deg \cV_i}+o(1)\\
        =&s+o(1),
    \end{align*}
    where the last equation comes from \cite[Equation (5.1)]{TYZ22} and the fact that $\abs{\cD_{p,V}}=\abs{\cB_{p,V}}$.
    By rearranging, we get
    \[W\leq \frac{\prod_{i=1}^{s}(\deg \cV_i)^{m_i}}{J_\omega^s\prod_{i=1}^{r}m_i^{m_i}}+o(1).\]
    Now by comparing the upper and lower bound of $W$ and taking $n$ to infinity, we get
    \[\frac{\prod_{i=1}^{s}(\deg \cV_i)^{m_i}}{J_\omega^s\prod_{i=1}^{r}m_i^{m_i}}\geq  \frac{\prod_{i=1}^{r}(k_i!)^{m_i}}{J_{\omega}\cdot d!}.\]
    Rearranging this inequality gives
    \[J_{\omega}\leq \left(\frac{d!}{\prod_{i=1}^{r}k_i!^{m_i}m_i^{m_i}}\prod_{i=1}^{r}(\deg \cV_i)^{m_i}\right)^{1/(s-1)},\]
    as desired.
    This proves the statement for connected joints configuration.
    For general joint configuration, it suffices to prove the inequality
    \[\sum_{i=1}^{N}\prod_{j=1}^{r}a_{i,j}^{m_j/(s-1)}\leq \prod_{j=1}^{r}\left(\sum_{i=1}^{N}a_{i,j}\right)^{m_j/(s-1)}\]
    where $a_{i,j}>0$ for any $i\in[N], j\in[r]$ by splitting the joints configuration into connected components.
    By using the reverse Minkowski inequality for $p=0$ and the power mean inequality, we get
    \[\prod_{j=1}^{r}\left(\sum_{i=1}^{N}a_{i,j}\right)^{m_j/(s-1)}\geq \left(\sum_{i=1}^{N}\prod_{j=1}^{r}a_{i,j}^{m_j/s}\right)^{s/(s-1)}\geq \sum_{i=1}^{N}\prod_{j=1}^{r}a_{i,j}^{m_j/(s-1)},\]
    as desired.
\end{proof}

\begin{remark}
    This shows that key ideas and arguments in \cite{TYZ22} already lead to an optimal upper bound once we count the multiplicity correctly.
    At the same time, this shows that if one believes that the upper bound could be improved by a constant when we do not count the joints with multiplicities, a completely new argument would be needed in order to really take advantage of the fact that we are really only counting each joint once.
    In \cite{CY23}, we have successfully utilized this fact for a very special case.
    We will discuss this more in \cref{subsection:MultijointGraphConst}.
\end{remark}

\subsection{Second definition of multiplicity}
By inspecting the proof of \cref{theorem:FirstMult}, in fact, we can come up with another notion of multiplicity that leads to a stronger result.
Let $\omega(p) =\nobreak \nu^*(p)^{s/(s-1)}$ where $\nu^*(p)$ is the new multiplicity that will be defined later.
For simplicity, for every $p\in\cJ$, every edge $S=S_1\cup\cdots \cup S_r\in E(G_p)$ every handicap $\vec{\alpha}$ and positive integer $n$, we write $w_p(S,\vec{\alpha},n)=\nobreak\prod_{V\in S_1\cup\cdots \cup S_r}\abs{\cD_{p,V}(\vec{\alpha},n)}/\binom{n}{\dim V}$.
Moreover, if $S=S_1\cup\cdots\cup S_r$ has the smallest $w_p(S,\vec{\alpha},n)$ for a fixed $p$, we set $V_i(p) = S_i$ for each $i\in[r]$.
Previously, for each $p$ we fixed a packing $\cM(p)$ in $G_p$ of size $\nu(p)$ and set $\vec{\alpha},W$ so that 
\[\abs{\frac{1}{\omega(p)}\prod_{S\in \cM(p)}w_p(S,\vec{\alpha},n)^{1/\nu(p)}-W}=o(1)\]
holds for each $p\in\cJ$.
This way, it is guaranteed that $W\geq \omega(p)^{-1}w_p(V_1(p)\cup\cdots \cup V_r(p),\vec{\alpha},n)$ and we can carry on from here to get a lower bound of $W$ in terms of $J_{\omega}$.
It is thus natural to generalize this by replacing the geometric mean over a maximum packing $\cM(p)$ with some weighted geometric mean over the entire edge set $E(G_p)$.
To this end, for every $p\in\cJ$, let $\mu_p$ be a probability mass function on $E(G_p)$.
We would like to find a handicap $\vec{\alpha}$ for each $n$ so that
\[\abs{\frac{1}{\omega(p)}\prod_{S\in E(G_p)}w_p(S,\vec{\alpha},n)^{\mu_p(S)}-W}=o(1)\]
holds for all $p\in\cJ$.
Assuming it is possible to choose such a handicap, the argument in \cref{theorem:FirstMult} still gives the same lower bound $W\geq \prod_{i=1}^{r}(k_i!)^{m_i}/(J_{\omega}\cdot d!)-o(1)$.
To make the whole proof still work, we need to choose $\omega(p)=\nu^*(p)^{s/(s-1)}$ cleverly now.

As in the proof, we would like to upper bound $\omega(p)W^{1/s}$ and then sum this over all $p\in\cJ$.
We have
\[\omega(p)W^{1/s}=\nu^*(p)\left[\prod_{i=1}^{r}m_i^{-m_i}(\deg\cV_i)^{m_i}\prod_{S=S_1\cup\cdots\cup S_r\in \cM(p)}\left(\prod_{V\in S_i}\frac{m_i\abs{\cD_{p,V}}}{\binom{n}{\dim V}\deg \cV_i}\right)^{\mu_p(S)}\right]^{1/s}+o(1).\]
For any $V\in \cV_1\cup\cdots\cup\cV_r$, with a slight abuse of notation, set $\mu_p(V)\eqdef \sum_{S\in E(G_p), S\ni V}\mu_p(S)$.
Then
\[\omega(p)W^{1/s}=\nu^*(p)\left[\prod_{i=1}^{r}m_i^{-m_i}(\deg\cV_i)^{m_i}\prod_{V\in \cV_i}\left(\frac{m_i\abs{\cD_{p,V}}}{\binom{n}{\dim V}\deg \cV_i}\right)^{\mu_p(V)}\right]^{1/s}+o(1)\]
To make the argument work, it suffices to upper bound this by
\[\frac{\prod_{i=1}^{r}(\deg \cV_i)^{m_i/s}}{s\prod_{i=1}^{r}m_i^{m_i/s}}\sum_{i=1}^{r}\sum_{\substack{V\in \cV_i,\\ \mu_p(V)>0}}\frac{m_i\abs{\cD_{p,V}}}{\binom{n}{\dim V}\deg \cV_i}+o(1)\]
as $\mu_p(V)>0$ implies $p\in V$, and so it suffices to choose $\nu^*(p)$ so that
\[\nu^*(p)\prod_{i=1}^{r}\prod_{V\in \cV_i}\left(\frac{m_i\abs{\cD_{p,V}}}{\binom{n}{\dim V}\deg \cV_i}\right)^{\mu_p(V)/s}\leq \frac{1}{s}\sum_{i=1}^{r}\sum_{\substack{V\in \cV_i\\ \mu_p(V)>0}}\frac{m_i\abs{\cD_{p,V}}}{\binom{n}{\dim V}\deg \cV_i}.\]
Note that $\sum_{i=1}^{r}\sum_{V\in \cV_i}\mu_p(V)/s = 1$ as each $S\in E(G_p)$ contributes $s\mu_p(S)$ to the sum.
Therefore, by weighted AM-GM, we may bound
\begin{align*}
    &\nu^*(p)\prod_{i=1}^{r}\prod_{V\in \cV_i}\left(\frac{m_i\abs{\cD_{p,V}}}{\binom{n}{\dim V}\deg \cV_i}\right)^{\mu_p(V)/s}\\
    =&\nu^*(p)\prod_{i=1}^{r}\prod_{\substack{V\in \cV_i\\ \mu_p(V)>0}}\left(\frac{\mu_p(V)}{s}\right)^{\mu_p(V)/s}\left(\frac{s}{\mu_p(V)}\cdot \frac{m_i\abs{\cD_{p,V}}}{\binom{n}{\dim V}\deg \cV_i}\right)^{\mu_p(V)/s}\\
    \leq&\left(\frac{\nu^*(p)}{s}\prod_{i=1}^{r}\prod_{\substack{V\in\cV_i\\ \mu_p(V)>0}}\mu_p(V)^{\mu_p(V)/s}\right)\sum_{i=1}^{r}\sum_{\substack{V\in \cV_i\\ \mu_p(V)>0}}\frac{m_i\abs{\cD_{p,V}}}{\binom{n}{\dim V}\deg \cV_i}.
\end{align*}
This motivates the following definition.

\begin{definition}
    For each $p\in\cJ$, recall that $G_p$ is an $(m_1+\cdots+m_r)$-uniform hypergraph on the vertex (multi)set $\cV_1\cup\cdots\cup\cV_r$, where we add an edge $S = S_1\cup\cdots\cup S_r$ for any $(S_1,\ldots, S_r)\in \binom{\cV_1}{m_1}\times\cdots\times\binom{\cV_r}{m_r}$ such that $S_1\cup\cdots\cup S_r$ form a joint at $p$.
    Define
    \[\nu^*(p)\eqdef \max_{\mu_p}\prod_{i=1}^{r}\prod_{V\in \cV_i}\mu_p(V)^{-\mu_p(V)/s},\]
    where the maximum is taken over all probability mass function $\mu_p$ on $E(G_p)$, the notation $\mu_p(V)$ means $\sum_{S\in E(G_p),\, S\ni V}\mu_p(S)$ and $\mu_p(V)^{-\mu_p(V)/s}$ is set to $1$ if $\mu_p(V)=0$.
\end{definition}

As this definition is artificial, we first connect it to the more natural notions of multiplicities.
The proof would also show another way of thinking about $\nu^*(p)$.
Recall from \cite{TYZ22} that $M(p)$ is the number of ways varieties from $\cV_1\cup\cdots\cup \cV_r$ form a joint at $p$.
Equivalently, in our notation, $M(p) = \abs{E(G_p)}$.

\begin{proposition}\label{proposition:mult}
    For any $(k_1,\ldots, k_r;m_1,\ldots, m_r)$-joints configuration $(\cJ,\cV_1,\ldots, \cV_r)$ and any $p\in\cJ$, we have
    \[\nu^*(p)\geq \nu(p)\]
    and
    \[\nu^*(p)\geq \left(\prod_{i=1}^{r}\frac{m_i!}{m_i^{m_i}}\cdot M(p)\right)^{1/s}.\]
\end{proposition}
\begin{proof}
    To show the first inequality, it suffices to take a packing $\cM(p)\subseteq E(G_p)$ of size $\nu(p)$ and take $\mu_p$ to be the uniform measure on $\cM(p)$.
    Then for every $V\in \cV_1\cup\cdots\cup \cV_r$ with $\mu_p(V)>0$, we have $\mu_p(V) = 1/\nu(v)$.
    Moreover, there are $s\nu(v)$ of such $V$'s.
    Therefore
    \[\nu^*(p)\geq \left(\frac{1}{\nu(p)}\right)^{-\frac{1}{s\nu(p)} \cdot s\nu(p)}=\nu(p),\]
    as desired.

    To show the second inequality, we will need another interpretation of the definition of $\nu^*$.
    Given any probability measure $\mu_p$ on $E(G_p)$, let $V_{\mu_p}$ be a random variable on $\cV_1\cup\cdots \cup\cV_r$ chosen by the following process: we first sample $S\in E(G_p)$ according to $\mu_p$, and then we pick $V_{\mu_p}$ uniformly from $S=S_1\cup\cdots\cup S_r$.
    Then
    \[\log_2\prod_{i=1}^{r}\prod_{V\in\cV_i}\mu_p(V)^{-\mu_p(V)/s} = H(V_{\mu_p})-\log_2s.\]
    This entropic interpretation is going to be useful in the proof.

    Now let $\mu_p$ be the probability mass function corresponding to the uniform distribution on $E(G_p)$.
    Choose $S$ uniformly from $E(G_p)$, and permute the elements in $S_1\cup \cdots \cup S_r$ randomly so that $S_1\cup\cdots\cup S_r = \{V_1,\ldots, V_{m_1}\}\cup\{V_{m_1+1},\ldots, V_{m_1+m_2}\}\cup\cdots\cup\{V_{s-m_r+1},\ldots, V_s\}.$
    Then
    \[H(V_1,\ldots, V_s) = \log_2M(p)+\sum_{i=1}^{r}\log_2m_i!.\]
    On the other hand, by subadditivity, we have
    \[H(V_1,\ldots, V_s) \leq \sum_{j=1}^{s}H(V_j).\]
    We will relate the right hand side to $sH(V_{\mu_p})$. 
    In order to do so, set $i(V_{\mu_p})$ to be the unique $i$ such that $V_{\mu_p}$ is chosen from $S_i\subseteq \cV_i$.
    Then
    \[H(V_{\mu_p}) = H(V_{\mu_p},i(V_{\mu_p})) = H(i(V_{\mu_p})) + H(V_{\mu_p}\mid i(V_{\mu_p})) = -\sum_{i=1}^{r}\frac{m_i}{s}\left(\log_2m_i-\log_2s\right)+\frac{1}{s}\sum_{j=1}^{s}H(V_j).\]
    As a consequence,
    \[\log_2M(p)+\sum_{i=1}^{r}\log_2m_i!\leq sH(V_{\mu_p})-s\log_2 s+\sum_{i=1}^{r}m_i\log_2m_i,\]
    and so
    \[s\log_2\nu^*(p)\geq sH(V_{\mu_p})-s\log_2s\geq \log_2M(p)+\sum_{i=1}^{r}\log_2(m_i!/m_i^{m_i}).\]
    We get the desired inequality after rearranging.
\end{proof}

\begin{theorem}\label{theorem:SecondMult}
    In the setting of \cref{theorem:FirstMult}, we have
    \[\sum_{p\in\cJ}\nu^*(p)^{(m_1+\cdots+m_r)/(m_1+\cdots+m_r-1)}\leq C^{\nu^*}_{k_1,\ldots, k_r;m_1,\ldots, m_r}\left(\left(\deg\cV_1\right)^{m_1}\cdots \left(\deg\cV_r\right)^{m_r}\right)^{1/(m_1+\cdots+m_r-1)},\]
    where
    \[C^{\nu^*}_{k_1,\ldots, k_r;m_1,\ldots, m_r}=C^{\nu}_{k_1,\ldots, k_r;m_1,\ldots, m_r}=\left(\frac{d!}{\prod_{i=1}^{r}k_i!^{m_i}m_i^{m_i}}\right)^{1/(m_1+\cdots+m_r-1)}.\]    
    Moreover, this constant cannot be improved.
\end{theorem}
\begin{proof}
    We first note that from \cref{theorem:FirstMult} and the fact that $\nu^*(p) \geq \nu(p)$ by \cref{proposition:mult}, we know that the constant cannot be improved.
    It remains to show the inequality.

    The proof of the inequality is mostly similar to the proof of the one in \cref{theorem:FirstMult}.
    We will point out several modifications needed to prove \cref{theorem:SecondMult}.
    First of all, we will set $\omega(p) = \nu^*(p)^{s/(s-1)}$.
    Moreover, the way we choose the handicap would be slightly different.
    For every $p\in\cJ$, let $\mu_p$ be a probability mass function on $E(G_p)$ that maximizes $\prod_{i=1}^{r}\prod_{V\in\cV_i}\mu_p(V)^{-\mu_p(V)/s}$ where $s=m_1+\cdots+m_r$.
    As in the motivation of the definition of $\nu^*(p)$, we would like to choose $\vec{\alpha}$ for each $n\in\NN$ such that
    \[\abs{\frac{1}{\omega(p)}\prod_{S=S_1\cup\cdots\cup S_r\in E(G_p)}\left[\prod_{V\in S_1\cup\cdots\cup S_r}\frac{\abs{\cD_{p,V}(\vec{\alpha},n)}}{\binom{n}{\dim V}}\right]^{\mu_p(S)}-W}=o(1).\]
    Though we cannot directly apply \cite[Lemma 5.10]{TYZ22} anymore, the proof there still applies verbatim once we replace $W_p(\vec{\alpha})$ with the first term in the absolute value.
    Therefore we may assume that for each $n\in\NN$ we have fixed a choice of handicap $\vec{\alpha}$ with the desired property.

    The rest of the proof is basically identical to the one of \cref{theorem:SecondMult}.
    We just need to replace $\cM(p)$ by $E(G_p)$ and also replace the AM-GM step when we upper bound $W$ with the new weighted AM-GM that we have done already.
\end{proof}

\begin{remark}
    Now \cite[Theorem 1.10 (b)]{TYZ22} follows as an immediate corollary of \cref{theorem:SecondMult} because of the inequality between $\nu^*(p)$ and $M(p)$ in \cref{proposition:mult}.
    As in the remark after \cref{theorem:FirstMult}, this means that the method in the proof of \cite[Theorem 1.10 (b)]{TYZ22} is already optimal in some sense, and improving the constant there would most likely require some new ideas.
    We make a brief remark here that the joints configurations that show the optimality of $C^{\nu^*}_{k_1,\ldots,k_r;m_1,\ldots, m_r}$ does not show the optimality of \cite[Theorem 1.10 (b)]{TYZ22} as \cref{proposition:mult} is not tight in that example.
    However, \cref{proposition:mult} can be tight, for example, when the tangent flats of the varieties at $p$ are in general position.
\end{remark}

To show how versatile \cref{theorem:SecondMult} is, we consider the following weighted version of the joints of varieties theorem.
Suppose that we do not count joints with multiplicities, i.e. for the joints configuration $(\cJ,\cV_1,\ldots, \cV_r)$ we pick for each $p\in\cJ$ subsets $V_1(p),\ldots, V_r(p)$ of $\cV_1,\ldots, \cV_r$, respectively, so that they form a joint at $p$.
We will consider nonnegative weights $w$ on $\cV_1\cup\cdots\cup \cV_r$, and for each $i\in [r]$ we set $w(\cV_i)\eqdef \sum_{V\in \cV_i}w(V)\deg V$.

\begin{corollary}\label{corollary:WeightedJoints}
    Let $k_1,\ldots,k_r,m_1,\ldots, m_r$ be some positive integers.
    For each $i=1,\ldots, r$, let $\cV_i$ be a finite set of $k_i$-dimensional varieties in $\FF^{m_1k_1+\cdots+m_rk_r}$.
    Let $\cJ$ be the set of joints where each joint $p$ is formed by the varieties $V_1(p)\cup\cdots\cup V_r(p)$, where $V_i(p)$ is a subset of $\cV_i$ of size $m_i$ for each $i=1,\ldots, r$.
    Let $w$ be a nonnegative weight on $\cV_1\cup\cdots\cup \cV_r$.
    Then
    \[\sum_{p\in\cJ}\prod_{V\in V_1(p)\cup\cdots\cup V_r(p)}w(V)^{1/(m_1+\cdots+m_r-1)}\leq C^{\nu^*}_{k_1,\ldots, k_r;m_1,\ldots, m_r}\left(w\left(\cV_1\right)^{m_1}\cdots w\left(\cV_r\right)^{m_r}\right)^{1/(m_1+\cdots+m_r-1)}.\]
\end{corollary}
\begin{proof}
    Note that both sides are continuous with respect to $w$, and the validity of the inequality remains unchanged when we scale $w$.
    Therefore it suffices to show the inequality when $w(V)\in\NN$ for any $V\in\cV_1,\ldots,\cV_r$.

    For each $i\in[r]$, let $\Tilde{\cV}_i$ be the multiset obtained by repeating each $V\in \cV_i$ for $w(V)$ times.
    For each $p\in\cJ$, let $\Tilde{H}_p$ be the associated $(m_1+\cdots+m_r)$-uniform hypergraph of $p$ in the joints configuration $(\cJ,\Tilde{\cV}_1,\ldots, \Tilde{\cV}_r)$.
    Let $\mu_p$ be the probability mass function on $\Tilde{H}_p$ obtained by choosing a random copy of $V$ for each $V\in V_1(p)\cup\cdots\cup V_r(p)$.
    Then it is easy to verify that
    \[\nu^*(p)\geq \prod_{V\in V_1(P)\cup\cdots\cup V_r(p)}w(V)^{1/(m_1+\cdots+m_r)}.\]
    Therefore the statement follows once we plug this into \cref{theorem:SecondMult}.
\end{proof}

In some sense, the weight $w$ allows us to focus more on a part of the joints configuration and less on the other parts of configuration.
It is natural to focus more on the denser part of the joints configuration, and this can be done by, for example, taking $w(V) = \#\{p\in V\}$ (where, with a slight abuse of notation, $p\in V$ if and only if $V\in V_1(p)\cup\cdots\cup V_r(p)$).
This leads to the following result.

\begin{corollary}\label{corollary:PointCountMult}
    Let $k_1,\ldots,k_r,m_1,\ldots, m_r$ be some positive integers.
    For each $i=1,\ldots, r$, let $\cV_i$ be a finite set of $k_i$-dimensional varieties in $\FF^{m_1k_1+\cdots+m_rk_r}$.
    Let $\cJ$ be the set of joints where each joint $p$ is formed by the varieties $V_1(p)\cup\cdots\cup V_r(p)$, where $V_i(p)$ is a subset of $\cV_i$ of size $m_i$ for each $i=1,\ldots, r$.
    We write $p\in V$ if and only if $V\in\cV_i$ and $V\in V_i(p)$ for some $i\in[r]$.
    For each $p\in \cJ$, let $\omega(p)=\prod_{V\in V_1(p)\cup\cdots\cup V_r(p)}\#\{p\in V\}.$
    Then
    \[\frac{1}{J}\sum_{p\in \cJ}\omega(p)^{1/(m_1+\cdots+m_r-1)}\leq \left(\frac{d!\cdot J}{k_1!^{m_1}\cdots k_r!^{m_r}}\right)^{1/(m_1+\cdots+m_r-1)}.\]
\end{corollary}
The weight $\omega(p)$ can be thought of as the importance of $p\in \cJ$, and the statement basically says that the $1/(m_1+\cdots+m_r-1)$-th power mean of the importances is at most $J$ times some constant.
Intuively, if the joints configuration is ``too dense'' at a joint, then it must be ``sparse'' somewhere else.

\begin{proof}
    In \cref{corollary:WeightedJoints}, set $w(V) = \#\{p\in V\}$ for every $V\in \cV_1\cup\cdots\cup\cV_r$.
    Then the left hand side becomes $\sum_{p\in \cJ}\omega(p)^{1/(m_1+\cdots+m_r-1)}.$
    To compute the right hand side, note that for each $i\in[r]$,
    \[w(\cV_i) = \sum_{V\in \cV_i}\#\{p\in V\} = \sum_{p\in \cJ}\abs{V_i(p)} = m_iJ.\] 
    Therefore the right hand side of the inequality is
    \begin{align*}C^{\nu^*}_{k_1,\ldots, k_r;m_1,\ldots, m_r}\left(w\left(\cV_1\right)^{m_1}\cdots w\left(\cV_r\right)^{m_r}\right)^{1/(m_1+\cdots+m_r-1)}=&\left(\frac{d!J^{m_1+\cdots+m_r}\prod_{i=1}^{r}m_i^{m_i}}{\prod_{i=1}^{r}k_i!^{m_i}m_i^{m_i}}\right)^{1/(m_1+\cdots+m_r-1)}\\
    =&\left(\frac{d!J^{m_1+\cdots+m_r}}{\prod_{i=1}^{r}k_i!^{m_i}}\right)^{1/(m_1+\cdots+m_r-1)}.
    \end{align*}
    Then after rearranging the inequality obtained from \cref{corollary:WeightedJoints}, we get the desired result.
\end{proof}

\begin{remark}
    As discussed earlier, the inequality implies that the $(1/(m_1+\cdots+m_r-1))$-th power mean of the $\omega(p)$'s is upper bounded by $J$ times some constant.
    It is natural to ask what is the largest $q$ such that this remains true when $1/(m_1+\cdots+m_r-1)$ is replaced by $q$, and indeed we can see that $1/(m_1+\cdots+m_r-1)$ is already optimal in this sense.
    To this end, let $N$ be a large positive integer, and consider the joints configuration where $0$ is formed by $V_1(0)\cup\cdots\cup V_r(0)$, and we artificially put $N$ more joints on each of the variety $V$ in $V_1(0)\cup\cdots\cup V_r(0)$ by adding more varieties around each of them.
    Then we have $J=(m_1+\cdots+m_r)N+1$ and $\omega(0) = (N+1)^{m_1+\cdots+m_r}$.
    Therefore, if the $q$-th power mean of the $\omega(p)$'s are bounded by $O(J)$, then
    \[\left(\frac{1}{J}(N+1)^{(m_1+\cdots+m_r)q}\right)^{1/q}=O(N),\]
    and so
    \[(m_1+\cdots+m_r)-\frac{1}{q}\leq 1\]
    by comparing the exponents of $N$ on both sides.
    This gives $q\leq 1/(m_1+\cdots+m_r-1)$, as desired.
\end{remark}

\section{Set-theoretic Version of the Joints Problem}\label{section:Graph}
In this section, we will investigate the set-theoretic versions of the joints problems.
This will allow us to connect the Kruskal--Katona theorem and the partial shadow problem to the joints problem.
With this connection, we will prove \cref{theorem:PartialShadow}.
We will also disprove some conjectures proposed in \cite{TYZ22} and \cite{YZ23}.
After that, we will show that $C_{k,k,k;1,1,1}$ grows exponentially in $k$, and also relate the joints problem to a closely related paper of ours \cite{CY23}.

If we consider $n$ hyperplanes $H_1,\dots,H_n$ in $\FF^d$ in general position, the intersection $\cap_{i\in F}H_i$ is a $k$-flat for any $(d-k)$-subset $F\subseteq [n]$. Similarly, each $d$-subset $P$ gives a point in $\FF^d$. Note that, for any $d$-subset $P\subseteq [n]$ and any subsets $F_1,\dots,F_m\subseteq [n]$, if $P\setminus F_1,\ldots,P\setminus F_m$ partition $P$, then the flats corresponding to $F_1,\dots,F_m$ form a joint at the point corresponding to $P$. In this way, we can state the set-theoretic version of the joints of flats problem.

\begin{definition}\label{definition:JointSet}
    Given positive integers $k_1,\dots,k_r,m_1,\dots,m_r$, let $d=k_1m_1+\dots+k_rm_r$ and $s=m_1+\dots+m_r$. Given $\bF_1,\dots,\bF_r$, where $\bF_i$ is a family of $(d-k_i)$-subsets of $\NN$, we say that $P\subseteq \NN$ is a \emph{$(k_1,\dots,k_r;m_1,\dots,m_r;0)$-joint set} if there exist sets $F_1,\dots,F_s$, where $m_i$ of them are from $\bF_i$, such that $P\setminus F_1,\dots ,P\setminus F_s$ partition $P$. 

    We say that $(\bJ;\bF_1,\dots,\bF_r)$ is a \emph{$(k_1,\dots,k_r;m_1,\dots,m_r;0)$-joint set system} if every $P\in\bJ$ is a $(k_1,\dots,k_r;m_1,\dots,m_r;0)$-joint set of $\bF_1,\dots,\bF_r$.
\end{definition}
For example, $(\bJ,\bF)$ is a $(1;d;0)$-joint set system if and only if $\bF$ contains $\partial \bJ$, where $\partial \bJ$ is the set of all the $(d-1)$-subsets of $P\in\bJ$. 
The Kruskal--Katona theorem gives a lower bound on the size of $\bF$ given the size of $\bJ$.
Note that every $(k_1,\dots,k_r;m_1,\dots,m_r;0)$-joint set system gives a $(k_1,\dots,k_r;m_1,\dots,m_r;0)$-joints configuration that is generically induced. 

As discussed in the introduction and \cref{example:BE}, sometimes it is better to consider projected generically induced configurations.
We may modify \cref{definition:JointSet} to accommodate for those configurations.

\begin{definition}
    Let $k_1,\ldots, k_r,m_1,\ldots, m_r,d,s$ be as mentioned in \cref{definition:JointSet}.
    Let $\delta$ be a nonnegative integer.
    Given $\bF_1,\dots,\bF_r$, where $\bF_i$ is a family of $(d+\delta-k_i)$-subsets of $\NN$, we say that $P\subseteq \NN$ is a \emph{$(k_1,\dots,k_r;m_1,\dots,m_r;\delta)$-joint set} if $\abs{P}=d+\delta$ and there exist sets $F_1,\dots,F_s$, where $m_i$ of them are from $\bF_i$, such that $P\setminus F_1,\dots ,P\setminus F_s$ are disjoint subsets of $P$. 

    Similar to \cref{definition:JointSet}, we say that $(\bJ;\bF_1,\dots,\bF_r)$ is a \emph{$(k_1,\dots,k_r;m_1,\dots,m_r;\delta)$-joint set system} if every $P\in\bJ$ is a $(k_1,\dots,k_r;m_1,\dots,m_r;\delta)$-joint set of $\bF_1,\dots,\bF_r$.
\end{definition}
We briefly remark that in the definition above, if we take $\delta=0$, we get exactly \cref{definition:JointSet}.

Note that every $(k_1,\ldots, k_r;m_1,\ldots, m_r;\delta)$-joint set system induces a generically induced configuration in $\FF^{d+\delta}$ that might not be a joints configuration.
However, when we take a general projection from $\FF^{d+\delta}$ to $\FF^d$, all $(k_1,\ldots, k_r;m_1,\ldots, m_r;\delta)$-joint sets now correspond to an actual joint.
Therefore, every $(k_1,\ldots, k_r;m_1,\ldots, m_r;\delta)$-joint set system induces a projected generically induced joints configuration in $\FF^d$, given that $\abs{\FF}$ is sufficiently large.

Thus, if we restrict ourselves to generically induced configurations or projected generically induced configurations, we are actually asking for upper bounds of $|\bJ|$ in terms of $|\bF_1|,\dots,|\bF_r|$, where $(\bJ;\bF_1,\ldots, \bF_r)$ is a $(k_1,\ldots, k_r;m_1,\ldots, m_r;\delta)$-joint set system for $\delta=0$ or arbitrary $\delta$, respectively. 
From \cite{TYZ22}, we know that \[|\bJ|=O((|\bF_1|^{m_1}\dots|\bF_r|^{m_r})^{1/(s-1)}).\] 
Hence, we may define $C^G_{k_1,\dots,k_r;m_1,\dots,m_r}$ and $C^P_{k_1,\dots,k_r;m_1,\dots,m_r}$ to be the smallest constant $C$ such that
\[|\bJ|\leq C(|\bF_1|^{m_1}\dots|\bF_r|^{m_r})^{1/(s-1)}\]
when $\delta=0$ and when $\delta$ can be an arbitrary nonnegative integer, respectively.
It follows that 
\[C^G_{k_1,\dots,k_r;m_1,\dots,m_r}\leq C^P_{k_1,\dots,k_r;m_1,\dots,m_r}\leq C_{k_1,\dots,k_r;m_1,\dots,m_r}.\]

\subsection{Connection to the partial shadow problem}
In this subsection, we will focus on $(1;d;\delta)$-joint set systems, which induce projected generically induced joints configurations for lines in $\FF^d$.
Note that in this case, a $(1;d;\delta)$-joint set is simply a set $P$ of size $d+\delta$ such that there are $d$ subsets of $P$ of size $d+\delta-1$ in $\bF$.
In the notation of the partial shadow problem defined in \cref{subsection:PartialShadow}, we have $r=d+\delta$, $k=\delta$ and $m=\abs{\bJ}$.
Therefore, $f(d+\delta,\abs{\bJ},\delta)$ is the minimum possible size of $\bF$.

Note that the configuration in \cref{example:BE} is the projected generically induced configuration induced by $\bJ=\binom{[6]}{5}$ and $\bF=\binom{[6]}{4}\backslash\{\{1,2,3,4\},\{1,2,5,6\},\{3,4,5,6\}\}$ when $d=4$ and $\delta=1$.
This is the example given by Bollob\'as and Eccles \cite{BE15}, showing that $f(5,6,1)\leq 12<13=f(4,6,0)$.
Here, $f(4,6,0)=13$ is the Kruskal--Katona bound.
This shows that, contrary to what Zhao and the second author suggested at the end of Section 3 of \cite{YZ23}, the optimal joints configuration might not be achieved by generically induced configurations.
However, it is still natural to ask the following.
\begin{question}
    Is the minimum number of lines in $\FF^d$ forming $J$ joints $\max_{\delta\geq 0}f(d+\delta,J,\delta)$?
    That is, can we always find a projected generically induced configuration that is optimal?
\end{question}
Fitch \cite[Chapter 3]{Fitch18} conjectured that when $J$ is sufficiently large with respect to $d$, the value $\max_{\delta\geq 0}f(d+\delta,J,\delta)$ is actually equal to the Kruskal--Katona bound $f(d,J,0)$.
Combined with this conjecture, we may ask the following.
\begin{question}
    Given a positive integer $d$, is it true that for all sufficiently large positive integers $J$, the number of lines in $\FF^d$ forming $J$ joints is at least $f(d,J,0)$?
\end{question}

Using this connection between the partial shadow problem and the joints problem, we can now prove \cref{theorem:PartialShadow}.

\begin{proof}[Proof of \cref{theorem:PartialShadow}]
    The case $r-k\leq 2$ is trivial, so we will deal with the case $r-k\geq 3$.
    Take $\FF$ to be an arbitrary infinite field, $\RR$ for example.
    Set $\delta = k$, $d=r-k$ and let $(\bJ,\bF)$ be a $(1;d;\delta)$-joint set system where $\abs{\bJ}=m$.
    By the discussion above, we know that $\bJ$ plays the role of $\cA$ and $\bF$ plays the role of $\cB$.
    Recall that $m=\binom{x}{r-k}=\binom{x}{d}$, and we need to show that $\abs{\bF}\geq \binom{x}{r-k-1}=\binom{x}{d-1}$.
    Since $\FF$ is infinite, the set system induces a projected generically induced joints configuration in $\FF^d$ with $\abs{\bJ}=\binom{x}{d}$ joints and $\abs{\bF}$ lines.
    By \cref{theorem:SharpJoints}, we get $\abs{\bF}\geq \binom{x}{d-1}$, which proves the first part.

    Now if $\abs{\bF}=\binom{x}{d-1}$, \cref{theorem:Structural} shows that the projected generically induced configuration is determined by $x$ hyperplanes $H'_1,\ldots,H'_x$, where the lines are the $(d-1)$-wise intersections and the joints are the $d$-wise intersections.
    Suppose that the projected generically induced configuration is obtained by projecting the generically induced configuration in $\FF^{d+\delta}$ determined by the hyperplanes $H_1, H_2,\ldots$ via the projection $\pi:\FF^{d+\delta}\to \FF^d$.
    For each $(1;d;\delta)$-joint set $P$, we can compute the hyperplanes containing the joints corresponding to $P$ as follows.
    Let $F_1,\ldots, F_d$ be the subsets of $P$ in $\bF$ of size $d+\delta-1$.
    Let $\textup{core}(P)$ be the core part $(P\backslash F_1)\cup\cdots\cup (P\backslash F_d)$ of $P$, and let $\textup{red}(P)\eqdef P\backslash \textup{core}(P)$ be the redundant part.
    Then the $(d-1)$-flats containing $(d-1)$ lines passing through $\bigcap_{i\in P}H_i$ are $\left(\bigcap_{i\in \textup{red}(P)}H_i\right)\cap H_j$ for all $j\in \textup{core}(P)$.
    When $\pi$ is sufficiently general, it maps two different $(d-1)$-flats determined above to different hyperplanes.
    We will fix such $\pi$.
    Note that for any two joints on the same line, they share all the $(d-1)\geq 2$ hyperplanes.
    In particular, for any $P,P'$ that correspond to two joints on the same line, there exists $j_1,j_2\in \textup{core}(P)$ and $j_1',j_2'\in \textup{core}(P')$ such that $\textup{red}(P)\cup \{j_i\}=\textup{red}(P')\cup\{j_i'\}$ for $i=1,2$.
    By taking the intersection, we get $\textup{red}(P)=\textup{red}(P')$.
    Since the joints configuration is connected, we see that the redundant parts of all sets in $\bJ$ are the same.
    Call this common redundant part $Y$.
    It is clear that $\abs{Y}=\delta$.
    Then each $H'_i$ is of the form $\pi\left(\left(\bigcap_{j\in Y}H_j\right)\cap H_{i'}\right)$, and we may relabel $H_1,H_2,\ldots$ so that $i'=i$.
    Take $X=\{1,\ldots,x\}$.
    Then for any $S\in \binom{X}{d}$, there is a unique set $P\in \bJ$ that corresponds to the joint $\cap_{i\in S}H_i'$.
    For each $i\in S$, the hyperplane $H_i'$ is one of the hyperplanes through the joint determined by $P$.
    By the characterization above, we have that $i\in \textup{core}(P)$.
    Therefore $S\subseteq \textup{core}(P)$, and as both have size $d$, we have $S=\textup{core}(P)$ and so $P=S\cup Y$.
    Note that there are $\binom{x}{d}$ such sets $S$, and different sets $S$ give rise to different elements of $\bJ$.
    Therefore, the family $\bJ$ contains precisely $S\cup Y$ for all $S\in \binom{X}{d}$.
    It is then clear that $\bF$ contains precisely $S\cup Y$ for all $S\in \binom{X}{d-1}$.
\end{proof}

\subsection{Disproving a conjecture in \cite{TYZ22}}
In \cite{TYZ22}, Tidor, Zhao and the second author conjectured that $N$ $2$-flats can form at most $(\sqrt{2}/3+o(1))N^{3/2}$ joints in $\FF^6$. 
We show in this section that this conjecture is false.

Indeed, we will show that $C^G_{2;3}\geq \sqrt{2/7}$ by constructing $(2;3)$-joint set systems $(\bJ,\bF)$ with $|\bJ|=28n^2,|\bF|=14n^2$.

We first construct the system for $n=1$.
Let $\bJ=\binom{[8]}{6}$ be all the $6$-subsets in $\{1,2,\dots,8\}$. Consider two copies of $K_4$ on vertex set $\{1,2,3,4\}$ and $\{5,6,7,8\}$ respectively. We color the edges in red, green, and blue in the following way. Let 
\begin{align*}
    R_1=\{\{1,2\},\{3,4\}\}, R_2=\{\{5,6\},\{7,8\}\},\\
    G_1=\{\{1,3\},\{2,4\}\}, G_2=\{\{5,7\},\{6,8\}\},\\
    B_1=\{\{1,4\},\{2,3\}\}, B_2=\{\{5,8\},\{6,7\}\}.
\end{align*}
Let $\bF$ be the family consists of $\{1,2,3,4\},\{5,6,7,8\}$, and all $4$-sets of the form $e_1\cup e_2$, where $e_1$ and $e_2$ are two edges that receive the same color, and $e_1, e_2$ are in the first and the second $K_4$ respectively. Note that we have $|\bJ|=\binom{8}{6}=28$ and $|\bF|=2+3\times 4=14$.

We shall check that $(\bJ,\bF)$ is a $(2;3)$-joint set system. Let $A\in\bJ$. If $A=\{a,b,5,6,7,8\}$, we can take $F_1=\{5,6,7,8\}, F_2=\{a,b,4+a,4+b\}, F_3=\{a,b,4+c,4+d\}$, where $\{a,b,c,d\}=\{1,2,3,4\}$. The case $A=\{1,2,3,4,4+a,4+b\}$ is similar. The remaining case is that $A\cap \{1,2,3,4\}$ and $A\cap \{5,6,7,8\}$ both have size $3$. In this case, both $A\cap \{1,2,3,4\}$ and $A\cap \{5,6,7,8\}$ are triangles that receive all three colors. Thus, we can pair the edges of the same color and get $F_1,F_2,F_3$.

We can then amplify this construction by copying each element in $\{1,2,\dots,8\}$ $n$ times.
For each element in $\bF$, we replace it with the $n^4$ copies, and each element in $\bJ$ is now replaced by its $n^6$ copies.
We end up with a $(2;3)$-joint set system with $|\bF|=14n^4$ and
\[|\bJ|=28n^6=\sqrt{\frac{2}{7}}|\bF|^{3/2}.\]
Therefore, $C_{2;3}\geq C^G_{2;3}\geq \sqrt{2/7}$.
The conjecture is also false by taking $n$ to infinity.
\begin{remark}
    In the last step of the construction, we see that we can always \emph{blow-up} a construction and get a larger construction with the same constant factor. That is, we may copy each element in our base set $n$ times. Each set $P\in \bJ$ becomes $n^{d}$ copies and each $F_i\in\bF_i$ becomes $n^{d-k_i}$ copies. The ratio $|\bJ|/(|\bF_1|^{m_1}\dots|\bF_r|^{m_r})^{1/(s-1)}$ in the blown-up construction remains the same.
\end{remark}

\subsection{Exponential growth of $C_{k,k,k;1,1,1}$}
From \cite[Theorem 1.10]{TYZ22}, we know that the constant $C_{k,k,k;1,1,1}$ for the multijoints problem of three $k$-flats is at most $((3k)!/(k!)^3)^{1/2}$, which is exponential in $k$. In this subsection, we will show that the constant is indeed exponential in $k$ by showing that $C^G_{k,k,k;1,1,1}$ is exponential in $k$.

We start with $k$ disjoint copies of $K_4$, where the vertex set of the $i$-th copy is $\{4i+1,\dots,4i+4\}$. For each $i\in [k]$, We color the edges in red, green, and blue by
\begin{align*}
    R_i=\{\{4i+1,4i+2\},\{4i+3,4i+4\}\}, \\
    G_i=\{\{4i+1,4i+3\},\{4i+2,4i+4\}\}, \\
    B_i=\{\{4i+1,4i+4\},\{4i+2,4i+3\}\}.
\end{align*}
Let $\bF_R$ be the collection of $2k$-sets of the form $\cup_{i\in [k]} e_i$, where $e_i\in R_i$, and define $\bF_G,\bF_B$ in the same way. Let 
\[\bJ=\{A\subseteq [4k]\mid |A\cap\{4i+1,\dots,4i+4\}|=3\text{ for all }i\in [k]\}.\]
We claim that $(\bJ;\bF_R,\bF_G,\bF_B)$ is a $(k,k,k;1,1,1)$-joint set system. Note that, for any $P\in\bJ$ and $i\in [k]$, the intersection $P\cap \{4i+1,\dots,4i+4\}$ is the vertex set of a triangle that receives all three colors. Thus, we can pair all the red, green, blue edges to get $F_R,F_G,F_B$, respectively. Note that $F_R,F_G,F_B$ certificate that $P$ is a joint set. We know that $|\bJ|=4^k$ and $|\bF_R|=|\bF_G|=|\bF_B|=2^k$. 

Again, we can then amplify this construction by copying each element in $[4k]$ $n$ times.
For each element in $\bF_R,\bF_G,\bF_B$, we replace it with the $n^{2k}$ copies, and each element in $\bJ$ is now replaced by its $n^{3k}$ copies.
We end up with a $(k,k,k;1,1,1)$-joint set system with $|\bF_R|=|\bF_G|=|\bF_B|=2^kn^{2k}$ and
\[|\bJ|=4^kn^{3k}=2^{k/2}|\bF_R|^{1/2}|\bF_G|^{1/2}|\bF_B|^{1/2}.\]
Therefore, $C_{k,k,k;1,1,1}\geq C^G_{k,k,k;1,1,1}\geq 2^{k/2}$.

For the $(k;3)$-joints problem and the $(k;3)$-joint set system problem, we can always combine the parts and get a construction. That is, if we set $\bF=\bF_R\cup\bF_G\cup\bF_B$, the pair $(\bJ;\bF)$ forms a $(k;3)$-joint set system. Since $|\bJ|=4^kn^{3k}$ and $|\bF|=3\cdot 2^{k}n^{2k}$, it follows that $C_{k;3}\geq C^G_{k;3}\geq 2^{k/2}/3\sqrt{3}$.
\subsection{Multijoints via multicolored hypergraphs}\label{subsection:MultijointGraphConst}
We are also interested in the set-theoretic version of the problem itself. 
Zhao and the second author conjectured in \cite{YZ23} that the best possible joints configuration can be given by a joint set system. 
We still believe so, and in terms of the notation in the paper, we make the following weaker conjecture.
\begin{conjecture}
For any positive integers $k_1,\dots,k_r,m_1,\dots,m_r$, we have
\[C^G_{k_1,\dots,k_r;m_1,\dots,m_r}=C_{k_1,\dots,k_r;m_1,\dots,m_r}.\]
\end{conjecture}
Therefore, the set-theoretic problems are also interesting to us, as they suggest what $C_{k_1,\ldots, k_r;m_1,\ldots, m_r}$ should be if the conjecture is true.

We have more discussion on the set-theoretic problems in \cite{CY23}. The main ingredient of the paper is entropy method. In \cite[Theorem 1.1]{CY23}, we show that $C^G_{1,1,1;1,1,1}=\sqrt{2}$. In comparison, we know that $C_{1,1,1;1,1,1}\leq \sqrt{6}$ from \cite{YZ23} and $C^{\nu}_{1,1,1;1,1,1}=\sqrt{6}$ from \cref{theorem:FirstMult}.
We also showed that 
\begin{align*}
   C^G_{\underbrace{1,\dots,1}_r;\underbrace{1,\dots,1}_r}\geq \begin{cases}
        (2^r/(r+1))^{1/(r-1)}&\textup{ if }r\textup{ is odd;}\\
        (2^{r-1}/r)^{1/(r-1)}&\textup{ if }r\textup{ is even.}
        \end{cases}
\end{align*}

For more discussion on some other cases, we also refer the readers to the paper.
\bibliographystyle{amsplain0}
\bibliography{ref_joints}

\providecommand{\bysame}{\leavevmode\hbox to3em{\hrulefill}\thinspace}
\providecommand{\MR}{\relax\ifhmode\unskip\space\fi MR }
\providecommand{\MRhref}[2]{%
  \href{http://www.ams.org/mathscinet-getitem?mr=#1}{#2}
}
\providecommand{\href}[2]{#2}
\begin{thebibliography}{10}

\bibitem{AS00}
Noga Alon and Joel~H. Spencer, \emph{The probabilistic method}, second ed.,
  Wiley-Interscience Series in Discrete Mathematics and Optimization,
  Wiley-Interscience [John Wiley \& Sons], New York, 2000, With an appendix on
  the life and work of Paul Erd\H{o}s.

\bibitem{BE15}
B\'{e}la Bollob\'{a}s and Tom Eccles, \emph{Partial shadows of set systems},
  Combin. Probab. Comput. \textbf{24} (2015), 825--828.

\bibitem{BC21}
Boris Bukh and Ting-Wei Chao, \emph{Sharp density bounds on the finite field
  {K}akeya problem}, Discrete Anal. (2021), Paper No. 26, 9pp.

\bibitem{CI14}
Anthony Carbery and Marina Iliopoulou, \emph{Counting joints in vector spaces
  over arbitrary fields}, arXiv:1403.6438.

\bibitem{CY23}
Ting-Wei Chao and Hung-Hsun~Hans Yu, \emph{Kruskal--{K}atona-type problems via
  entropy method}, arXiv:2307.15379.

\bibitem{CEGPSSS92}
Bernard Chazelle, Herbert Edelsbrunner, Leonidas~J. Guibas, Richard Pollack,
  Raimund Seidel, Micha Sharir, and Jack Snoeyink, \emph{Counting and cutting
  cycles of lines and rods in space}, Comput. Geom. \textbf{1} (1992),
  305--323.

\bibitem{Dvir09}
Zeev Dvir, \emph{On the size of {K}akeya sets in finite fields}, J. Amer. Math.
  Soc. \textbf{22} (2009), 1093--1097.

\bibitem{Dvir10}
Zeev Dvir, \emph{Incidence theorems and their applications}, Found. Trends
  Theor. Comput. Sci. \textbf{6} (2010), 257--393 (2012).

\bibitem{DKSS13}
Zeev Dvir, Swastik Kopparty, Shubhangi Saraf, and Madhu Sudan, \emph{Extensions
  to the method of multiplicities, with applications to {K}akeya sets and
  mergers}, SIAM J. Comput. \textbf{42} (2013), 2305--2328.

\bibitem{Fitch18}
Matthew Fitch, \emph{Applications of entropy to extremal problems}, 2018, Ph.D.
  Thesis, University of Warwick, https://wrap.warwick.ac.uk/132608/.

\bibitem{GT13}
Ben Green and Terence Tao, \emph{On sets defining few ordinary lines}, Discrete
  Comput. Geom. \textbf{50} (2013), 409--468.

\bibitem{Guth10}
Larry Guth, \emph{The endpoint case of the {B}ennett-{C}arbery-{T}ao
  multilinear {K}akeya conjecture}, Acta Math. \textbf{205} (2010), 263--286.

\bibitem{Guth-book}
Larry Guth, \emph{Polynomial methods in combinatorics}, University Lecture
  Series, vol.~64, American Mathematical Society, Providence, RI, 2016.

\bibitem{GK10}
Larry Guth and Nets~Hawk Katz, \emph{Algebraic methods in discrete analogs of
  the {K}akeya problem}, Adv. Math. \textbf{225} (2010), 2828--2839.

\bibitem{GK15}
Larry Guth and Nets~Hawk Katz, \emph{On the {Erd\H{o}s} distinct distances
  problem in the plane}, Ann. of Math. (2) \textbf{181} (2015), 155--190.

\bibitem{KSS10}
Haim Kaplan, Micha Sharir, and Eugenii Shustin, \emph{On lines and joints},
  Discrete Comput. Geom. \textbf{44} (2010), 838--843.

\bibitem{Lovasz93}
L\'{a}szl\'{o} Lov\'{a}sz, \emph{Combinatorial problems and exercises}, second
  ed., North-Holland Publishing Co., Amsterdam, 1993.

\bibitem{Nach-book}
Asaf Nachmias, \emph{Planar maps, random walks and circle packing}, Lecture
  Notes in Mathematics, vol. 2243, Springer, Cham, 2020, \'{E}cole
  d'\'{e}t\'{e} de probabilit\'{e}s de Saint-Flour XLVIII---2018.

\bibitem{Qui09}
Ren\'{e} Quilodr\'{a}n, \emph{The joints problem in {$\mathbb R^n$}}, SIAM J.
  Discrete Math. \textbf{23} (2009/10), 2211--2213.

\bibitem{Tao14}
Terence Tao, \emph{Algebraic combinatorial geometry: the polynomial method in
  arithmetic combinatorics, incidence combinatorics, and number theory}, EMS
  Surv. Math. Sci. \textbf{1} (2014), 1--46.

\bibitem{TYZ22}
Jonathan Tidor, Hung-Hsun Hans~Yu, and Yufei Zhao, \emph{Joints of varieties},
  Geom. Funct. Anal. \textbf{32} (2022), 302--339.

\bibitem{Wol99}
Thomas Wolff, \emph{Recent work connected with the {K}akeya problem}, Prospects
  in mathematics ({P}rinceton, {NJ}, 1996), Amer. Math. Soc., Providence, RI,
  1999, pp.~129--162.

\bibitem{YZ23}
Hung-Hsun~Hans Yu and Yufei Zhao, \emph{Joints tightened}, Amer. J. Math.
  \textbf{145} (2023), 569--583.

\end{thebibliography}
\end{document}